\documentclass{article}

\usepackage{fullpage}
\usepackage[english]{babel}
\usepackage[utf8]{inputenc}
\usepackage{amsmath,amsfonts,amssymb,amsthm}
\usepackage{enumitem}
\usepackage[all]{xy}
\usepackage{graphicx}
\usepackage{tikz}

\newcommand{\ZZ}{\mathbb{Z}}
\newcommand{\RR}{\mathbb{R}}

\newcommand{\eps}{\varepsilon}
\newcommand{\lbd}{\lambda}
\newcommand{\gam}{\gamma}
\newcommand{\dd}{\partial}
\newcommand{\SL}{\mathcal{SL}_n}
\newcommand{\dl}[1]{#1\!\!-\!\!\mathcal{SL}_n}
\newcommand{\dr}[1]{#1\!\!-\!\!r\mathcal{SL}_n}

\newcommand{\impc}{\overset{c}{\Rightarrow}}
\newcommand{\impw}{\overset{w}{\Rightarrow}}

\newcommand{\phidSL}{\Phi_{2\!-\!\SL}}

\newcommand\xqed[1]{%
  \leavevmode\unskip\penalty9999 \hbox{}\nobreak\hfill
  \quad\hbox{#1}}
\newcommand\demo{\xqed{$\square$}}

\theoremstyle{definition}
\newtheorem{de}{\normalsize Definition}[section]
\newtheorem{lem}[de]{\normalsize Lemma}
\newtheorem{prop}[de]{\normalsize Proposition}
\newtheorem{theo}[de]{\normalsize Theorem}
\newtheorem{cor}[de]{\normalsize Corollary}

\makeatletter
\newcommand{\subjclass}[2][1991]{%
	\let\@oldtitle\@title%
	\gdef\@title{\@oldtitle\footnotetext{#1 \emph{Mathematics subject classification.} #2.}}%
}
\newcommand{\keywords}[1]{%
	\let\@@oldtitle\@title%
	\gdef\@title{\@@oldtitle\footnotetext{\emph{Key words and phrases.} #1.}}%
}
\makeatother

\author{Boris Colombari}

\title{Welded extensions and ribbon restrictions of diagrammatical moves}
\date{}

\begin{document}
\maketitle
\begin{abstract}
In this paper, we consider local moves on classical and welded diagrams of string links, and the notion of welded extension of a classical move. Such extensions being non-unique in general, the idea is to find a topological criterion which could isolate one extension from the others. To that end, we turn to the relation between welded string links and knotted surfaces in $\RR^4$, and the ribbon subclass of these surfaces. This provides a topological interpretation of classical local moves as surgeries on surfaces, and of virtual local moves as surgeries on ribbon surfaces. Comparing these surgeries leads to the notion of ribbon residue of a classical local move, and we show that up to some broad conditions there can be at most one welded extension which is a ribbon residue. We provide three examples of ribbon residues, for the self-crossing change, the Delta and the band-pass moves. However, for the latter, we note that the given residue is actually not an extension of the band-pass move, showing that a classical move may have a ribbon residue and a welded extension, but no ribbon residue which is an extension.
\end{abstract}

\section*{Introduction}

Knot theory aims at studying embeddings of circles in $\RR^3$ up to ambient isotopies and, more generally, embeddings of codimension $2$ submanifolds in $\RR^n$.
As shown by K. Reidemeister in dimension $3$, and then extended to dimension $4$ by D. Roseman, this topology--based study can be translated in a combinatorial way through the use of diagrams, which are generic projections of the submanifold onto $\RR^{n-1}\times\{0\}\subset \RR^n$, up to local moves corresponding to some elementary local isotopies. Using this diagrammatic approach, one can extend the classical notion of knotted objects to the notion of welded objects, first defined for braids by R. Fenn, R. Rim\'anyi and C. Rourke in \cite{FenRimRou}. This welded theory is a quotient of the virtual extension, defined independently by L. Kauffman in \cite{Kauffman} and M. Goussarov, M. Polyak and O. Viro in \cite{GouPolVir} by allowing a new type of, so-called virtual, crossings on diagrams and a new type of, so-called detour, moves which freely trade any piece of strand supporting only virtual crossings for any other such virtual strand with same extremities. For welded objects, strands are also allowed to pass above (but not under) virtual crossings. Whereas the works of Kauffman \cite{Kauffman} and Goussarov--Polyak--Viro \cite{GouPolVir} are motivated by combinatorial aspects of link diagrams description, T. Brendle and A. Hatcher showed in \cite{BreHat} that the work of Fenn--Rim\'anyi--Rourke is much more related to 4--dimensional topology, as their braid-permutation groups are closely related to paths of circles configurations, which are surfaces in $\RR^4$ just like paths of points configurations are topological braids. As a matter of fact, welded link theory can be seen as an intermediary step between classical links and knotted surfaces. For general welded objects, the connection with knotted surfaces was made clear by S. Satoh \cite{Satoh} who extended the Tube map, first defined for classical objects by T. Yajima \cite{Yajima} by, roughly speaking, inflating strands into knotted tubes in $\RR^4$, to any such welded object. In particular, as it is the codomain of the Tube map, it emphasized the important role played by the ribbon subclass of knotted surfaces, corresponding to embedded surfaces which are the boundary of immersed solid handlebodies with only ribbon singularities. \\

Besides ambient isotopies, other topological quotients were combinatorially modelled using additional local diagrammatical moves. In \cite{fused}, B. Audoux, P. Bellingeri, J-B. Meilhan and E. Wagner started to study the question of potential welded extensions for such additional moves. Even though motivated by topological quotients in dimension 4, their study remained close to the classical knot theory side of the welded theory, a local move $M_w$ on welded diagrams being indeed said to extend a given local move $M_c$ on classical diagrams if two classical diagrams are related up to $M_w$ and welded Reidemeister moves if and only if they are related up to $M_c$ and classical Reidemeister moves. Surprisingly enough, it appeared that there exists classical local moves, {\it e.g.} the $\Delta$ move (see Figure 8), admitting multiple non-equivalent welded extensions. \\

The main goal of the present paper is to resolve such ambiguities by making the study of  welded extensions closer to topology, using the knotted surface theory side of the welded theory. Indeed, another (actually equivalent) way to relate classical knots with knotted surfaces is to spin a 1--dimensional knotted object in $\RR^3\subset\RR^4$ around a plane to obtain a surface. When similarly spinning a classical diagram, one obtains a broken surface diagram, the 4--dimensional counterpart of link diagrams, and when spinning a classical local move $M_c$, one obtains a surgery operation Spun($M_c$) which modifies in an explicit way broken surface diagrams inside some solid torus, and hence knotted surfaces inside some $S^1\times B^3\subset \RR^4$. A local move $M_w$ is then said to be a ribbon residue of $M_c$ if two ribbon surfaces $S_1$ and $S_2$ are related by Spun($M_c$) surgeries as knotted surfaces if and only if they can be realized as $S_1=\text{Tube}(L_1)$ and $S_2=\text{Tube}(L_2)$ with $L_1$ and $L_2$ being related by $M_w$ and welded Reidemeister moves. \\

As in \cite{fused}, we focus on the string link case, which are embedded intervals with prescribed fixed ends, and consider specifically three local moves, namely $SC$ which models link-homotopy, $\Delta$  which models link-homology and $BP$ which models band-passing (see Figure 8). More precisely:
\begin{itemize}
\item for $SC$, it was proven in \cite{fused} that the self-virtualization move $SV$, which turns any classical crossing involving portions of the same strand to a virtual crossing, is a welded extension. Without surprise, we prove that it is also a ribbon residue (Theorem~\ref{SVresSC});
\item for $\Delta$, it was proven in \cite{fused} that both the fused move $F$, which allows any strand to pass under classical crossings, and the virtual conjugation move $VC$, which surrounds a classical crossing by two virtual ones, are welded extensions. We prove that $F$ is a ribbon residue while $VC$ is not (Theorem~\ref{FresDelta}). This provides a way to designate $F$ as a preferred welded extension, carrying more topological meaning;
\item for $BP$, we prove (Theorem~\ref{VresBP}) that the action of Spun($BP$) of knotted surfaces is much stronger than that of $BP$ on classical string links, as its ribbon residue is the $V$ move, which trivializes welded string links. Even after restricting the action of Spun($BP$) to avoid some artefacts, the ribbon residue obtained (Theorem~\ref{BVpsresBP}) still trivializes classical string links. As a result, no welded extension of $BP$ which w-generates the $SV$ move can be a ribbon residue, including the one given in \cite{fused}.
\end{itemize}

The paper is organized as follows. In Section 1, we set the global background: the general notation is set in Section 1.1, welded knot theory and its relationship with ribbon surfaces are presented in Section 1.2, and local moves, welded extensions and ribbon residues are defined in Section 1.3. Section 2 is devoted to the above-mentionned moves, $SC$ in Section 2.1, $\Delta$ in Section 2.2 and $BP$ in Section 2.3. \\

\textit{Aknowledgements.} This article was inspired by results obtained during the redaction of my master's degree thesis. I would like to thank my thesis tutor B. Audoux for his guidance and helpful advice in the writing of this paper. I am also grateful to l'Institut de Mathématiques de Marseille for hosting me during my master research project, which was supported by l'\'{E}cole Normale Supérieure de Cachan.

\section{Settings}

\subsection{Notation}

We begin by introducing some notation. Let $n$ and $d$ be positive integers. We denote by $I:=[0,1]$ the unit interval, $B^d$ the closed unit ball in $\RR^d$ and $S^{d}=\dd B^{d+1}$ the $d$-dimensional sphere. Let $B^{d,1}:=B^d\times I$, $S^{d,1}:=S^d\times I$, and $\dd_{\eps}B^{d,1}:=B^d\times\{\eps\}$, $\dd_{\eps}S^{d,1}:=S^d\times\{\eps\}$ for $\eps=0,1$. The manifolds $I$ and $B^d$ are given their usual orientation (induced by the canonical orientation of $\RR^d$), $S^d$ is oriented as the boundary of $B^{d+1}$, and $B^{d,1}$, $S^{d,1}$ are given the product orientation. Manifolds and maps are always in the smooth category. \\

We will work with submanifolds of $B^{d,1}$ which have a fixed cartesian product structure near $\dd_0B^{d,1}\cup\dd_1B^{d,1}$. More precisely, let $X$ be a manifold, and $b:X\to\mathring{B}^d$ an embedding. We will consider embeddings (resp. immersions) $f:X\times I\to B^{d,1}$ for which there exists $\delta\in(0,1)$ such that:
\begin{itemize}
\item $f(x,t)=(b(x),t)$ for $t\in[0,\delta)\cup(1-\delta,1]$;
\item $f(X\times[\delta,1-\delta])\subset\mathring{B}^{d,1}$.
\end{itemize}
We call the image $Y=f(X\times I)$ an embedded (resp. immersed) submanifold of $B^{d,1}$. We denote by $\dd_{\eps}Y:=f(X\times\{\eps\})$ for $\eps=0,1$ and $\dd_{\ast}Y:=f(\dd X\times I)$ the lower, upper and lateral boundaries of $Y$ respectively. \\

In what follows, we will consider sets of such submanifolds for a fixed oriented $X$ (typically a disjoint union of balls or spheres) and a fixed embedding $b$. Thanks to the boundary condition, we can define the stacking product $Y_1\bullet Y_2$ for $Y_i=f_i(X\times I)$, $i=1,2$, by:
$$Y_1\bullet Y_2=f(X\times I),\quad f(x,t)=\left\{ \begin{array}{ll} f_1(x,2t) & \textrm{if } t\in[0,\frac{1}{2}], \\ f_2(x,2t-1) & \textrm{if } t\in[\frac{1}{2},1]. \end{array} \right.$$
To preserve the cartesian product structure of submanifols near the lower and upper boundaries, we will only consider isotopies of $B^{d,1}$ which are the identity in a neighborhood of $\dd_0B^{d,1}\cup\dd_1B^{d,1}$. Up to these isotopies, the stacking product is associative. \\

Let $p_1<\cdots<p_n$ be $n$ ordered points in the interval $(-1,1)$, which are fixed once and for all (for example take $p_i=(2i-1-n)/n$). Moreover, let $b_D:B_1^2\sqcup\cdots\sqcup B_n^2\to\mathring{B}^3$ be an embedding of $n$ disjoint disks, and $b_C:S_1^1\sqcup\cdots\sqcup S_n^1\to\mathring{B^3}$ its restriction to the circles $S_i^1:=\dd B_i^2$. We will use these as our fixed $b$ in Definitions~\ref{string2link} and~\ref{3ribbon}. \\

We will also make use of some algebraic notions: for a group $G$ normally generated by some elements $x_1,\ldots,x_n$, we denote by $RG$ the reduced group defined as the quotient of $G$ by the normal subgroup generated by the commutators $[x_i,gx_ig^{-1}]$ for $1\leq i\leq n$ and $g\in G$. It is the biggest quotient of $G$ in which the $x_i$'s commute with their conjugates. \\

For a group $G$ normally generated by $x_1,\ldots,x_n$, we denote by $\text{End}_C(G)$ (resp. $\text{Aut}_C(G)$) the set of conjugating endomorphisms (resp. automorphisms) of $G$, i.e. the subset of $\text{End}(G)$ (resp. $\text{Aut}(G)$) whose elements send each $x_i$ to one of its conjugates. We also define $\text{Aut}_C^0(G)$ as the subset of $\text{Aut}_C(G)$  whose elements send the product $x_1\cdots x_n$ to itself.

\subsection{Welded theory}

\subsubsection{Definition}

\begin{de}\label{stringlink}
An \emph{$n$--component string link} is an embedding of $\underset{1\leq i\leq n}{\sqcup}(I_i=\{i\}\times I)$ in $B^{2,1}$ with $b(i)=(0,p_i)\in B^2$. The $I_i$ are called the strands of the string link, and are oriented from $\dd_0I_i$ to $\dd_1I_i$. We denote by $\SL$ the set of string links up to isotopy. It is given a monoid structure by the stacking product.
\end{de}

A string link can be represented in two dimensions by taking a generic projection on a plane, where the singularities are transverse double points, called crossings. These crossings are represented by erasing part of the lower strand, and given a sign according to the orientation of the strands as indicated below: \\

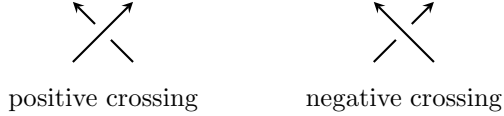
\begin{center}
\begin{tikzpicture}
\draw [thick] [>=stealth,->] (-1.6,0) -- (-2.4,0.8) ;
\fill [white] (-2,0.4) circle (0.15) ;
\draw [thick] [>=stealth,->] (-2.4,0) -- (-1.6,0.8) ;
\draw (-2,-0.5) node{positive crossing} ;
\draw [thick] [>=stealth,->] (1.6,0) -- (2.4,0.8) ;
\fill [white] (2,0.4) circle (0.15) ;
\draw [thick] [>=stealth,->] (2.4,0) -- (1.6,0.8) ;
\draw (2,-0.5) node{negative crossing} ;
\end{tikzpicture}
\end{center}

Welded string links can be defined using this diagrammatic approach. First, we need to consider a third type of crossing, called \emph{virtual crossing}:

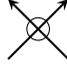
\begin{center}
\begin{tikzpicture}
\draw [thick] [>=stealth,->] (0.8,0) -- (0,0.8) ;
\draw (0.4,0.4) circle (0.15) ;
\draw [thick] [>=stealth,->] (0,0) -- (0.8,0.8) ;
\end{tikzpicture}
\end{center}

\begin{de}\label{sld}
An \emph{$n$--component virtual string link diagram} is an immersion of $\underset{1\leq i\leq n}{\sqcup}I_i$ in $B^{1,1}$ such that:
\begin{itemize}
\item $b(i)=p_i$, and $I_i$ is oriented from $(p_i,0)$ to $(p_i,1)$;
\item there is a finite number of singularities, which are transverse double points;
\item each double point is labelled to indicate a positive, negative or virtual crossing.
\end{itemize}
We denote by $vSLD_n$ the set of $n$--component virtual string link diagrams up to isotopy and reparametrization. It is given a monoid structure by the stacking product. We denote by $SLD_n$ the subset of $vSLD_n$ composed of diagrams with no virtual crossing, which are called classical diagrams. \\

For a diagram $D\in vSLD_n$, the portions of strands delimited by the undercrossings are called \emph{arcs}. If $D$ is a classical diagram, the arcs are simply the connected components obtained after erasing parts of the lower strands as described above. The arcs connected to $\dd_0B^{1,1}$ are called the bottom arcs, and the ones connected to $\dd_1B^{1,1}$ are called the top arcs.
\end{de}

As proven by Reidemeister (see \cite{Reidemeister} in the case of knots and links, which extends to string links), two classical diagrams represent the same string link if and only if one can be obtained from the other by applying some local moves, called Reidemeister moves, illustrated in Figure 1. \\

\begin{figure}[h]
\centering
\begin{tikzpicture}
\draw [thick] (0.5,0.5) .. controls +(0,-0.3) and +(0.2,-0.4) .. (0.1,0.5) ;
\draw [thick] (0.1,0.5) .. controls +(-0.1,0.2) and +(0,-0.1) .. (0,1) ;
\draw [white,fill=white] (0.1,0.5) circle (0.15) ;
\draw [thick] (0,0) .. controls +(0,0.1) and +(-0.1,-0.2) .. (0.1,0.5) ;
\draw [thick] (0.1,0.5) .. controls +(0.2,0.4) and +(0,0.3) .. (0.5,0.5) ;
\draw [<->] (1,0.5) -- (1.8,0.5) ;
\draw (1.4,0.8) node{R1} ;
\draw [thick] (2.3,0) -- (2.3,1) ;
\draw [<->] (2.8,0.5) -- (3.6,0.5) ;
\draw (3.2,0.8) node{R1} ;
\draw [thick] (4.1,0) .. controls +(0,0.1) and +(-0.1,-0.2) .. (4.2,0.5) ;
\draw [thick] (4.2,0.5) .. controls +(0.2,0.4) and +(0,0.3) .. (4.6,0.5) ;
\draw [white,fill=white] (4.2,0.5) circle (0.15) ;
\draw [thick] (4.6,0.5) .. controls +(0,-0.3) and +(0.2,-0.4) .. (4.2,0.5) ;
\draw [thick] (4.2,0.5) .. controls +(-0.1,0.2) and +(0,-0.1) .. (4.1,1) ;

\begin{scope}[xshift=6cm]
\draw [thick] (0,0) -- (0,1) ;
\draw [thick] (0.5,0) -- (0.5,1) ;
\draw [<->] (1,0.5) -- (1.8,0.5) ;
\draw (1.4,0.8) node{R2} ;
\draw [thick] (2.8,0) .. controls +(-0.2,0.1) and +(0,-0.2) .. (2.4,0.5) ;
\draw [thick] (2.4,0.5) .. controls +(0,0.2) and +(-0.2,-0.1) .. (2.8,1) ;
\draw [white,fill=white] (2.55,0.8) circle (0.15) ;
\draw [white,fill=white] (2.55,0.2) circle (0.15) ;
\draw [thick] (2.3,0) .. controls +(0.2,0.1) and +(0,-0.2) .. (2.7,0.5) ;
\draw [thick] (2.7,0.5) .. controls +(0,0.2) and +(0.2,-0.1) .. (2.3,1) ;
\end{scope}

\begin{scope}[xshift=10.2cm,yshift=0.75cm]
\draw [thick] (0.85,0) -- ++(-0.5,0.87) ;
\draw [white,fill=white] (0.71,0.25) circle (0.1) ;
\draw [white,fill=white] (0.5,0.61) circle (0.1) ;
\draw [thick] (0.15,0) -- ++(0.5,0.87) ;
\draw [white,fill=white] (0.29,0.25) circle (0.1) ;
\draw [thick] (0,0.25) -- (1,0.25) ;
\draw [<->] (1.5,0.43) -- (2.3,0.43) ;
\draw (1.9,0.73) node{R3} ;
\draw [thick] (2.95,0.87) -- ++(0.5,-0.87) ;
\draw [white,fill=white] (3.09,0.62) circle (0.1) ;
\draw [white,fill=white] (3.3,0.26) circle (0.1) ;
\draw [thick] (3.65,0.87) -- ++(-0.5,-0.87) ;
\draw [white,fill=white] (3.51,0.62) circle (0.1) ;
\draw [thick] (2.8,0.62) -- (3.8,0.62) ;
\end{scope}

\begin{scope}[xshift=10.2cm,yshift=-0.62cm]
\draw [thick] (0.15,0) -- ++(0.5,0.87) ;
\draw [white,fill=white] (0.5,0.61) circle (0.1) ;
\draw [thick] (0.85,0) -- ++(-0.5,0.87) ;
\draw [white,fill=white] (0.71,0.25) circle (0.1) ;
\draw [white,fill=white] (0.29,0.25) circle (0.1) ;
\draw [thick] (0,0.25) -- (1,0.25) ;
\draw [<->] (1.5,0.43) -- (2.3,0.43) ;
\draw (1.9,0.73) node{R3} ;
\draw [thick] (3.65,0.87) -- ++(-0.5,-0.87) ;
\draw [white,fill=white] (3.3,0.26) circle (0.1) ;
\draw [thick] (2.95,0.87) -- ++(0.5,-0.87) ;
\draw [white,fill=white] (3.09,0.62) circle (0.1) ;
\draw [white,fill=white] (3.51,0.62) circle (0.1) ;
\draw [thick] (2.8,0.62) -- (3.8,0.62) ;
\end{scope}
\end{tikzpicture}
\caption{Reidemeister moves}
\end{figure}
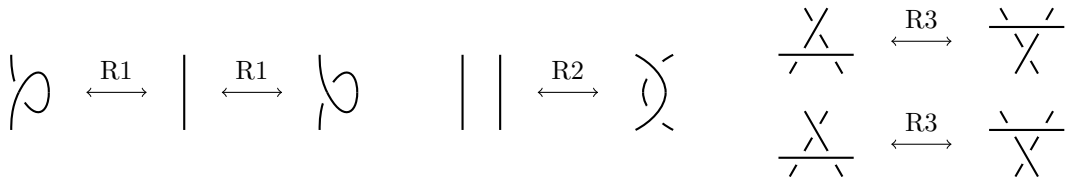

We use this diagrammatical approach to define welded string links, by considering virtual string link diagrams up to some local moves. We keep the (classical) Reidemeister moves, but also add new moves involving virtual crossings, illustrated in Figure 2. \\

\begin{figure}[h]
\centering
\begin{tikzpicture}
\draw [thick] (0,0) -- (0,1) ;
\draw [<->] (0.5,0.5) -- (1.3,0.5) ;
\draw (0.9,0.8) node{vR1} ;
\draw [thick] (1.8,0) .. controls +(0,0.1) and +(-0.1,-0.2) .. (1.9,0.5) ;
\draw [thick] (1.9,0.5) .. controls +(0.2,0.4) and +(0,0.3) .. (2.3,0.5) ;
\draw (1.9,0.5) circle (0.15) ;
\draw [thick] (2.3,0.5) .. controls +(0,-0.3) and +(0.2,-0.4) .. (1.9,0.5) ;
\draw [thick] (1.9,0.5) .. controls +(-0.1,0.2) and +(0,-0.1) .. (1.8,1) ;

\begin{scope}[xshift=3.7cm]
\draw [thick] (0,0) -- (0,1) ;
\draw [thick] (0.5,0) -- (0.5,1) ;
\draw [<->] (1,0.5) -- (1.8,0.5) ;
\draw (1.4,0.8) node{vR2} ;
\draw [thick] (2.8,0) .. controls +(-0.2,0.1) and +(0,-0.2) .. (2.4,0.5) ;
\draw [thick] (2.4,0.5) .. controls +(0,0.2) and +(-0.2,-0.1) .. (2.8,1) ;
\draw (2.55,0.8) circle (0.15) ;
\draw (2.55,0.2) circle (0.15) ;
\draw [thick] (2.3,0) .. controls +(0.2,0.1) and +(0,-0.2) .. (2.7,0.5) ;
\draw [thick] (2.7,0.5) .. controls +(0,0.2) and +(0.2,-0.1) .. (2.3,1) ;
\end{scope}

\begin{scope}[xshift=7.9cm]
\draw [thick] (0,0.29) -- (1.16,0.29) ;
\draw [thick] (0.17,0) -- ++(0.58,1);
\draw [thick] (0.99,0) -- ++(-0.58,1);
\draw (0.34,0.29) circle (0.12) ;
\draw (0.82,0.29) circle (0.12) ;
\draw (0.58,0.71) circle (0.12) ;
\draw [<->] (1.66,0.5) -- (2.46,0.5) ;
\draw (2.06,0.8) node{vR3} ;
\draw [thick] (2.96,0.71) -- (4.12,0.71) ;
\draw [thick] (3.95,1) -- ++(-0.58,-1);
\draw [thick] (3.13,1) -- ++(0.58,-1);
\draw (3.78,0.71) circle (0.12) ;
\draw (3.3,0.71) circle (0.12) ;
\draw (3.54,0.29) circle (0.12) ;
\end{scope}

\begin{scope}[xshift=1.15cm,yshift=-2cm]
\draw [thick] (0,0.29) -- (1.16,0.29) ;
\draw [thick] (0.99,0) -- ++(-0.58,1);
\draw [white,fill=white] (0.58,0.71) circle (0.12) ;
\draw [thick] (0.17,0) -- ++(0.58,1);
\draw (0.34,0.29) circle (0.12) ;
\draw (0.82,0.29) circle (0.12) ;
\draw [<->] (1.66,0.5) -- (2.46,0.5) ;
\draw (2.06,0.8) node{Mixed} ;
\draw [thick] (2.96,0.71) -- (4.12,0.71) ;
\draw [thick] (3.13,1) -- ++(0.58,-1);
\draw [white,fill=white] (3.54,0.29) circle (0.12) ;
\draw [thick] (3.95,1) -- ++(-0.58,-1);
\draw (3.78,0.71) circle (0.12) ;
\draw (3.3,0.71) circle (0.12) ;
\end{scope}

\begin{scope}[xshift=6.75cm,yshift=-2cm]
\draw [thick] (0.99,0) -- ++(-0.58,1);
\draw (0.58,0.71) circle (0.12) ;
\draw [thick] (0.17,0) -- ++(0.58,1);
\draw [white,fill=white] (0.34,0.29) circle (0.12) ;
\draw [white,fill=white] (0.82,0.29) circle (0.12) ;
\draw [thick] (0,0.29) -- (1.16,0.29) ;
\draw [<->] (1.66,0.5) -- (2.46,0.5) ;
\draw (2.06,0.8) node{OC} ;
\draw [thick] (3.13,1) -- ++(0.58,-1);
\draw (3.54,0.29) circle (0.12) ;
\draw [thick] (3.95,1) -- ++(-0.58,-1);
\draw [white,fill=white] (3.78,0.71) circle (0.12) ;
\draw [white,fill=white] (3.3,0.71) circle (0.12) ;
\draw [thick] (2.96,0.71) -- (4.12,0.71) ;
\end{scope}
\end{tikzpicture}
\caption{Additional moves on virtual diagrams}
\end{figure}
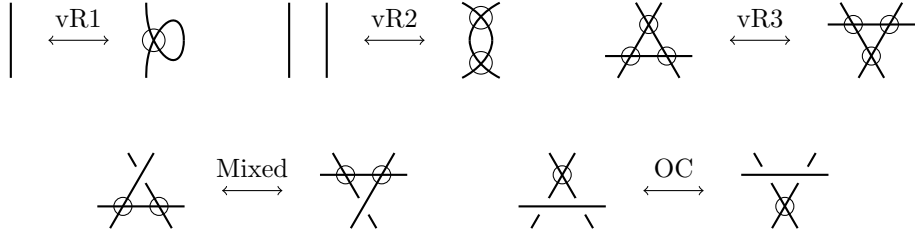

We denote by Reid (resp. vReid) the classical (resp. virtual) Reidemeister moves $R1$, $R2$, $R3$ (resp. $vR1$, $vR2$, $vR3$), and by wReid the welded Reidemeister moves, consisting of Reid, vReid, Mixed and $OC$. These welded Reidemeister moves enable what is called the \emph{detour move}: if a portion of a strand only involves virtual crossings, it can be changed for any other portion of strand involving only virtual crossings and having the same extremities.

\begin{de}
An \emph{$n$--component welded string link} is an equivalence class of $vSLD_n$ under wReid. We denote by $w\SL$ the monoid of $n$--component welded string links.
\end{de}

Welded string links can be represented in a more combinatorial way by Gauss diagrams.

\begin{de}
A \emph{Gauss diagram} on $n$ strands $\underset{1\leq i\leq n}{\sqcup}I_i$ is a finite set of triplets $(t,h,\eps)\in(\underset{1\leq i\leq n}{\sqcup}I_i)^2\times\{\pm 1\}$ such that the $t$'s and the $h$'s are all distinct. These triplets are called \emph{arrows}, with a \emph{tail} $t$ and a \emph{head} $h$ positionned on the $n$ strands, and a \emph{sign} $\eps$. We denote by $GD_n$ the set of Gauss diagrams up to isotopy.
\end{de}

A virtual diagram can be described by a Gauss diagram by associating an arrow to each classical crossing, the tail (resp. the head) indicating the position of the preimage on the upper (resp. lower) strand, and the sign indicating the type of crossing. Virtual crossings are not represented. \\

\begin{figure}[h]
\centering
\begin{tikzpicture}
\draw [thick] (0.75,0.5) .. controls +(-0.7,0) and +(0,-0.5) .. (0,2) ;
\draw [white,fill=white] (0.14,1) circle (0.15) ;
\draw [thick] (0,0) .. controls +(0,0.5) and +(-0.7,0) .. (0.75,1.5) ;
\draw [thick] (0.75,1.5) .. controls +(0.7,0) and +(0.7,0) .. (0.75,0.5) ;
\draw [thick,>=stealth,->] (0,2) -- (0,2.5) ;
\draw [white,fill=white] (1.12,0.64) circle (0.12) ;
\draw [thick] (1.5,0) .. controls +(0,0.5) and +(0,-0.5) .. (0.8,1.2) ;
\draw [white,fill=white] (0.86,0.965) circle (0.12) ;
\draw [thick] (0.75,0) .. controls +(0,0.5) and +(-0.4,-0.2) .. (1.1,1.2) ;
\draw (0.75,0.5) circle (0.12) ;
\draw [thick] (1.3,1.3) .. controls +(0.1,0.05) and +(0.1,-0.2) .. (1.4,1.7) ;
\draw [thick,>=stealth,->] (0.8,1.2) .. controls +(0,0.5) and +(0,-0.5) .. (1.5,2.5) ;
\draw [white,fill=white] (1.22,1.92) circle (0.12) ;
\draw [thick,>=stealth,->] (1.4,1.7) .. controls +(-0.2,0.4) and +(0,-0.5) .. (0.75,2.5) ;
\draw (0.87,1.48) circle (0.12) ;

\begin{scope}[xshift=4cm]
\draw [ultra thick,>=stealth,->] (0,0) -- (0,2.5) ;
\draw [ultra thick,>=stealth,->] (1,0) -- (1,2.5) ;
\draw [ultra thick,>=stealth,->] (2,0) -- (2,2.5) ;
\draw [>=stealth,->] (0,0.3) .. controls +(-0.5,0) and +(-0.5,0) .. (0,2) ;
\draw (-0.3,2.2) node{$+$} ;
\draw [>=stealth,->] (0,0.75) -- (1,0.75) ;
\draw (0.75,0.5) node{$+$} ;
\draw [>=stealth,->] (1,2) -- (2,2) ;
\draw (1.7,2.2) node{$-$} ;
\draw [>=stealth,->] (2,0.5) .. controls +(-0.5,0) and +(1,0) .. (0,1.4) ;
\draw (0.3,1.6) node{$+$} ;
\draw [>=stealth,->] (1,0.3) .. controls +(0.5,0) and +(-0.5,0) .. (2,1.3) ;
\draw (1.7,1.5) node{$+$} ;
\end{scope}
\end{tikzpicture}
\caption{A virtual diagram, and the associated Gauss diagram}
\end{figure}
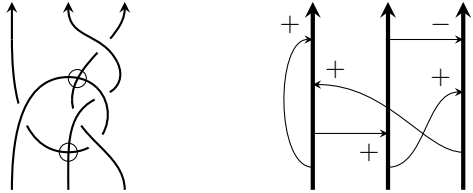

Similarly to virtual string link diagrams, we need to allow local moves on Gauss diagrams in order to obtain a one-to-one correspondance with welded string links. The moves on Gauss diagrams corresponding to the welded Reidemeister moves are illustrated in Figure 4. Since the virtual Reidemeister moves only involve virtual crossings, they do not affect Gauss diagrams, and neither does the Mixed move. The $R3$ move is labelled with a $(\ast)$ to indicate that it must satisfy some sign conditions: to apply $R3$, we must have $\delta_1\eps_1=\delta_2\eps_2=\delta_3\eps_3$, where $\delta_i=1$ or $-1$ depending on whether the $i^{th}$ portion of strand is oriented upward or downward. \\

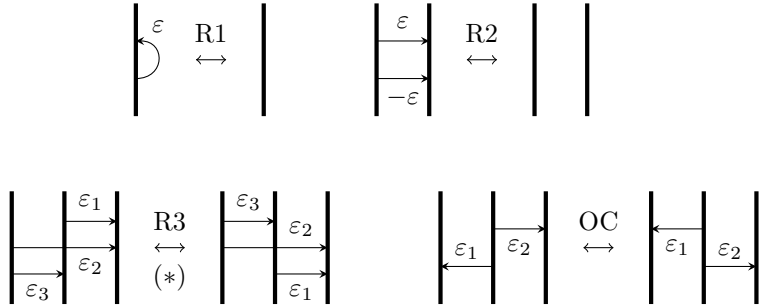
\begin{figure}[h]
\centering
\begin{tikzpicture}
\begin{scope}[xshift=1.65cm]
\draw [ultra thick] (0,0) -- ++(0,1.5) ;
\draw [>=stealth,->] (0,0.5) .. controls +(0.4,0) and +(0.4,0) .. (0,1) ;
\draw (0.3,1.2) node{$\eps$} ;
\draw [<->] (0.8,0.75) -- (1.2,0.75) ;
\draw (1,1.1) node{R1} ;
\draw [ultra thick] (1.7,0) -- ++(0,1.5) ;
\end{scope}

\begin{scope}[xshift=4.85cm]
\draw [ultra thick] (0,0) -- ++(0,1.5) ;
\draw [ultra thick] (0.7,0) -- ++(0,1.5) ;
\draw [>=stealth,->] (0,0.5) -- ++(0.7,0) ;
\draw (0.35,0.25) node{$-\eps$} ;
\draw [>=stealth,->] (0,1) -- ++(0.7,0) ;
\draw (0.35,1.25) node{$\eps$} ;
\draw [<->] (1.2,0.75) -- (1.6,0.75) ;
\draw (1.4,1.1) node{R2} ;
\draw [ultra thick] (2.1,0) -- ++(0,1.5) ;
\draw [ultra thick] (2.8,0) -- ++(0,1.5) ;
\end{scope}

\begin{scope}[xshift=0cm,yshift=-2.5cm]
\draw [ultra thick] (0,0) -- ++(0,1.5) ;
\draw [ultra thick] (0.7,0) -- ++(0,1.5) ;
\draw [ultra thick] (1.4,0) -- ++(0,1.5) ;
\draw [>=stealth,->] (0,0.4) -- ++(0.7,0) ;
\draw (0.35,0.15) node{$\eps_3$} ;
\draw [>=stealth,->] (0,0.75) -- ++(1.4,0) ;
\draw (1.05,0.5) node{$\eps_2$} ;
\draw [>=stealth,->] (0.7,1.1) -- ++(0.7,0) ;
\draw (1.05,1.35) node{$\eps_1$} ;
\draw [<->] (1.9,0.75) -- (2.3,0.75) ;
\draw (2.1,1.1) node{R3} ;
\draw (2.1,0.35) node{$(\ast)$} ;
\draw [ultra thick] (2.8,0) -- ++(0,1.5) ;
\draw [ultra thick] (3.5,0) -- ++(0,1.5) ;
\draw [ultra thick] (4.2,0) -- ++(0,1.5) ;
\draw [>=stealth,->] (2.8,1.1) -- ++(0.7,0) ;
\draw (3.15,1.35) node{$\eps_3$} ;
\draw [>=stealth,->] (2.8,0.75) -- ++(1.4,0) ;
\draw (3.85,1) node{$\eps_2$} ;
\draw [>=stealth,->] (3.5,0.4) -- ++(0.7,0) ;
\draw (3.85,0.15) node{$\eps_1$} ;
\end{scope}

\begin{scope}[xshift=5.7cm,yshift=-2.5cm]
\draw [ultra thick] (0,0) -- ++(0,1.5) ;
\draw [ultra thick] (0.7,0) -- ++(0,1.5) ;
\draw [ultra thick] (1.4,0) -- ++(0,1.5) ;
\draw [>=stealth,->] (0.7,0.5) -- ++(-0.7,0) ;
\draw (0.35,0.7) node{$\eps_1$} ;
\draw [>=stealth,->] (0.7,1) -- ++(0.7,0) ;
\draw (1.05,0.75) node{$\eps_2$} ;
\draw [<->] (1.9,0.75) -- (2.3,0.75) ;
\draw (2.1,1.1) node{OC} ;
\draw [ultra thick] (2.8,0) -- ++(0,1.5) ;
\draw [ultra thick] (3.5,0) -- ++(0,1.5) ;
\draw [ultra thick] (4.2,0) -- ++(0,1.5) ;
\draw [>=stealth,->] (3.5,1) -- ++(-0.7,0) ;
\draw (3.15,0.75) node{$\eps_1$} ;
\draw [>=stealth,->] (3.5,0.5) -- ++(0.7,0) ;
\draw (3.85,0.7) node{$\eps_2$} ;
\end{scope}
\end{tikzpicture}
\caption{Local moves on Gauss diagrams}
\end{figure}

It is well known and straightforwardly checked that up to these local moves, virtual diagrams and Gauss diagrams are faithful representations of string links:

\begin{prop}
The following monoid isomorphisms hold:
\begin{itemize}
\item $\SL\simeq SLD_n/\{\text{Reid}\}$;
\item $vSLD_n/\{\text{vReid},\text{Mixed}\}\simeq GD_n$;
\item $w\SL=vSLD_n/\{\text{wReid}\}\simeq GD_n/\{\text{Reid},OC\}$.
\end{itemize}
\end{prop}

\subsubsection{Relation with knotted surfaces}

String links can be related to knotted surfaces through two maps, called \emph{Spun} and \emph{Tube}. The Spun map spins a classical string link around a plane in 4 dimensions to obtain a surface, while the Tube map, first defined for classical knots by T. Yajima in \cite{Yajima} and then extended to the welded case by S. Satoh in \cite{Satoh}, inflates a welded string link and takes the boundary to obtain \textquotedblleft tubes". One important fact is that the Tube map sends welded string links to the ribbon subclass (see Definition~\ref{ribbon2SL} below) of the surfaces considered here. \\

A knotted surface is an embedding of a surface in $\RR^4$. As in the case of knots, such surfaces can be projected on a hyperplane, in order to obtain a surface in $\RR^3$ with three types of singularities (see \cite{Roseman} or \cite{S4S}): lines of double points (where it is locally the intersection of two planes), isolated triple points (locally the intersection of three planes in a point) and isolated branch points (where the projection is not an immersion).

\begin{center}
\begin{tikzpicture}
\draw [gray,dashed] (0,0,-1) -- (-0.78,0,-1) ;
\draw [gray,dashed] (0,0,-1) -- (0,-0.78,-1) ;
\draw (-1,0,-1) -- (-1,0,1) ;
\draw (-1,0,1) -- (1,0,1) ;
\draw (1,0,1) -- (1,0,-1) ;
\draw (1,0,-1) -- (0,0,-1) ;
\draw (-0.78,0,-1) -- (-1,0,-1) ;
\draw (0,-1,-1) -- (0,-1,1) ;
\draw (0,-1,1) -- (0,1,1) ;
\draw (0,1,1) -- (0,1,-1) ;
\draw (0,1,-1) -- (0,0,-1) ;
\draw (0,-0.78,-1) -- (0,-1,-1) ;
\draw [thick] (0,0,1) -- (0,0,-1) ;
\draw (0,-2) node{line of double points} ;

\begin{scope}[xshift=5cm]
\draw [gray,dashed] (0.63,0,-1) -- (0,0,-1) ;
\draw [gray!50,dashed] (0,0,-1) -- (-0.78,0,-1) ;
\draw [gray,dashed] (-0.78,0,-1) -- (-1,0,-1) ;
\draw [gray,dashed] (-1,0,-1) -- (-1,0,0) ;
\draw [gray,dashed] (0,-1,0) -- (0,-1,-1) ;
\draw [gray,dashed] (0,-1,-1) -- (0,-0.78,-1) ;
\draw [gray!50,dashed] (0,-0.78,-1) -- (0,0,-1) ;
\draw [gray,dashed] (0,0,-1) -- (0,0.62,-1) ;
\draw [gray,dashed] (0,-1,0) -- (-0.37,-1,0) ;
\draw [gray,dashed] (-1,-0.38,0) -- (-1,0,0) ;
\draw [gray,dashed,thick] (0,0,0) -- (0,0,-1) ;
\draw [gray,dashed,thick] (0,0,0) -- (0,-0.38,0) ;
\draw [gray,dashed,thick] (0,0,0) -- (-0.38,0,0) ;
\draw (-1,0,1) -- (1,0,1) ;
\draw (1,0,1) -- (1,0,-1) ;
\draw (1,0,-1) -- (0.63,0,-1) ;
\draw (-1,0,0) -- (-1,0,1) ;
\draw (0,1,-1) -- (0,1,1) ;
\draw (0,1,1) -- (0,-1,1) ;
\draw (0,-1,1) -- (0,-1,0) ;
\draw (0,0.62,-1) -- (0,1,-1) ;
\draw (-1,1,0) -- (1,1,0) ;
\draw (1,1,0) -- (1,-1,0) ;
\draw (1,-1,0) -- (0,-1,0) ;
\draw (-0.37,-1,0) -- (-1,-1,0) ;
\draw (-1,-1,0) -- (-1,-0.38,0) ;
\draw (-1,0,0) -- (-1,1,0) ;
\draw [thick] (0,0,1) -- (0,0,0) ;
\draw [thick] (0,1,0) -- (0,0,0) ;
\draw [thick] (0,-0.38,0) -- (0,-1,0) ;
\draw [thick] (1,0,0) -- (0,0,0) ;
\draw [thick] (-0.38,0,0) -- (-1,0,0) ;
\fill (0,0,0) circle (0.07) ;
\draw (0,-2) node{triple point} ;
\end{scope}

\begin{scope}[xshift=10.5cm,yshift=0.25cm]
\draw [gray,dashed] (0.5,-0.75,-1) -- (-1.5,-0.75,-1) ;
\draw [gray,dashed] (-1.5,-0.75,-1) .. controls +(0,0.2,0.25) and +(0,-0.2,-0.5) .. (-1.5,0,0) ;
\draw (1,-0.75,1) -- (1,0,0) ;
\draw (1,0,0) .. controls +(0,0.975,-1.3) and +(0,0.975,1.3) .. (1,0,0) ;
\draw (1,0,0) -- (1,-0.75,-1) ;
\draw [thick] (1,0,0) -- (-0.5,0,0) ;
\draw (-0.5,0,0) -- (-1.5,0,0) ;
\draw (1,-0.75,1) -- (-1.5,-0.75,1) ;
\draw (1,-0.75,-1) -- (0.5,-0.75,-1) ;
\draw (-1.5,-0.75,1) .. controls +(0,0.2,-0.25) and +(0,0.2,0.5) .. (-1.5,0,0) ;
\draw (1.05,0.775,0) -- (-0.5,0,0) ;
\fill (-0.5,0,0) circle (0.06) ;
\draw (-0.25,-2.25) node{branch point} ;
\end{scope}
\end{tikzpicture}
\end{center}

Such a projection, together with the information of upper/lower parts of the surface at singularities, is called a \emph{broken surface diagram} of the knotted surface. As in the 1--dimensional case, a broken surface diagram can be represented by deleting thin bands around the lines of double points on the lower part of the diagram. This is illustrated below, where the upper and lower parts of the diagrams have been chosen arbitrarily:

\begin{center}
\begin{tikzpicture}
\draw [gray,dashed] (-0.3,0,-1) -- (-0.78,0,-1) ;
\draw [gray,dashed] (-0.3,0,0.2) -- (-0.3,0,-1) ;
\draw [gray,dashed] (0,-0.3,-1) -- (0,-0.78,-1) ;
\draw (-1,0,-1) -- (-1,0,1) ;
\draw (-1,0,1) -- (-0.3,0,1) ;
\draw (-0.3,0,1) -- (-0.3,0,0.2) ;
\draw (0.3,0,1) -- (1,0,1) -- (1,0,-1) -- (0.3,0,-1) -- cycle ;
\draw (-0.78,0,-1) -- (-1,0,-1) ;
\draw (0,-1,-1) -- (0,-1,1) ;
\draw (0,-1,1) -- (0,1,1) ;
\draw (0,1,1) -- (0,1,-1) ;
\draw (0,1,-1) -- (0,-0.3,-1) ;
\draw (0,-0.78,-1) -- (0,-1,-1) ;
\draw (0,-2) node{line of double points} ;

\begin{scope}[xshift=5cm]
\draw (0,-0.25,-1) -- (0,-0.14,-1) ;
\draw [gray,dashed] (0,-0.14,-1) -- (0,0.62,-1) ;
\draw (0,0.62,-1) -- (0,1,-1) -- (0,1,1) -- (0,-1,1) -- (0,-1,-0.63) ;
\draw [gray,dashed] (0,-1,-0.63) -- (0,-1,-1) -- (0,-0.77,-1) ;
\draw [gray!50,dashed] (0,-0.77,-1) -- (0,-0.64,-1) ;
\draw [gray,dashed] (0,-0.64,-1) -- (0,-0.25,-1) ;
\draw (0.25,0,-0.65) -- (0.25,0,1) -- (1,0,1) -- (1,0,-1) -- (0.62,0,-1) ;
\draw [gray,dashed] (0.62,0,-1) -- (0.25,0,-1) -- (0.25,0,-0.65) ;
\draw (-1,0,-0.65) -- (-1,0,1) -- (-0.25,0,1) -- (-0.25,0,0.35) ;
\draw [gray,dashed] (-0.25,0,0.35) -- (-0.25,0,-1) -- (-0.63,0,-1) ;
\draw [gray!50,dashed] (-0.63,0,-1) -- (-0.77,0,-1) ;
\draw [gray,dashed] (-0.77,0,-1) -- (-1,0,-1) -- (-1,0,-0.65) ;
\draw (-0.39,0.25,0) -- (-1,0.25,0) -- (-1,1,0) -- (-0.25,1,0) -- (-0.25,0.75,0) ;
\draw [gray,dashed] (-0.25,0.75,0) -- (-0.25,0.25,0) -- (-0.39,0.25,0) ;
\draw (0.25,1,0) -- (1,1,0) -- (1,0.25,0) -- (0.25,0.25,0) -- cycle ;
\draw (0.75,-0.25,0) -- (1,-0.25,0) -- (1,-1,0) -- (0.25,-1,0) -- (0.25,-0.38,0) ;
\draw [gray,dashed] (0.75,-0.25,0) -- (0.25,-0.25,0) -- (0.25,-0.38,0) ;
\draw (-0.39,-1,0) -- (-1,-1,0) -- (-1,-0.38,0) ;
\draw [gray,dashed] (-1,-0.38,0) -- (-1,-0.25,0) -- (-0.5,-0.25,0) ;
\draw (-0.5,-0.25,0) -- (-0.39,-0.25,0) ;
\draw [gray,dashed] (-0.39,-0.25,0) -- (-0.25,-0.25,0) -- (-0.25,-1,0) -- (-0.39,-1,0) ;
\draw (0,-2) node{triple point} ;
\end{scope}

\begin{scope}[xshift=10.5cm,yshift=0.25cm]
\draw [gray,dashed] (0.5,-0.75,-1) -- (-1.5,-0.75,-1) ;
\draw [gray,dashed] (-1.5,-0.75,-1) .. controls +(0,0.2,0.25) and +(0,-0.2,-0.5) .. (-1.5,0,0) ;
\draw (1,-0.75,1) -- (1,-0.15,0.2) ;
\draw (1,0.15,-0.2) .. controls +(0,0.75,-1) and +(0,0.975,1.3) .. (1,0,0) ;
\draw (1,0,0) -- (1,-0.75,-1) ;
\draw (-0.5,0,0) -- (1,-0.15,0.2) ;
\draw [gray,dashed] (-0.5,0,0) -- (0.82,0.132,-0.176) ;
\draw (0.82,0.132,-0.176) -- (1,0.15,-0.2) ;
\draw (-0.5,0,0) -- (-1.5,0,0) ;
\draw (1,-0.75,1) -- (-1.5,-0.75,1) ;
\draw (1,-0.75,-1) -- (0.5,-0.75,-1) ;
\draw (-1.5,-0.75,1) .. controls +(0,0.2,-0.25) and +(0,0.2,0.5) .. (-1.5,0,0) ;
\draw (1.05,0.77,0) -- (-0.5,0,0) ;
\draw (-0.25,-2.25) node{branch point} ;
\end{scope}
\end{tikzpicture}
\end{center}

As in the case of knots, there are local moves on broken surface diagrams which identify different diagrams associated to the same knotted surface. These are called Roseman moves (see \cite{Roseman} for a detailed description). \\

We now define an analogue of string links in the case of surfaces.

\begin{de}\label{string2link}
A \emph{string 2--link} is an embedding of $\underset{1\leq i\leq n}{\sqcup}S_i^{1,1}=\underset{1\leq i\leq n}{\sqcup}(S_i^1\times I)$ in $B^{3,1}$ with $b=b_C$ as our fixed boundary embedding. We denote by $\dl{2}$ the set of string 2--links up to isotopy. It is given a monoid structure by the stacking product.
\end{de}

We will also consider the subclass of ribbon string 2--links, which is the string 2--link analogue of the ribbon subclass of knots.

\begin{de}\label{3ribbon}
A \emph{3--ribbon} is an immersion of $\underset{1\leq i\leq n}{\sqcup}B_i^{2,1}=\underset{1\leq i\leq n}{\sqcup}(B_i^2\times I)$ in $B^{3,1}$ with $b=b_D$, and a singular set composed of a finite number of ribbon singularities, which are defined as follows: a connected singularity is \emph{ribbon} if it is a disk $\delta$ given by a transverse intersection of the images of two components $B_i^{2,1}$ and $B_j^{2,1}$ (with possibly $i=j$), with preimages $\delta_c\subset B_i^{2,1}$ and $\delta_{ess}\subset B_j^{2,1}$ of $\delta$ satisfying the following conditions:
\begin{itemize}
\item $\delta_c\subset\mathring{B}_i^{2,1}$;
\item $\mathring{\delta}_{ess}\subset\mathring{B}_j^{2,1}$ and $\dd\delta_{ess}\subset\dd B_j^2\times I$ is non-trivial in $H_1(\dd B_j^2\times I)$.
\end{itemize}
We call $\delta_c$ the contractible preimage and $\delta_{ess}$ the essential preimage. \\

By considering the images of tangent vectors at the points $x_c\in\delta_c$ and $x_{ess}\in\delta_{ess}$ in the preimage of a point $x\in\delta$, we can associate a sign to a ribbon singularity. See {\cite[\S 3.2.1]{AudHdR}} for more details.
\end{de}

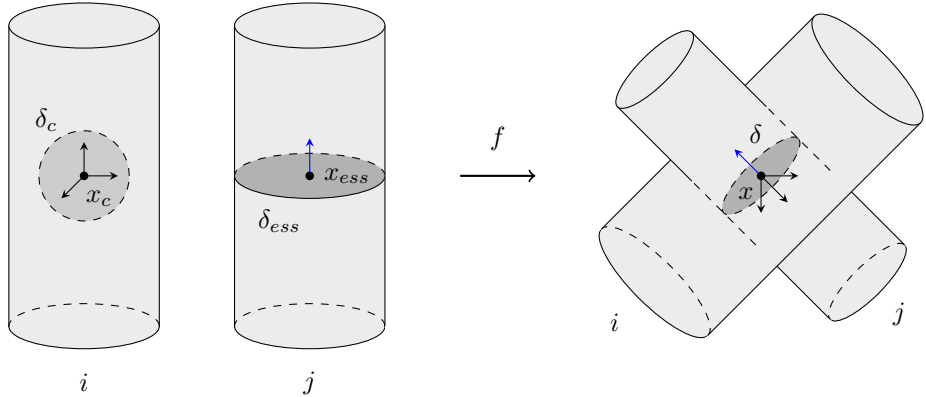
\begin{figure}[h]
\centering
\begin{tikzpicture}
\fill [gray!15] (0,-2) -- (0,2) -- (2,2) -- (2,-2) -- cycle ;
\begin{scope}
\clip (0,-2) -- (2,-2) -- (2,-3) -- (0,-3) -- cycle ;
\draw [fill=gray!15] (1,-2) ellipse (1 and 0.3) ;
\end{scope}
\begin{scope}
\clip (0,-2) -- (2,-2) -- (2,-1) -- (0,-1) -- cycle ;
\draw [dashed] (1,-2) ellipse (1 and 0.3) ;
\end{scope}
\draw (0,-2) -- (0,2) ;
\draw (2,-2) -- (2,2) ;
\draw [fill=gray!15] (1,2) ellipse (1 and 0.3) ;
\draw (1,-2.75) node{$i$} ;
\draw [dashed,fill=gray!40] (1,0) circle (0.6) ;
\draw (0.5,0.75) node{$\delta_c$} ;
\draw [fill=black] (1,0) circle (0.05) ;
\draw (1.2,-0.3) node{$x_c$} ;
\draw [>=stealth,->] (1,0) -- ++(-0.3,-0.3) ;
\draw [>=stealth,->] (1,0) -- ++(0.45,0) ;
\draw [>=stealth,->] (1,0) -- ++(0,0.45) ;

\begin{scope}[xshift=3cm]
\fill [gray!15] (0,-2) -- (0,2) -- (2,2) -- (2,-2) -- cycle ;
\begin{scope}
\clip (0,-2) -- (2,-2) -- (2,-3) -- (0,-3) -- cycle ;
\draw [fill=gray!15] (1,-2) ellipse (1 and 0.3) ;
\end{scope}
\begin{scope}
\clip (0,-2) -- (2,-2) -- (2,-1) -- (0,-1) -- cycle ;
\draw [dashed] (1,-2) ellipse (1 and 0.3) ;
\end{scope}
\draw (0,-2) -- (0,2) ;
\draw (2,-2) -- (2,2) ;
\draw [fill=gray!15] (1,2) ellipse (1 and 0.3) ;
\draw (1,-2.75) node{$j$} ;
\begin{scope}
\clip (0,0) -- (2,0) -- (2,-1) -- (0,-1) -- cycle ;
\draw [fill=gray!60] (1,0) ellipse (1 and 0.3) ;
\end{scope}
\begin{scope}
\clip (0,0) -- (2,0) -- (2,1) -- (0,1) -- cycle ;
\draw [dashed,fill=gray!60] (1,0) ellipse (1 and 0.3) ;
\end{scope}
\draw (0.6,-0.6) node{$\delta_{ess}$} ;
\draw (1.5,0) node{$x_{ess}$} ;
\draw [blue,>=stealth,->] (1,0) -- ++(0,0.5) ;
\draw [fill=black] (1,0) circle (0.05) ;
\end{scope}

\begin{scope}[xshift=6cm]
\draw [thick,->] (0,0) -- ++(1,0) ;
\draw (0.5,0.5) node{$f$} ;
\end{scope}

\begin{scope}[xshift=10cm,rotate=-45]
\fill [gray!15] (-1,-2) -- (-1,2) -- (1,2) -- (1,-2) -- cycle ;
\draw [fill=gray!15] (0,2) ellipse (1 and 0.3) ;
\end{scope}
\begin{scope}[xshift=10cm,rotate=-45]
\clip (-1,-2) -- (-1,-1) -- (1,-1) -- (1,-2) -- cycle ;
\draw [dashed,fill=gray!15] (0,-2) ellipse (1 and 0.3) ;
\end{scope}
\begin{scope}[xshift=10cm,rotate=-45]
\clip (-1,-2) -- (-1,-3) -- (1,-3) -- (1,-2) -- cycle ;
\draw [fill=gray!15] (0,-2) ellipse (1 and 0.3) ;
\draw (0,-2.75) node{$i$} ;
\end{scope}

\begin{scope}[xshift=10cm,rotate=45]
\fill [gray!15] (-0.7,-2) -- (-0.7,2) -- (0.7,2) -- (0.7,-2) -- cycle ;
\draw [fill=gray!15] (0,2) ellipse (0.7 and 0.21) ;
\end{scope}
\begin{scope}[xshift=10cm,rotate=45]
\clip (-0.7,-2) -- (-0.7,-1) -- (0.7,-1) -- (0.7,-2) -- cycle ;
\draw [dashed,fill=gray!15] (0,-2) ellipse (0.7 and 0.21) ;
\end{scope}
\begin{scope}[xshift=10cm,rotate=45]
\clip (-0.7,-2) -- (-0.7,-3) -- (0.7,-3) -- (0.7,-2) -- cycle ;
\draw [fill=gray!15] (0,-2) ellipse (0.7 and 0.21) ;
\draw (0,-2.6) node{$j$} ;
\end{scope}

\begin{scope}[xshift=10cm,rotate=-45]
\draw (-1,-2) -- (-1,-0.7) ;
\draw (-1,0.7) -- (-1,2) ;
\draw (1,-2) -- (1,2) ;
\end{scope}

\begin{scope}[xshift=10cm,rotate=45]
\draw (-0.7,-2) -- (-0.7,-1) ;
\draw [dashed] (-0.7,-0.6) -- (-0.7,0.6) ;
\draw (-0.7,0.6) -- (-0.7,2) ;
\draw (0.7,-2) -- (0.7,-1) ;
\draw [dashed] (0.7,-0.6) -- (0.7,0.6) ;
\draw (0.7,0.6) -- (0.7,2) ;
\draw [dashed,fill=gray!60] (0,0) ellipse (0.7 and 0.21) ;
\draw (0.35,0.45) node{$\delta$} ;
\draw (-0.3,-0.02) node{$x$} ;
\draw [>=stealth,->] (0,0) -- ++(-0.35,-0.35) ;
\draw [>=stealth,->] (0,0) -- ++(0,-0.5) ;
\draw [>=stealth,->] (0,0) -- ++(0.35,-0.35) ;
\draw [blue,>=stealth,->] (0,0) -- ++(0,0.5) ;
\draw [fill=black] (0,0) circle (0.05) ;
\end{scope}
\end{tikzpicture}
\caption{Ribbon singularity}
\end{figure}

\begin{de}\label{ribbon2SL}
A \emph{ribbon string 2--link} is a string 2--link $L$ which is the lateral boundary of a 3--ribbon $R$: $L=\dd_{\ast}R$. We say that $R$ is a \emph{ribbon filling} of $L$. We denote by $\dr{2}$ the monoid of ribbon string 2--links.
\end{de}

As described in Section 3.2 of \cite{AudHdR}, we can associate a Gauss diagram to each 3--ribbon, with arrows corresponding to ribbon singularities. This induces a one-to-one monoid homomorphism between 3--ribbons up to istotopy and Gauss diagrams up to the $OC$ move. The inverse of this homomorphism becomes invariant under Reidemeister moves when composed with the \textquotedblleft lateral boundary" map $\dd_{\ast}:\{\text{3--ribbon}\}/\{\text{isotopy}\}\to\dr{2}$, and induces a surjective homomorphism $\text{Tube}:GD_n/\{\text{Reid},OC\}\to\dr{2}$. \\

\begin{figure}[h]
\centering
\begin{tikzpicture}
\begin{scope}[scale=0.9]
\begin{scope}[rotate=45]
\draw (0,-1.1) ellipse (0.5 and 0.2) ;
\fill [white] (-0.5,-1) -- (-0.5,-0.5) -- (0.5,-0.5) -- (0.5,-1) -- cycle ;
\draw (-0.5,-2) -- (-0.5,-1.1) ;
\draw (0.5,-2) -- (0.5,-1.1) ;
\end{scope}
\begin{scope}[rotate=45]
\clip (-0.5,-2) -- (-0.5,-2.5) -- (0.5,-2.5) -- (0.5,-2) -- cycle ;
\draw (0,-2) ellipse (0.5 and 0.2) ;
\end{scope}
\begin{scope}[rotate=45]
\clip (-0.5,-2) -- (-0.5,-1.5) -- (0.5,-1.5) -- (0.5,-2) -- cycle ;
\draw [dashed] (0,-2) ellipse (0.5 and 0.2) ;
\end{scope}

\begin{scope}[rotate=-45]
\draw (0,2) ellipse (1 and 0.4) ;
\draw (-1,-2) -- (-1,2) ;
\draw (1,-2) -- (1,2) ;
\end{scope}
\begin{scope}[rotate=-45]
\clip (-1,-2) -- (-1,-3) -- (1,-3) -- (1,-2) -- cycle ;
\draw (0,-2) ellipse (1 and 0.4) ;
\end{scope}
\begin{scope}[rotate=-45]
\clip (-1,-2) -- (-1,-1) -- (1,-1) -- (1,-2) -- cycle ;
\draw [dashed] (0,-2) ellipse (1 and 0.4) ;
\end{scope}

\begin{scope}[rotate=45]
\draw (0,0.5) ellipse (0.75 and 0.3) ;
\fill [white] (-0.5,0.5) -- (-0.5,1.2) -- (0.5,1.2) -- (0.5,0.5) -- cycle ;
\draw [dashed] (0,0.5) ellipse (0.25 and 0.1) ;
\draw (0,2) ellipse (0.5 and 0.2) ;
\draw [dashed] (-0.5,-0.5) -- (-0.5,0.28) ;
\draw (-0.5,0.28) -- (-0.5,2) ;
\draw [dashed] (0.5,-0.5) -- (0.5,0.28) ;
\draw (0.5,0.28) -- (0.5,2) ;
\draw [dashed] (0,-0.5) ellipse (0.5 and 0.2) ;
\end{scope}
\end{scope}
\end{tikzpicture}
\caption{Broken surface diagram of a ribbon singularity}
\end{figure}
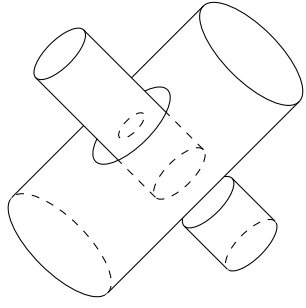

Figure 6 illustrates a broken surface diagram of a ribbon string 2--link at a ribbon singularity. We can give a more geometric definition of the Tube map in terms of broken surface diagrams:

\begin{de}\label{defTube}
For a welded string link $L\in w\SL$, let $D\in vSLD_n$ be a diagram of $L$, which we place in $B^{2,1}$ using the embedding $(y,z)\in B^{1,1}\mapsto(0,y,z)\in B^{2,1}$. Let $N$ be a tubular neighborhood of $D$ in $B^{2,1}$, and $\dd_{\ast}N$ its lateral boundary. At each crossing of $D$, we modify $\dd_{\ast}N$ as indicated in Figure 7: a positive (resp. negative) crossing gives a broken surface diagram of a positive (resp. negative) ribbon singularity, and a virtual crossing gives two disjoint tubes. The Tube map is then defined as sending $L$ to the element of $\dr{2}$ represented by this broken surface diagram.
\end{de}

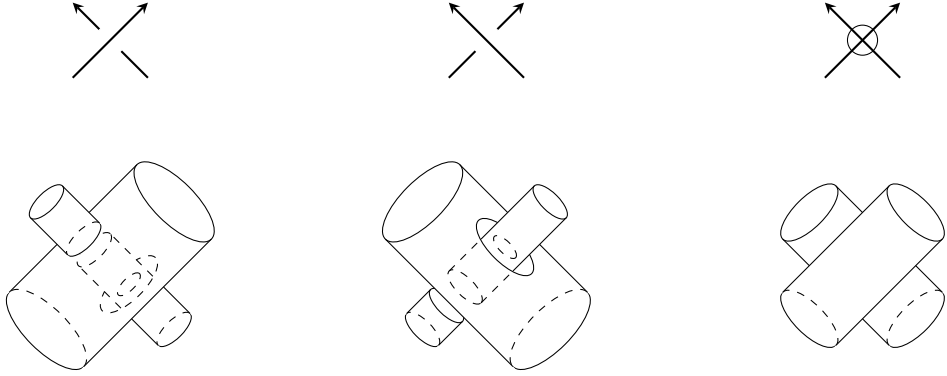
\begin{figure}[h]
\centering
\begin{tikzpicture}
\begin{scope}[xshift=-5cm,yshift=0cm]
\draw [thick,>=stealth,->] (0.5,-0.5) -- ++(-1,1) ;
\draw [white,fill=white] (0,0) circle (0.2) ;
\draw [thick,>=stealth,->] (-0.5,-0.5) -- ++(1,1) ;
\end{scope}

\begin{scope}[xshift=0cm,yshift=0cm]
\draw [thick,>=stealth,->] (-0.5,-0.5) -- ++(1,1) ;
\draw [white,fill=white] (0,0) circle (0.2) ;
\draw [thick,>=stealth,->] (0.5,-0.5) -- ++(-1,1) ;
\end{scope}

\begin{scope}[xshift=5cm,yshift=0cm]
\draw [thick,>=stealth,->] (-0.5,-0.5) -- ++(1,1) ;
\draw [thick,>=stealth,->] (0.5,-0.5) -- ++(-1,1) ;
\draw (0,0) circle (0.2) ;
\end{scope}

\begin{scope}[xshift=-5cm,yshift=-3cm]
\begin{scope}[rotate=45]
\clip (-0.3,-1.2) -- (-0.3,-1.5) -- (0.3,-1.5) -- (0.3,-1.2) -- cycle ;
\draw (0,-1.2) ellipse (0.3 and 0.12) ;
\end{scope}
\begin{scope}[rotate=45]
\clip (-0.3,-1.2) -- (-0.3,-1) -- (0.3,-1) -- (0.3,-1.2) -- cycle ;
\draw [dashed] (0,-1.2) ellipse (0.3 and 0.12) ;
\end{scope}
\begin{scope}[rotate=-45]
\draw (0,1.2) ellipse (0.7 and 0.28) ;
\draw (-0.7,-1.2) -- (-0.7,-0.3) ;
\draw (-0.7,0.3) -- (-0.7,1.2) ;
\draw (0.7,-1.2) -- (0.7,1.2) ;
\end{scope}
\begin{scope}[rotate=-45]
\clip (-0.7,-1.2) -- (-0.7,-1.5) -- (0.7,-1.5) -- (0.7,-1.2) -- cycle ;
\draw (0,-1.2) ellipse (0.7 and 0.28) ;
\end{scope}
\begin{scope}[rotate=-45]
\clip (-0.7,-1.2) -- (-0.7,-0.9) -- (0.7,-0.9) -- (0.7,-1.2) -- cycle ;
\draw [dashed] (0,-1.2) ellipse (0.7 and 0.28) ;
\end{scope}
\begin{scope}[rotate=45]
\draw [dashed] (0,-0.35) ellipse (0.5 and 0.2) ;
\fill [white] (-0.3,-0.35) -- (-0.3,0) -- (0.3,0) -- (0.3,-0.35) -- cycle ;
\draw [dashed] (0,-0.35) ellipse (0.2 and 0.08) ;
\draw [dashed] (-0.3,-0.5) -- (-0.3,0.3) ;
\draw [dashed] (0.3,-0.5) -- (0.3,0.3) ;
\draw [dashed] (0,0.3) ellipse (0.3 and 0.12) ;
\fill [white] (-0.17,0.3) -- (-0.17,0.5) -- (0.17,0.5) -- (0.17,0.3) -- cycle ;
\draw (-0.3,-1.2) -- (-0.3,-0.7) ;
\draw (0.3,-1.2) -- (0.3,-0.7) ;
\draw (-0.3,0.5) -- (-0.3,1.2) ;
\draw (0.3,0.5) -- (0.3,1.2) ;
\draw (0,1.2) ellipse (0.3 and 0.12) ;
\end{scope}
\begin{scope}[rotate=45]
\clip (-0.3,0.5) -- (-0.3,0.3) -- (0.3,0.3) -- (0.3,0.5) -- cycle ;
\draw (0,0.5) ellipse (0.3 and 0.12) ;
\end{scope}
\begin{scope}[rotate=45]
\clip (-0.3,0.5) -- (-0.3,0.7) -- (0.3,0.7) -- (0.3,0.5) -- cycle ;
\draw [dashed] (0,0.5) ellipse (0.3 and 0.12) ;
\end{scope}
\end{scope}

\begin{scope}[xshift=0cm,yshift=-3cm]
\begin{scope}[rotate=-45]
\draw (0,-0.75) ellipse (0.3 and 0.12) ;
\fill [white] (-0.3,-0.7) -- (-0.3,-0.5) -- (0.3,-0.5) -- (0.3,-0.7) -- cycle ;
\draw (-0.3,-1.2) -- (-0.3,-0.75) ;
\draw (0.3,-1.2) -- (0.3,-0.75) ;
\end{scope}
\begin{scope}[rotate=-45]
\clip (-0.3,-1.2) -- (-0.3,-1.5) -- (0.3,-1.5) -- (0.3,-1.2) -- cycle ;
\draw (0,-1.2) ellipse (0.3 and 0.12) ;
\end{scope}
\begin{scope}[rotate=-45]
\clip (-0.3,-1.2) -- (-0.3,-1) -- (0.3,-1) -- (0.3,-1.2) -- cycle ;
\draw [dashed] (0,-1.2) ellipse (0.3 and 0.12) ;
\end{scope}
\begin{scope}[rotate=45]
\draw (0,1.2) ellipse (0.7 and 0.28) ;
\draw (-0.7,-1.2) -- (-0.7,1.2) ;
\draw (0.7,-1.2) -- (0.7,1.2) ;
\end{scope}
\begin{scope}[rotate=45]
\clip (-0.7,-1.2) -- (-0.7,-1.5) -- (0.7,-1.5) -- (0.7,-1.2) -- cycle ;
\draw (0,-1.2) ellipse (0.7 and 0.28) ;
\end{scope}
\begin{scope}[rotate=45]
\clip (-0.7,-1.2) -- (-0.7,-0.9) -- (0.7,-0.9) -- (0.7,-1.2) -- cycle ;
\draw [dashed] (0,-1.2) ellipse (0.7 and 0.28) ;
\end{scope}
\begin{scope}[rotate=-45]
\draw (0,0.35) ellipse (0.5 and 0.2) ;
\fill [white] (-0.3,0.35) -- (-0.3,0.8) -- (0.3,0.8) -- (0.3,0.35) -- cycle ;
\draw (0,1.2) ellipse (0.3 and 0.12) ;
\draw [dashed] (-0.3,-0.4) -- (-0.3,0.2) ;
\draw (-0.3,0.2) -- (-0.3,1.2) ;
\draw [dashed] (0.3,-0.4) -- (0.3,0.2) ;
\draw [dashed] (0,-0.4) ellipse (0.3 and 0.12) ;
\draw (0.3,0.2) -- (0.3,1.2) ;
\draw [dashed] (0,0.35) ellipse (0.2 and 0.08) ;
\end{scope}
\end{scope}

\begin{scope}[xshift=5cm,yshift=-3cm]
\begin{scope}[rotate=-45]
\draw (0,1) ellipse (0.5 and 0.2) ;
\draw (-0.5,-1) -- (-0.5,1) ;
\draw (0.5,-1) -- (0.5,1) ;
\end{scope}
\begin{scope}[rotate=-45]
\clip (-0.5,-1) -- (-0.5,-1.5) -- (0.5,-1.5) -- (0.5,-1) -- cycle ;
\draw (0,-1) ellipse (0.5 and 0.2) ;
\end{scope}
\begin{scope}[rotate=-45]
\clip (-0.5,-1) -- (-0.5,-0.5) -- (0.5,-0.5) -- (0.5,-1) -- cycle ;
\draw [dashed] (0,-1) ellipse (0.5 and 0.2) ;
\end{scope}
\begin{scope}[rotate=45]
\draw (0,1) ellipse (0.5 and 0.2) ;
\draw (-0.5,-1) -- (-0.5,-0.5) ;
\draw (-0.5,0.5) -- (-0.5,1) ;
\draw (0.5,-1) -- (0.5,-0.5) ;
\draw (0.5,0.5) -- (0.5,1) ;
\end{scope}
\begin{scope}[rotate=45]
\clip (-0.5,-1) -- (-0.5,-1.5) -- (0.5,-1.5) -- (0.5,-1) -- cycle ;
\draw (0,-1) ellipse (0.5 and 0.2) ;
\end{scope}
\begin{scope}[rotate=45]
\clip (-0.5,-1) -- (-0.5,-0.5) -- (0.5,-0.5) -- (0.5,-1) -- cycle ;
\draw [dashed] (0,-1) ellipse (0.5 and 0.2) ;
\end{scope}
\end{scope}
\end{tikzpicture}
\caption{Image under Tube of each crossing}
\end{figure}

We now define the Spun map, which gives another way to obtain a knotted surface from a string link. Note however that this map is only defined on classical string links, while the Tube map is defined on welded objects.

\begin{de}\label{defSpun}
Let $L\in\SL$ be given by a parametrization $(x_i(t),y_i(t),z_i(t))\in B^{2,1}$ for $1\leq i\leq n$ and $t\in[0,1]$. Then Spun($L$) is defined to be the string 2--link parametrized by :
$$\left(\frac{x_i(t)}{2},\frac{y_i(t)-1}{2}\cos(\theta),\frac{y_i(t)-1}{2}\sin(\theta),z_i(t)\right)\in B^{3,1},$$
for $1\leq i\leq n$, $t\in[0,1]$, and $0\leq\theta\leq 2\pi$. In other words, we start by placing $L$ in $B^2\big((0,-\frac{1}{2}),\frac{1}{2}\big)\times\{0\}\times[0,1]\subset B^{3,1}\subset\RR^4$ by applying the map $(x,y,z)\mapsto(x/2,(y-1)/2,0,z)$, then we take its trace under a complete rotation around the plane $\RR\times\{0\}^2\times\RR$. This defines a map $\text{Spun}:\SL\to\dl{2}$. Indeed, an isotopy of $B^{2,1}$ induces an isotopy of $B^{3,1}$ by the same construction.
\end{de}

Suppose $L\in\SL$ is represented by a diagram $D\in SLD_n$, which is obtained by projecting $L$ onto the $(yOz)$ plane. Then Spun($L$) is represented by the broken surface diagram obtained by rotating $D$, parametrized by:
$$\left(\frac{y_i(t)-1}{2}\cos(\theta),\frac{y_i(t)-1}{2}\sin(\theta),z_i(t)\right)\in B^{2,1},$$
for $1\leq i\leq n$, $t\in[0,1]$, and $0\leq\theta\leq 2\pi$, with the upper/lower information being given by the value of $x_i(t)$. Note that the only singularities of this broken surface diagram are lines of double points, obtained by rotating the crossings of $D$.

\subsection{Local moves}

We will now consider local moves in a general way.

\subsubsection{Definition}

\begin{de}\label{localmove}
A \emph{classical local move} is a local move which only involves classical crossings. A \emph{virtual local move} can involve any type of crossings.
\end{de}

In particular, classical local moves are virtual local moves (in the same way as classical diagrams are also virtual). The latter will also be simply referred to as local moves, with the mention virtual being added when necessary to distinguish between classical and virtual moves.

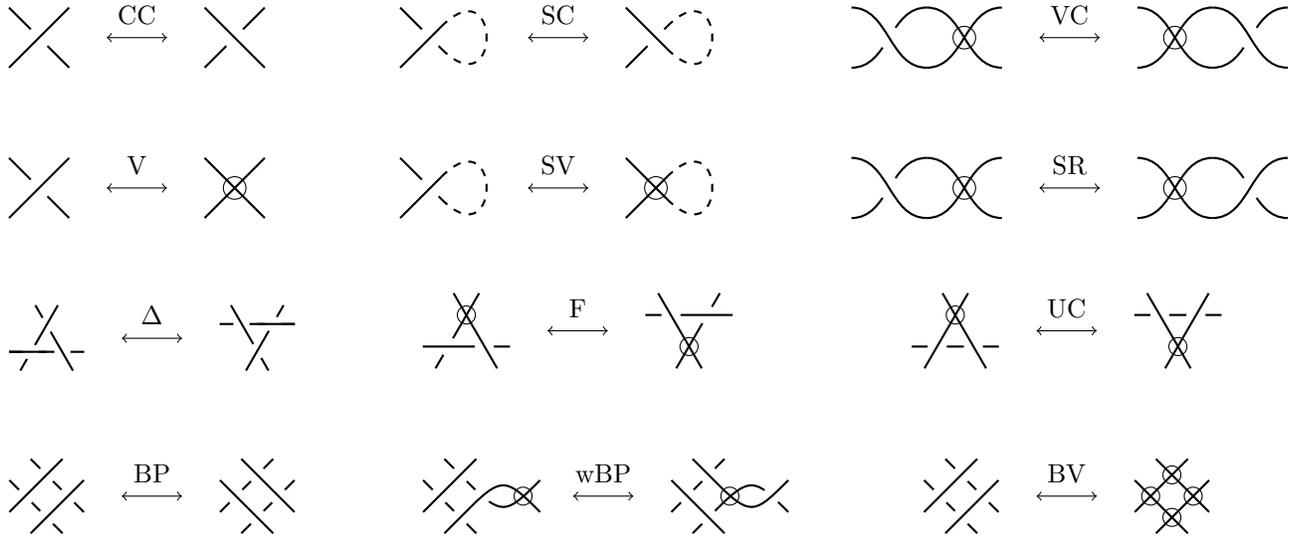
\begin{figure}[h]
\centering
\begin{tikzpicture}
\draw [white] (0,0) -- (17,0) -- (17,-7.2) -- (0,-7.2) -- cycle ;

\begin{scope}[xshift=0cm,yshift=-1cm]
\draw [thick] (0,0.8) -- (0.8,0) ;
\draw [white,fill=white] (0.4,0.4) circle (0.15) ;
\draw [thick] (0,0) -- (0.8,0.8) ;
\draw [<->] (1.3,0.4) -- (2.1,0.4) ;
\draw (1.7,0.7) node{CC} ;
\draw [thick] (2.6,0) -- (3.4,0.8) ;
\draw [white,fill=white] (3,0.4) circle (0.15) ;
\draw [thick] (2.6,0.8) -- (3.4,0) ;
\end{scope}

\begin{scope}[xshift=5.2cm,yshift=-1cm]
\draw [thick] (0,0.8) -- (0.55,0.25) ;
\draw [thick,dashed] (0.55,0.25) .. controls +(0.8,-0.8) and +(0.8,0.8) .. (0.55,0.55) ;
\draw [white,fill=white] (0.4,0.4) circle (0.15) ;
\draw [thick] (0,0) -- (0.55,0.55) ;
\draw [<->] (1.7,0.4) -- (2.5,0.4) ;
\draw (2.1,0.7) node{SC} ;
\draw [thick] (3,0) -- (3.55,0.55) ;
\draw [thick,dashed] (3.55,0.25) .. controls +(0.8,-0.8) and +(0.8,0.8) .. (3.55,0.55) ;
\draw [white,fill=white] (3.4,0.4) circle (0.15) ;
\draw [thick] (3,0.8) -- (3.55,0.25) ;
\end{scope}

\begin{scope}[xshift=11.2cm,yshift=-1cm]
\draw [thick] (0,0) .. controls +(0.5,0) and +(-0.5,0) .. (1,0.8) ;
\draw [thick] (1,0.8) .. controls +(0.5,0) and +(-0.5,0) .. (2,0) ;
\draw [white,fill=white] (0.5,0.4) circle (0.15) ;
\draw [thick] (0,0.8) .. controls +(0.5,0) and +(-0.5,0) .. (1,0) ;
\draw [thick] (1,0) .. controls +(0.5,0) and +(-0.5,0) .. (2,0.8) ;
\draw (1.5,0.4) circle (0.15) ;
\draw [<->] (2.5,0.4) -- (3.3,0.4) ;
\draw (2.9,0.7) node{VC} ;
\draw [thick] (3.8,0.8) .. controls +(0.5,0) and +(-0.5,0) .. (4.8,0) ;
\draw [thick] (4.8,0) .. controls +(0.5,0) and +(-0.5,0) .. (5.8,0.8) ;
\draw [white,fill=white] (5.3,0.4) circle (0.15) ;
\draw [thick] (3.8,0) .. controls +(0.5,0) and +(-0.5,0) .. (4.8,0.8) ;
\draw [thick] (4.8,0.8) .. controls +(0.5,0) and +(-0.5,0) .. (5.8,0) ;
\draw (4.3,0.4) circle (0.15) ;
\end{scope}

\begin{scope}[xshift=0cm,yshift=-3cm]
\draw [thick] (0,0.8) -- (0.8,0) ;
\draw [white,fill=white] (0.4,0.4) circle (0.15) ;
\draw [thick] (0,0) -- (0.8,0.8) ;
\draw [<->] (1.3,0.4) -- (2.1,0.4) ;
\draw (1.7,0.7) node{V} ;
\draw [thick] (2.6,0) -- (3.4,0.8) ;
\draw (3,0.4) circle (0.15) ;
\draw [thick] (2.6,0.8) -- (3.4,0) ;
\end{scope}

\begin{scope}[xshift=5.2cm,yshift=-3cm]
\draw [thick] (0,0.8) -- (0.55,0.25) ;
\draw [thick,dashed] (0.55,0.25) .. controls +(0.8,-0.8) and +(0.8,0.8) .. (0.55,0.55) ;
\draw [white,fill=white] (0.4,0.4) circle (0.15) ;
\draw [thick] (0,0) -- (0.55,0.55) ;
\draw [<->] (1.7,0.4) -- (2.5,0.4) ;
\draw (2.1,0.7) node{SV} ;
\draw [thick] (3,0) -- (3.55,0.55) ;
\draw [thick,dashed] (3.55,0.25) .. controls +(0.8,-0.8) and +(0.8,0.8) .. (3.55,0.55) ;
\draw (3.4,0.4) circle (0.15) ;
\draw [thick] (3,0.8) -- (3.55,0.25) ;
\end{scope}

\begin{scope}[xshift=11.2cm,yshift=-3cm]
\draw [thick] (0,0) .. controls +(0.5,0) and +(-0.5,0) .. (1,0.8) ;
\draw [thick] (1,0.8) .. controls +(0.5,0) and +(-0.5,0) .. (2,0) ;
\draw [white,fill=white] (0.5,0.4) circle (0.15) ;
\draw [thick] (0,0.8) .. controls +(0.5,0) and +(-0.5,0) .. (1,0) ;
\draw [thick] (1,0) .. controls +(0.5,0) and +(-0.5,0) .. (2,0.8) ;
\draw (1.5,0.4) circle (0.15) ;
\draw [<->] (2.5,0.4) -- (3.3,0.4) ;
\draw (2.9,0.7) node{SR} ;
\draw [thick] (3.8,0) .. controls +(0.5,0) and +(-0.5,0) .. (4.8,0.8) ;
\draw [thick] (4.8,0.8) .. controls +(0.5,0) and +(-0.5,0) .. (5.8,0) ;
\draw [white,fill=white] (5.3,0.4) circle (0.15) ;
\draw [thick] (3.8,0.8) .. controls +(0.5,0) and +(-0.5,0) .. (4.8,0) ;
\draw [thick] (4.8,0) .. controls +(0.5,0) and +(-0.5,0) .. (5.8,0.8) ;
\draw (4.3,0.4) circle (0.15) ;
\end{scope}

\begin{scope}[xshift=0cm,yshift=-5cm]
\draw [thick] (0,0.22) -- (1,0.22) ;
\draw [white,fill=white] (0.71,0.22) circle (0.1) ;
\draw [thick] (0.85,-0.03) -- ++(-0.5,0.87) ;
\draw [white,fill=white] (0.5,0.58) circle (0.1) ;
\draw [thick] (0.15,-0.03) -- ++(0.5,0.87) ;
\draw [white,fill=white] (0.29,0.22) circle (0.1) ;
\draw [thick] (0,0.22) -- (0.5,0.22) ;
\draw [<->] (1.5,0.4) -- (2.3,0.4) ;
\draw (1.9,0.7) node{$\Delta$} ;
\draw [thick] (2.8,0.59) -- (3.8,0.59) ;
\draw [white,fill=white] (3.09,0.59) circle (0.1) ;
\draw [thick] (2.95,0.84) -- ++(0.5,-0.87) ;
\draw [white,fill=white] (3.3,0.23) circle (0.1) ;
\draw [thick] (3.65,0.84) -- ++(-0.5,-0.87) ;
\draw [white,fill=white] (3.51,0.59) circle (0.1) ;
\draw [thick] (3.3,0.59) -- (3.8,0.59) ;
\end{scope}

\begin{scope}[xshift=5.5cm,yshift=-5cm]
\draw [thick] (0.17,0) -- ++(0.58,1) ;
\draw (0.58,0.71) circle (0.12) ;
\draw [white,fill=white] (0.34,0.29) circle (0.12) ;
\draw [thick] (0,0.29) -- (1.16,0.29) ;
\draw [white,fill=white] (0.82,0.29) circle (0.12) ;
\draw [thick] (0.99,0) -- ++(-0.58,1) ;
\draw [<->] (1.66,0.5) -- (2.46,0.5) ;
\draw (2.06,0.8) node{F} ;
\draw [thick] (3.95,1) -- ++(-0.58,-1) ;
\draw (3.54,0.29) circle (0.12) ;
\draw [white,fill=white] (3.78,0.71) circle (0.12) ;
\draw [thick] (2.96,0.71) -- (4.12,0.71) ;
\draw [white,fill=white] (3.3,0.71) circle (0.12) ;
\draw [thick] (3.13,1) -- ++(0.58,-1) ;
\end{scope}

\begin{scope}[xshift=12cm,yshift=-5cm]
\draw [thick] (0,0.29) -- (1.16,0.29) ;
\draw [white,fill=white] (0.34,0.29) circle (0.12) ;
\draw [white,fill=white] (0.82,0.29) circle (0.12) ;
\draw [thick] (0.99,0) -- ++(-0.58,1);
\draw (0.58,0.71) circle (0.12) ;
\draw [thick] (0.17,0) -- ++(0.58,1);
\draw [<->] (1.66,0.5) -- (2.46,0.5) ;
\draw (2.06,0.8) node{UC} ;
\draw [thick] (2.96,0.71) -- (4.12,0.71) ;
\draw [white,fill=white] (3.78,0.71) circle (0.12) ;
\draw [white,fill=white] (3.3,0.71) circle (0.12) ;
\draw [thick] (3.13,1) -- ++(0.58,-1);
\draw (3.54,0.29) circle (0.12) ;
\draw [thick] (3.95,1) -- ++(-0.58,-1);
\end{scope}

\begin{scope}[xshift=0.5cm,yshift=-6.7cm]
\begin{scope}[rotate=45]
\draw [thick] (-0.2,-0.5) -- (-0.2,0.5) ;
\draw [thick] (0.2,-0.5) -- (0.2,0.5) ;
\fill [white] (-0.2,-0.2) circle (0.12) ;
\fill [white] (-0.2,0.2) circle (0.12) ;
\fill [white] (0.2,-0.2) circle (0.12) ;
\fill [white] (0.2,0.2) circle (0.12) ;
\draw [thick] (-0.5,-0.2) -- (0.5,-0.2) ;
\draw [thick] (-0.5,0.2) -- (0.5,0.2) ;
\end{scope}
\begin{scope}[xshift=1cm]
\draw [<->] (0,0) -- (0.8,0) ;
\draw (0.4,0.3) node{BP} ;
\end{scope}
\begin{scope}[xshift=2.8cm,rotate=45]
\draw [thick] (-0.5,-0.2) -- (0.5,-0.2) ;
\draw [thick] (-0.5,0.2) -- (0.5,0.2) ;
\fill [white] (-0.2,-0.2) circle (0.12) ;
\fill [white] (-0.2,0.2) circle (0.12) ;
\fill [white] (0.2,-0.2) circle (0.12) ;
\fill [white] (0.2,0.2) circle (0.12) ;
\draw [thick] (-0.2,-0.5) -- (-0.2,0.5) ;
\draw [thick] (0.2,-0.5) -- (0.2,0.5) ;
\end{scope}
\end{scope}

\begin{scope}[xshift=6cm,yshift=-6.7cm]
\begin{scope}[rotate=45]
\draw [thick] (-0.2,-0.5) -- (-0.2,0.5) ;
\draw [thick] (0.2,0) -- (0.2,0.5) ;
\draw [thick] (0.2,0) .. controls +(0,-0.6) and +(-0.6,0) .. (0.9,-0.6) ;
\fill [white] (-0.2,-0.2) circle (0.12) ;
\fill [white] (-0.2,0.2) circle (0.12) ;
\fill [white] (0.2,-0.2) circle (0.12) ;
\fill [white] (0.2,0.2) circle (0.12) ;
\draw (0.59,-0.59) circle (0.12) ;
\draw [thick] (-0.5,-0.2) -- (0,-0.2) ;
\draw [thick] (0,-0.2) .. controls +(0.6,0) and +(0,0.6) .. (0.6,-0.9) ;
\draw [thick] (-0.5,0.2) -- (0.5,0.2) ;
\end{scope}
\begin{scope}[xshift=1.5cm]
\draw [<->] (0,0) -- (0.8,0) ;
\draw (0.4,0.3) node{wBP} ;
\end{scope}
\begin{scope}[xshift=3.3cm,rotate=45]
\draw [thick] (-0.5,-0.2) -- (0,-0.2) ;
\draw [thick] (0,-0.2) .. controls +(0.6,0) and +(0,0.6) .. (0.6,-0.9) ;
\draw [thick] (-0.5,0.2) -- (0.5,0.2) ;
\fill [white] (-0.2,-0.2) circle (0.12) ;
\fill [white] (-0.2,0.2) circle (0.12) ;
\draw (0.2,-0.2) circle (0.12) ;
\fill [white] (0.2,0.2) circle (0.12) ;
\fill [white] (0.59,-0.59) circle (0.12) ;
\draw [thick] (-0.2,-0.5) -- (-0.2,0.5) ;
\draw [thick] (0.2,0) -- (0.2,0.5) ;
\draw [thick] (0.2,0) .. controls +(0,-0.6) and +(-0.6,0) .. (0.9,-0.6) ;
\end{scope}
\end{scope}

\begin{scope}[xshift=12.66cm,yshift=-6.7cm]
\begin{scope}[xshift=0cm,rotate=45]
\draw [thick] (-0.2,-0.5) -- (-0.2,0.5) ;
\draw [thick] (0.2,-0.5) -- (0.2,0.5) ;
\fill [white] (-0.2,-0.2) circle (0.12) ;
\fill [white] (-0.2,0.2) circle (0.12) ;
\fill [white] (0.2,-0.2) circle (0.12) ;
\fill [white] (0.2,0.2) circle (0.12) ;
\draw [thick] (-0.5,-0.2) -- (0.5,-0.2) ;
\draw [thick] (-0.5,0.2) -- (0.5,0.2) ;
\end{scope}
\begin{scope}[xshift=1cm]
\draw [<->] (0,0) -- (0.8,0) ;
\draw (0.4,0.3) node{BV} ;
\end{scope}
\begin{scope}[xshift=2.8cm,rotate=45]
\draw [thick] (-0.5,-0.2) -- (0.5,-0.2) ;
\draw [thick] (-0.5,0.2) -- (0.5,0.2) ;
\draw (-0.2,-0.2) circle (0.12) ;
\draw (-0.2,0.2) circle (0.12) ;
\draw (0.2,-0.2) circle (0.12) ;
\draw (0.2,0.2) circle (0.12) ;
\draw [thick] (-0.2,-0.5) -- (-0.2,0.5) ;
\draw [thick] (0.2,-0.5) -- (0.2,0.5) ;
\end{scope}
\end{scope}
\end{tikzpicture}
\caption{Some local moves on virtual string link diagrams}
\end{figure}

A few examples are illustrated on Figure 8, where a dotted line indicates that the crossing occurs between two portions of the same strand. Some classical local moves are derived from topological operations on links. For example, the $CC$ (for \textquotedblleft crossing change") move corresponds to homotopy, which allows strands to cross each other, while $SC$ (for \textquotedblleft self-crossing change") corresponds to link homotopy, only allowing each strand to cross itself. The $BP$ (for \textquotedblleft band-pass") move represents the crossing of two "bands", delimited by parallel strands. \\

In an effort to extend the effect of some classical local moves to welded string links, certain local moves on virtual diagrams, which closely resemble their classical counterparts, are introduced. For example, $V$ (for \textquotedblleft virtualization") and $SV$ (for \textquotedblleft self-virtualization") are derived from $CC$ and $SC$ respectively, while $F$ (for \textquotedblleft fused") is derived from $\Delta$, and $BV$ (for \textquotedblleft band virtualization") is derived from $BP$. This notion of extension will be discussed in the next section. \\

Finally, some virtual local moves come naturally from Gauss diagrams. For example, $VC$ (for \textquotedblleft virtual conjugation") reverses the orientation of the arrow representing the classical crossing, while $SR$ (for \textquotedblleft sign reversal") changes its sign. \\

We also give their version on Gauss diagrams in Figure 9. The dotted lines in $SC$ and $SV$ indicate that the extremities of the arrows belong to the same strand, but are not necessarily adjacent on this strand. These are called self-arrows. As before, the $(\ast)$ indicates the presence of some sign conditions: the $\Delta$ move must verify the same conditions as the $R3$ move, while the $BP$, $wBP$ and $BV$ moves must verify $\eps_{ij}\eps_{kl}=\delta_i\delta_j\delta_k\delta_l$, where the $\delta$'s are defined as in Section 1.2.1.

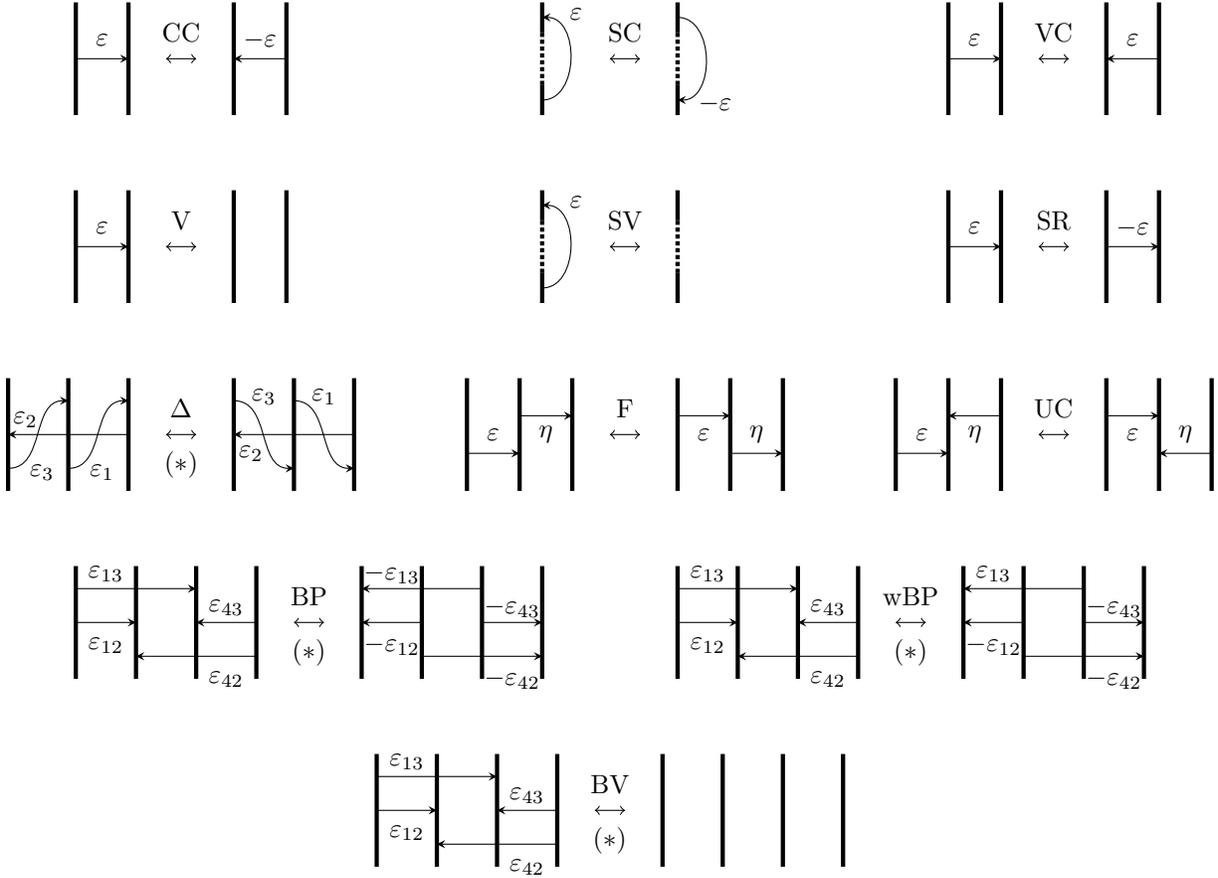
\begin{figure}[h]
\centering
\begin{tikzpicture}
\draw [white] (0,0.1) -- (16,0.1) -- (16,-11.7) -- (0,-11.7) -- cycle ;

\begin{scope}[xshift=0.9cm,yshift=-1.5cm]
\draw [ultra thick] (0,0) -- ++(0,1.5) ;
\draw [ultra thick] (0.7,0) -- ++(0,1.5) ;
\draw [>=stealth,->] (0,0.75) -- ++(0.7,0) ;
\draw (0.35,1) node{$\eps$} ;
\draw [<->] (1.2,0.75) -- (1.6,0.75) ;
\draw (1.4,1.1) node{CC} ;
\draw [ultra thick] (2.1,0) -- ++(0,1.5) ;
\draw [ultra thick] (2.8,0) -- ++(0,1.5) ;
\draw [>=stealth,->] (2.8,0.75) -- ++(-0.7,0) ;
\draw (2.45,1) node{$-\eps$} ;
\end{scope}

\begin{scope}[xshift=7.1cm,yshift=-1.5cm]
\draw [ultra thick] (0,0) -- ++(0,0.4) ;
\draw [ultra thick,densely dotted] (0,0.4) -- ++(0,0.7) ;
\draw [ultra thick] (0,1.1) -- ++(0,0.4) ;
\draw [>=stealth,->] (0,0.2) .. controls +(0.5,0) and +(0.5,0) .. (0,1.3) ;
\draw (0.45,1.35) node{$\eps$} ;
\draw [<->] (0.9,0.75) -- (1.3,0.75) ;
\draw (1.1,1.1) node{SC} ;
\draw [ultra thick] (1.8,0) -- ++(0,0.4) ;
\draw [ultra thick,densely dotted] (1.8,0.4) -- ++(0,0.7) ;
\draw [ultra thick] (1.8,1.1) -- ++(0,0.4) ;
\draw [>=stealth,->] (1.8,1.3) .. controls +(0.5,0) and +(0.5,0) .. (1.8,0.2) ;
\draw (2.3,0.15) node{$-\eps$} ;
\end{scope}

\begin{scope}[xshift=12.5cm,yshift=-1.5cm]
\draw [ultra thick] (0,0) -- ++(0,1.5) ;
\draw [ultra thick] (0.7,0) -- ++(0,1.5) ;
\draw [>=stealth,->] (0,0.75) -- ++(0.7,0) ;
\draw (0.35,1) node{$\eps$} ;
\draw [<->] (1.2,0.75) -- (1.6,0.75) ;
\draw (1.4,1.1) node{VC} ;
\draw [ultra thick] (2.1,0) -- ++(0,1.5) ;
\draw [ultra thick] (2.8,0) -- ++(0,1.5) ;
\draw [>=stealth,->] (2.8,0.75) -- ++(-0.7,0) ;
\draw (2.45,1) node{$\eps$} ;
\end{scope}

\begin{scope}[xshift=0.9cm,yshift=-4cm]
\draw [ultra thick] (0,0) -- ++(0,1.5) ;
\draw [ultra thick] (0.7,0) -- ++(0,1.5) ;
\draw [>=stealth,->] (0,0.75) -- ++(0.7,0) ;
\draw (0.35,1) node{$\eps$} ;
\draw [<->] (1.2,0.75) -- (1.6,0.75) ;
\draw (1.4,1.1) node{V} ;
\draw [ultra thick] (2.1,0) -- ++(0,1.5) ;
\draw [ultra thick] (2.8,0) -- ++(0,1.5) ;
\end{scope}

\begin{scope}[xshift=7.1cm,yshift=-4cm]
\draw [ultra thick] (0,0) -- ++(0,0.4) ;
\draw [ultra thick,densely dotted] (0,0.4) -- ++(0,0.7) ;
\draw [ultra thick] (0,1.1) -- ++(0,0.4) ;
\draw [>=stealth,->] (0,0.2) .. controls +(0.5,0) and +(0.5,0) .. (0,1.3) ;
\draw (0.45,1.35) node{$\eps$} ;
\draw [<->] (0.9,0.75) -- (1.3,0.75) ;
\draw (1.1,1.1) node{SV} ;
\draw [ultra thick] (1.8,0) -- ++(0,0.4) ;
\draw [ultra thick,densely dotted] (1.8,0.4) -- ++(0,0.7) ;
\draw [ultra thick] (1.8,1.1) -- ++(0,0.4) ;
\end{scope}

\begin{scope}[xshift=12.5cm,yshift=-4cm]
\draw [ultra thick] (0,0) -- ++(0,1.5) ;
\draw [ultra thick] (0.7,0) -- ++(0,1.5) ;
\draw [>=stealth,->] (0,0.75) -- ++(0.7,0) ;
\draw (0.35,1) node{$\eps$} ;
\draw [<->] (1.2,0.75) -- (1.6,0.75) ;
\draw (1.4,1.1) node{SR} ;
\draw [ultra thick] (2.1,0) -- ++(0,1.5) ;
\draw [ultra thick] (2.8,0) -- ++(0,1.5) ;
\draw [>=stealth,->] (2.1,0.75) -- ++(0.7,0) ;
\draw (2.45,1) node{$-\eps$} ;
\end{scope}

\begin{scope}[xshift=0cm,yshift=-6.5cm]
\draw [ultra thick] (0,0) -- ++(0,1.5) ;
\draw [ultra thick] (0.8,0) -- ++(0,1.5) ;
\draw [ultra thick] (1.6,0) -- ++(0,1.5) ;
\draw [>=stealth,->] (0,0.3) .. controls +(0.5,0) and +(-0.5,0) .. (0.8,1.2) ;
\draw (0.45,0.25) node{$\eps_3$} ;
\draw [>=stealth,->] (0.8,0.3) .. controls +(0.5,0) and +(-0.5,0) .. (1.6,1.2) ;
\draw (1.25,0.25) node{$\eps_1$} ;
\draw [>=stealth,->] (1.6,0.75) -- ++(-1.6,0) ;
\draw (0.23,0.95) node{$\eps_2$} ;
\draw [<->] (2.1,0.75) -- (2.5,0.75) ;
\draw (2.3,1.1) node{$\Delta$} ;
\draw (2.3,0.35) node{$(\ast)$} ;
\draw [ultra thick] (3,0) -- ++(0,1.5) ;
\draw [ultra thick] (3.8,0) -- ++(0,1.5) ;
\draw [ultra thick] (4.6,0) -- ++(0,1.5) ;
\draw [>=stealth,->] (3,1.2) .. controls +(0.5,0) and +(-0.5,0) .. (3.8,0.3) ;
\draw (3.4,1.25) node{$\eps_3$} ;
\draw [>=stealth,->] (3.8,1.2) .. controls +(0.5,0) and +(-0.5,0) .. (4.6,0.3) ;
\draw (4.2,1.25) node{$\eps_1$} ;
\draw [>=stealth,->] (4.6,0.75) -- ++(-1.6,0) ;
\draw (3.23,0.5) node{$\eps_2$} ;
\end{scope}

\begin{scope}[xshift=6.1cm,yshift=-6.5cm]
\draw [ultra thick] (0,0) -- ++(0,1.5) ;
\draw [ultra thick] (0.7,0) -- ++(0,1.5) ;
\draw [ultra thick] (1.4,0) -- ++(0,1.5) ;
\draw [>=stealth,->] (0,0.5) -- ++(0.7,0) ;
\draw (0.35,0.75) node{$\eps$} ;
\draw [>=stealth,->] (0.7,1) -- ++(0.7,0) ;
\draw (1.05,0.75) node{$\eta$} ;
\draw [<->] (1.9,0.75) -- (2.3,0.75) ;
\draw (2.1,1.1) node{F} ;
\draw [ultra thick] (2.8,0) -- ++(0,1.5) ;
\draw [ultra thick] (3.5,0) -- ++(0,1.5) ;
\draw [ultra thick] (4.2,0) -- ++(0,1.5) ;
\draw [>=stealth,->] (2.8,1) -- ++(0.7,0) ;
\draw (3.15,0.75) node{$\eps$} ;
\draw [>=stealth,->] (3.5,0.5) -- ++(0.7,0) ;
\draw (3.85,0.75) node{$\eta$} ;
\end{scope}

\begin{scope}[xshift=11.8cm,yshift=-6.5cm]
\draw [ultra thick] (0,0) -- ++(0,1.5) ;
\draw [ultra thick] (0.7,0) -- ++(0,1.5) ;
\draw [ultra thick] (1.4,0) -- ++(0,1.5) ;
\draw [>=stealth,->] (0,0.5) -- ++(0.7,0) ;
\draw (0.35,0.75) node{$\eps$} ;
\draw [>=stealth,->] (1.4,1) -- ++(-0.7,0) ;
\draw (1.05,0.75) node{$\eta$} ;
\draw [<->] (1.9,0.75) -- (2.3,0.75) ;
\draw (2.1,1.1) node{UC} ;
\draw [ultra thick] (2.8,0) -- ++(0,1.5) ;
\draw [ultra thick] (3.5,0) -- ++(0,1.5) ;
\draw [ultra thick] (4.2,0) -- ++(0,1.5) ;
\draw [>=stealth,->] (2.8,1) -- ++(0.7,0) ;
\draw (3.15,0.75) node{$\eps$} ;
\draw [>=stealth,->] (4.2,0.5) -- ++(-0.7,0) ;
\draw (3.85,0.75) node{$\eta$} ;
\end{scope}

\begin{scope}[xshift=0.9cm,yshift=-9cm]
\draw [ultra thick] (0,0) -- ++(0,1.5) ;
\draw [ultra thick] (0.8,0) -- ++(0,1.5) ;
\draw [ultra thick] (1.6,0) -- ++(0,1.5) ;
\draw [ultra thick] (2.4,0) -- ++(0,1.5) ;
\draw [>=stealth,->] (0,1.2) -- (1.6,1.2) ;
\draw (0.4,1.4) node{$\eps_{13}$} ;
\draw [>=stealth,->] (0,0.75) -- (0.8,0.75) ;
\draw (0.4,0.45) node{$\eps_{12}$} ;
\draw [>=stealth,->] (2.4,0.75) -- (1.6,0.75) ;
\draw (2,0.95) node{$\eps_{43}$} ;
\draw [>=stealth,->] (2.4,0.3) -- (0.8,0.3) ;
\draw (2,0) node{$\eps_{42}$} ;
\draw [<->] (2.9,0.75) -- (3.3,0.75) ;
\draw (3.1,1.1) node{BP} ;
\draw (3.1,0.35) node{$(\ast)$} ;
\draw [ultra thick] (3.8,0) -- ++(0,1.5) ;
\draw [ultra thick] (4.6,0) -- ++(0,1.5) ;
\draw [ultra thick] (5.4,0) -- ++(0,1.5) ;
\draw [ultra thick] (6.2,0) -- ++(0,1.5) ;
\draw [>=stealth,->] (5.4,1.2) -- (3.8,1.2) ;
\draw (4.2,1.4) node{$-\eps_{13}$} ;
\draw [>=stealth,->] (4.6,0.75) -- (3.8,0.75) ;
\draw (4.2,0.45) node{$-\eps_{12}$} ;
\draw [>=stealth,->] (5.4,0.75) -- (6.2,0.75) ;
\draw (5.8,0.95) node{$-\eps_{43}$} ;
\draw [>=stealth,->] (4.6,0.3) -- (6.2,0.3) ;
\draw (5.8,0) node{$-\eps_{42}$} ;
\end{scope}

\begin{scope}[xshift=8.9cm,yshift=-9cm]
\draw [ultra thick] (0,0) -- ++(0,1.5) ;
\draw [ultra thick] (0.8,0) -- ++(0,1.5) ;
\draw [ultra thick] (1.6,0) -- ++(0,1.5) ;
\draw [ultra thick] (2.4,0) -- ++(0,1.5) ;
\draw [>=stealth,->] (0,1.2) -- (1.6,1.2) ;
\draw (0.4,1.4) node{$\eps_{13}$} ;
\draw [>=stealth,->] (0,0.75) -- (0.8,0.75) ;
\draw (0.4,0.45) node{$\eps_{12}$} ;
\draw [>=stealth,->] (2.4,0.75) -- (1.6,0.75) ;
\draw (2,0.95) node{$\eps_{43}$} ;
\draw [>=stealth,->] (2.4,0.3) -- (0.8,0.3) ;
\draw (2,0) node{$\eps_{42}$} ;
\draw [<->] (2.9,0.75) -- (3.3,0.75) ;
\draw (3.1,1.1) node{wBP} ;
\draw (3.1,0.35) node{$(\ast)$} ;
\draw [ultra thick] (3.8,0) -- ++(0,1.5) ;
\draw [ultra thick] (4.6,0) -- ++(0,1.5) ;
\draw [ultra thick] (5.4,0) -- ++(0,1.5) ;
\draw [ultra thick] (6.2,0) -- ++(0,1.5) ;
\draw [>=stealth,->] (5.4,1.2) -- (3.8,1.2) ;
\draw (4.2,1.4) node{$\eps_{13}$} ;
\draw [>=stealth,->] (4.6,0.75) -- (3.8,0.75) ;
\draw (4.2,0.45) node{$-\eps_{12}$} ;
\draw [>=stealth,->] (5.4,0.75) -- (6.2,0.75) ;
\draw (5.8,0.95) node{$-\eps_{43}$} ;
\draw [>=stealth,->] (4.6,0.3) -- (6.2,0.3) ;
\draw (5.8,0) node{$-\eps_{42}$} ;
\end{scope}

\begin{scope}[xshift=4.9cm,yshift=-11.5cm]
\draw [ultra thick] (0,0) -- ++(0,1.5) ;
\draw [ultra thick] (0.8,0) -- ++(0,1.5) ;
\draw [ultra thick] (1.6,0) -- ++(0,1.5) ;
\draw [ultra thick] (2.4,0) -- ++(0,1.5) ;
\draw [>=stealth,->] (0,1.2) -- (1.6,1.2) ;
\draw (0.4,1.4) node{$\eps_{13}$} ;
\draw [>=stealth,->] (0,0.75) -- (0.8,0.75) ;
\draw (0.4,0.45) node{$\eps_{12}$} ;
\draw [>=stealth,->] (2.4,0.75) -- (1.6,0.75) ;
\draw (2,0.95) node{$\eps_{43}$} ;
\draw [>=stealth,->] (2.4,0.3) -- (0.8,0.3) ;
\draw (2,0) node{$\eps_{42}$} ;
\draw [<->] (2.9,0.75) -- (3.3,0.75) ;
\draw (3.1,1.1) node{BV} ;
\draw (3.1,0.35) node{$(\ast)$} ;
\draw [ultra thick] (3.8,0) -- ++(0,1.5) ;
\draw [ultra thick] (4.6,0) -- ++(0,1.5) ;
\draw [ultra thick] (5.4,0) -- ++(0,1.5) ;
\draw [ultra thick] (6.2,0) -- ++(0,1.5) ;
\end{scope}
\end{tikzpicture}
\caption{Some local moves on Gauss diagrams}
\end{figure}

\begin{de}\label{generate}\cite{fused}
Let $M_1$ and $M_2$ be local moves on classical (resp. virtual) string link diagrams. We say that $M_2$ \emph{c-generates} (resp. \emph{w-generates}) $M_1$ if $M_1$ can be realized using $M_2$ and classical (resp. welded) Reidemeister moves. We denote it by $M_2\impc M_1$ (resp. $M_2\impw M_1$). If $M_1\impc M_2$ and $M_2\impc M_1$ (resp. $M_1\impw M_2$ and $M_2\impw M_1$), we say that $M_1$ and $M_2$ are \emph{c-equivalent} (resp. \emph{w-equivalent}).
\end{de}

\paragraph{Examples:}
As proven in \cite{fused}, we have the following relations:
$$V\impw CC,\quad SV\impw SC,\quad F\overset{w}{\Leftrightarrow}UC,\quad VC\impw F\impw\Delta,\quad wBP\impw SR,F,BP.$$

\subsubsection{Welded extension}

For a classical (resp. virtual) local move $M$, we denote by $\SL^M$ (resp. $w\SL^M$) the quotient of classical (resp. welded) string links by the equivalence relation induced by this move. We then have $M_2\overset{c}{\Rightarrow}M_1$ (resp. $M_2\overset{w}{\Rightarrow}M_1$) if and only if the identity map of $\SL$ (resp. $w\SL$) induces a well defined map $\SL^{M_1}\rightarrow\SL^{M_2}$ (resp. $w\SL^{M_1}\rightarrow w\SL^{M_2}$). \\

If $M_c$ is a classical local move and $M_w$ is a local move such that $M_w\impw M_c$, then the inclusion $SLD_n\to vSLD_n$ induces a map $\SL^{M_c}\to w\SL^{M_w}$.

\begin{de}\cite{fused}
When a local move $M_w$ w-generates a classical move $M_c$, we say that $M_w$  is a \emph{welded extension} of $M_c$ if the induced map $\SL^{M_c}\to w\SL^{M_w}$ is injective.
\end{de}
The local move $M_w$ extends $M_c$ in the sense that if two classical diagrams are related by welded Reidemeister moves and $M_w$, then they are also related by classical Reidemeister moves and $M_c$.

\begin{de}\cite{fused}
Let $M$ be a classical (resp. virtual) local move, $A$ a monoid and $\phi:SLD_n\to A$ (resp. $\phi:vSLD_n\to A$) a monoid homomorphism. We say that $\phi$ \emph{c-classifies} (resp. \emph{w-classifies}) $M$ if it is preserved by $M$  and classical (resp. welded) Reidemeister moves, and the induced map $\overline{\phi}:\SL^M\to A$ (resp. $\overline{\phi}:w\SL^M\to A$) is an isomorphism.
\end{de}

\begin{de}\cite{GouPolVir}
For $i,j\in\{1,\dots,n\}$, $i\neq j$, we define the \emph{virtual linking number} $\text{vlk}_{ij}:vSLD_n\to\ZZ$ by counting, with signs, the number of crossings where the $i^{th}$ component passes over the $j^{th}$ one. If $D\in vSLD_n$ is represented by a Gauss diagram, $\text{vlk}_{ij}(D)$ is the number of arrows going from the $i^{th}$ to $j^{th}$ strands, counted with their sign. We note that the (classical) linking number $\text{lk}_{ij}:SLD_n\to\ZZ$ is the restriction of $\text{vlk}_{ij}$ to classical string link diagrams.
\end{de}

The (virtual) linking numbers are preserved by classical and welded Reidemeister moves, hence they are well-defined on classical and welded string links. By taking combinations of the linking numbers, we obtain classifying invariants for certain local moves. We will make use of $\text{vlk}_{i\ast}:=\sum_{j\neq i}\text{vlk}_{ij}$, $\text{vlk}_{\ast i}:=\sum_{j\neq i}\text{vlk}_{ji}$ and $\text{vlk}_{ij}^{\text{mod}}:=\text{vlk}_{ij}\;\text{mod}\;2\in\ZZ_2$. The same notation is used on classical linking numbers.

\begin{prop}\label{classification}
We have the following classification:
\begin{itemize}
\item \cite{MurNak} $(\text{lk}_{ij})_{1\leq i<j\leq n}:SLD_n\to\ZZ^{n(n-1)/2}$ c-classifies $\Delta$ ;
\item \cite{fused} $(\text{vlk}_{ij}-\text{vlk}_{ji})_{1\leq i<j\leq n}:vSLD_n\to\ZZ^{n(n-1)/2}$ w-classifies $CC$ ;
\item \cite{fused} $(\text{vlk}_{ij})_{1\leq i\neq j\leq n}:vSLD_n\to\ZZ^{n(n-1)}$ w-classifies $F$ ;
\item \cite{fused} $(\text{vlk}_{ij}+\text{vlk}_{ji})_{1\leq i<j\leq n}:vSLD_n\to\ZZ^{n(n-1)/2}$ w-classifies $VC$ ;
\item \cite{Kaw} and \cite{MurNak} $(\text{lk}_{i\ast}^{\text{mod}})_{1\leq i\leq n-1}:SLD_n\to\ZZ_2^{n-1}$ c-classifies $BP$ ;
\item \cite{fused} $(\text{vlk}_{ij}^{\text{mod}}+\text{vlk}_{ji}^{\text{mod}})_{1\leq i<j\leq n}\oplus(\text{vlk}_{i\ast}^{\text{mod}})_{1\leq i\leq n-1}:vSLD_n\to\ZZ_2^{n(n-1)/2}\oplus\ZZ_2^{n-1}$ w-classifies $wBP$.
\end{itemize}
\end{prop}

From this classification, we obtain the following extensions:

\begin{cor}\label{weldedext}\cite{fused}
The $F$ and $VC$ moves both extend $\Delta$, and the $wBP$ move extends $BP$.
\end{cor}

It can be noted that, as in the case of $\Delta$, a classical move can have several welded extensions. When so, we will try to use the relation between string links and knotted surfaces to isolate one extension among the others.

\subsubsection{Ribbon residue}

In this section, we provide two ways to extend the action of the local moves on (welded) string links to (ribbon) string 2--links using the Tube and Spun maps. The final goal is to compare the action of a move Tube($M_w$) with the restriction on ribbon string 2--links of the action of Spun($M_c$) when $M_w$ is a welded extension of $M_c$.

\begin{de}
Let $M$ be a local move. Let us consider the binary relation on $\dr{2}$ which identifies two elements $R_1,R_2\in\dr{2}$ if there exists $L_1\in\text{Tube}^{-1}(R_1),L_2\in\text{Tube}^{-1}(R_2)$ such that $L_1$ and $L_2$ are equivalent under $M$. This relation is reflexive and symmetric, but not necessarily transitive, so we define Tube($M$) as the transitive closure of this relation, which is an equivalence relation on $\dr{2}$.
\end{de}

\paragraph{Remark:}
If $M$ w-generates the $SV$ move, we can show (see the last remark in Section 2.1.3) that the binary relation defined above is already an equivalence relation, so there is no need to take its transitive closure. In this case, the induced map $\text{Tube}:w\SL^M\to\dr{2}/\text{Tube}(M)$ is bijective. Hence for two local moves $M_1$ and $M_2$, each w-generating $SV$, we have $\text{Tube}(M_1)=\text{Tube}(M_2)$ if and only if $M_1$ and $M_2$ are w-equivalent.

\begin{de}
Let $L$ and $L'$ be string 2--links represented by broken surface diagrams $D$ and $D'$, respectively. For a classical local move $M$, we say that $L$ and $L'$ are related by a Spun($M$) move if there exists a solid torus $T\subset B^{2,1}$ such that $D$ and $D'$ are identical except in $T$, where they differ by the spun of the move $M$. We denote by Spun($M$) the equivalence relation on $\dl{2}$ which identifies two such string 2--links $L$ and $L'$.
\end{de}

For a classical local move $M$, we can consider both Spun($M$) and Tube($M$). In general, the restriction of Spun($M$) to ribbon string 2--links is not equal to Tube($M$). For example, for $M=CC$, every ribbon string 2--link is trivial up to Spun($CC$) (see the example below), while $\dr{2}/\text{Tube}(CC)$ is not trivial. This can be proved using the classification of $CC$ on welded string links and the generalization of linking numbers for string 2--links developed in Section 2.2.

\begin{de}\label{residue}
Let $M_c$ be a classical local move. We say that a local move $M_w$ is a \emph{ribbon residue} of $M_c$ if Tube($M_w$) is the restriction to $\dr{2}$ of the equivalence relation Spun($M_c$) on $\dl{2}$. This amounts to say that two ribbon string 2--links are equivalent under Tube($M_w$) if and only if they are equivalent under Spun($M_c$) in the set of string 2--links.
\end{de}

\paragraph{Notation:}
Considering an equivalence relation $\mathcal{R}$ on a set $X$ as its defining subset $\{(x,x')\,|\,x\,\mathcal{R}\,x'\}$ of $X\times X$, for $Y\subset X$ we denote by $\mathcal{R}|_Y:=\mathcal{R}\cap(Y\times Y)$ the restriction of $\mathcal{R}$ to $Y$. If $\mathcal{R}$ and $\mathcal{R}'$ are two equivalence relations on $X$, $\mathcal{R}\subset\mathcal{R}'$ means that $x\,\mathcal{R}\,x'$ implies $x\,\mathcal{R}'\,x'$ for $x,x'\in X$. \\

With this notation, a local move $M_w$ is a ribbon residue of $M_c$ if and only if $\text{Spun}(M_c)|_{2\!-\!r\SL}=\text{Tube}(M_w)$. \\

From the remark above, it follows that a classical move can have at most one ribbon residue which w-generates the $SV$ move. As we will see in Section 2.2, the $F$ and $VC$ moves both w-generate $SV$, and by their classification they are not w-equivalent, so only one of them (if any) can be a ribbon residue of $\Delta$.

\paragraph{Example:}
It is not difficult to see that $V$ is a ribbon residue of $CC$: since $w\SL^{V}$ is trivial, so is $\dr{2}/\text{Tube}(V)$, so it is enough to verify that Tube($V$) can be performed using a Spun($CC$) move. Figure 10 illustrates how this can be done, with the Spun($CC$) move being used in a torus neighborhood of the top right line of double points. We can then use Roseman moves to separate the two tubes. \\

\begin{figure}[h]
\centering
\begin{tikzpicture}
\begin{scope}[scale=0.75]
\begin{scope}[xshift=-4cm,yshift=4cm]
\draw [thick,>=stealth,->] (0,0) -- ++(0.8,0.8) ;
\draw [white,fill=white] (0.4,0.4) circle (0.15) ;
\draw [thick,>=stealth,->] (0.8,0) -- ++(-0.8,0.8) ;
\end{scope}

\begin{scope}[xshift=-4.2cm,yshift=3.1cm]
\draw [thick,->] (0,0) -- ++(-0.7,-1) ;
\draw (-1,-0.3) node{Tube} ;
\end{scope}

\begin{scope}[xshift=-6cm,yshift=0cm]
\begin{scope}[rotate=-45]
\draw (0,-0.75) ellipse (0.3 and 0.12) ;
\fill [white] (-0.3,-0.7) -- (-0.3,-0.5) -- (0.3,-0.5) -- (0.3,-0.7) -- cycle ;
\draw (-0.3,-1.2) -- (-0.3,-0.75) ;
\draw (0.3,-1.2) -- (0.3,-0.75) ;
\end{scope}
\begin{scope}[rotate=-45]
\clip (-0.3,-1.2) -- (-0.3,-1.5) -- (0.3,-1.5) -- (0.3,-1.2) -- cycle ;
\draw (0,-1.2) ellipse (0.3 and 0.12) ;
\end{scope}
\begin{scope}[rotate=-45]
\clip (-0.3,-1.2) -- (-0.3,-1) -- (0.3,-1) -- (0.3,-1.2) -- cycle ;
\draw [dashed] (0,-1.2) ellipse (0.3 and 0.12) ;
\end{scope}
\begin{scope}[rotate=45]
\draw (0,1.2) ellipse (0.7 and 0.28) ;
\draw (-0.7,-1.2) -- (-0.7,1.2) ;
\draw (0.7,-1.2) -- (0.7,1.2) ;
\end{scope}
\begin{scope}[rotate=45]
\clip (-0.7,-1.2) -- (-0.7,-1.5) -- (0.7,-1.5) -- (0.7,-1.2) -- cycle ;
\draw (0,-1.2) ellipse (0.7 and 0.28) ;
\end{scope}
\begin{scope}[rotate=45]
\clip (-0.7,-1.2) -- (-0.7,-0.9) -- (0.7,-0.9) -- (0.7,-1.2) -- cycle ;
\draw [dashed] (0,-1.2) ellipse (0.7 and 0.28) ;
\end{scope}
\begin{scope}[rotate=-45]
\draw (0,0.35) ellipse (0.5 and 0.2) ;
\fill [white] (-0.3,0.35) -- (-0.3,0.8) -- (0.3,0.8) -- (0.3,0.35) -- cycle ;
\draw (0,1.2) ellipse (0.3 and 0.12) ;
\draw [dashed] (-0.3,-0.4) -- (-0.3,0.2) ;
\draw (-0.3,0.2) -- (-0.3,1.2) ;
\draw [dashed] (0.3,-0.4) -- (0.3,0.2) ;
\draw [dashed] (0,-0.4) ellipse (0.3 and 0.12) ;
\draw (0.3,0.2) -- (0.3,1.2) ;
\draw [dashed] (0,0.35) ellipse (0.2 and 0.08) ;
\end{scope}
\end{scope}

\begin{scope}[xshift=-3.6cm,yshift=0cm]
\draw [thick,->] (0,0) -- ++(1,0) ;
\draw (0.5,0.5) node{Spun(CC)} ;
\end{scope}

\begin{scope}[xshift=0cm,yshift=0cm]
\begin{scope}[rotate=-45]
\draw (0,-0.75) ellipse (0.3 and 0.12) ;
\fill [white] (-0.3,-0.7) -- (-0.3,-0.5) -- (0.3,-0.5) -- (0.3,-0.7) -- cycle ;
\draw (-0.3,-1.2) -- (-0.3,-0.75) ;
\draw (0.3,-1.2) -- (0.3,-0.75) ;
\end{scope}
\begin{scope}[rotate=-45]
\clip (-0.3,-1.2) -- (-0.3,-1.5) -- (0.3,-1.5) -- (0.3,-1.2) -- cycle ;
\draw (0,-1.2) ellipse (0.3 and 0.12) ;
\end{scope}
\begin{scope}[rotate=-45]
\clip (-0.3,-1.2) -- (-0.3,-1) -- (0.3,-1) -- (0.3,-1.2) -- cycle ;
\draw [dashed] (0,-1.2) ellipse (0.3 and 0.12) ;
\end{scope}
\begin{scope}[rotate=45]
\draw (0,1.2) ellipse (0.7 and 0.28) ;
\draw (-0.7,-1.2) -- (-0.7,1.2) ;
\draw (0.7,-1.2) -- (0.7,1.2) ;
\end{scope}
\begin{scope}[rotate=45]
\clip (-0.7,-1.2) -- (-0.7,-1.5) -- (0.7,-1.5) -- (0.7,-1.2) -- cycle ;
\draw (0,-1.2) ellipse (0.7 and 0.28) ;
\end{scope}
\begin{scope}[rotate=45]
\clip (-0.7,-1.2) -- (-0.7,-0.9) -- (0.7,-0.9) -- (0.7,-1.2) -- cycle ;
\draw [dashed] (0,-1.2) ellipse (0.7 and 0.28) ;
\end{scope}
\begin{scope}[rotate=-45]
\fill [white] (-0.3,0.35) -- (-0.3,0.8) -- (0.3,0.8) -- (0.3,0.35) -- cycle ;
\draw (0,1.2) ellipse (0.3 and 0.12) ;
\draw [dashed] (0,0.3) ellipse (0.3 and 0.12) ;
\fill [white] (-0.17,0.3) -- (-0.17,0.5) -- (0.17,0.5) -- (0.17,0.3) -- cycle ;
\draw [dashed] (-0.3,-0.4) -- (-0.3,0.3) ;
\draw (-0.3,0.5) -- (-0.3,1.2) ;
\draw [dashed] (0.3,-0.4) -- (0.3,0.3) ;
\draw [dashed] (0,-0.4) ellipse (0.3 and 0.12) ;
\draw (0.3,0.5) -- (0.3,1.2) ;
\end{scope}
\begin{scope}[rotate=-45]
\clip (-0.3,0.5) -- (-0.3,0.3) -- (0.3,0.3) -- (0.3,0.5) -- cycle ;
\draw (0,0.5) ellipse (0.3 and 0.12) ;
\end{scope}
\begin{scope}[rotate=-45]
\clip (-0.3,0.5) -- (-0.3,0.7) -- (0.3,0.7) -- (0.3,0.5) -- cycle ;
\draw [dashed] (0,0.5) ellipse (0.3 and 0.12) ;
\end{scope}
\end{scope}

\begin{scope}[xshift=2.5cm,yshift=0cm]
\draw [thick,->] (0,0) -- ++(1,0) ;
\draw (0.5,1) node{Roseman} ;
\draw (0.5,0.5) node{moves} ;
\end{scope}

\begin{scope}[xshift=6cm,yshift=0cm]
\begin{scope}[rotate=-45]
\draw (0,1.2) ellipse (0.3 and 0.12) ;
\draw (-0.3,-1.2) -- (-0.3,1.2) ;
\draw (0.3,-1.2) -- (0.3,1.2) ;
\end{scope}
\begin{scope}[rotate=-45]
\clip (-0.3,-1.2) -- (-0.3,-1.5) -- (0.3,-1.5) -- (0.3,-1.2) -- cycle ;
\draw (0,-1.2) ellipse (0.3 and 0.12) ;
\end{scope}
\begin{scope}[rotate=-45]
\clip (-0.3,-1.2) -- (-0.3,-0.5) -- (0.3,-0.5) -- (0.3,-1.2) -- cycle ;
\draw [dashed] (0,-1.2) ellipse (0.3 and 0.12) ;
\end{scope}
\begin{scope}[rotate=45]
\draw (0,1.2) ellipse (0.7 and 0.28) ;
\draw (-0.7,-1.2) -- (-0.7,-0.3) ;
\draw (-0.7,0.3) -- (-0.7,1.2) ;
\draw (0.7,-1.2) -- (0.7,-0.3) ;
\draw (0.7,0.3) -- (0.7,1.2) ;
\end{scope}
\begin{scope}[rotate=45]
\clip (-0.7,-1.2) -- (-0.7,-1.5) -- (0.7,-1.5) -- (0.7,-1.2) -- cycle ;
\draw (0,-1.2) ellipse (0.7 and 0.28) ;
\end{scope}
\begin{scope}[rotate=45]
\clip (-0.7,-1.2) -- (-0.7,-0.5) -- (0.7,-0.5) -- (0.7,-1.2) -- cycle ;
\draw [dashed] (0,-1.2) ellipse (0.7 and 0.28) ;
\end{scope}
\end{scope}

\begin{scope}[xshift=4cm,yshift=3.1cm]
\draw [thick,->] (0,0) -- ++(0.7,-1) ;
\draw (1,-0.3) node{Tube} ;
\end{scope}

\begin{scope}[xshift=3.2cm,yshift=4cm]
\draw [thick,>=stealth,->] (0,0) -- ++(0.8,0.8) ;
\draw [thick,>=stealth,->] (0.8,0) -- ++(-0.8,0.8) ;
\draw (0.4,0.4) circle (0.15) ;
\end{scope}

\begin{scope}[xshift=0cm,yshift=4.4cm]
\draw [thick,->] (-1.5,0) -- (1.5,0) ;
\draw (0,0.5) node{V} ;
\end{scope}
\end{scope}
\end{tikzpicture}
\caption{Performing Tube($V$) using Spun($CC$)}
\end{figure}
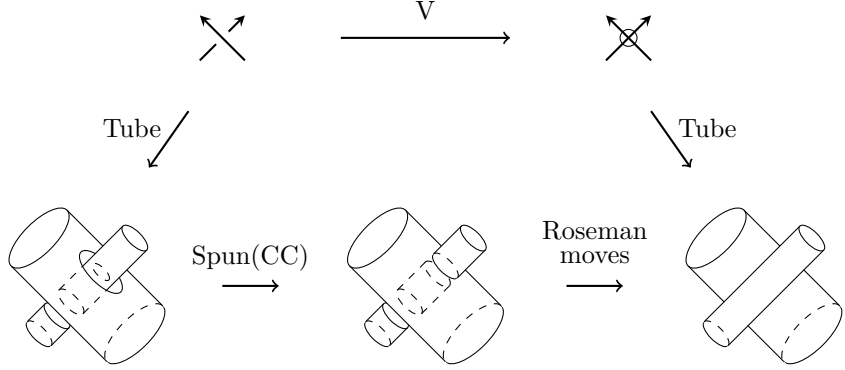

In what follows, we focus on three different cases. In order to prove that a local move $M_w$ is a residue of a classical move $M_c$, we will use the following strategy: first we find a w-classifying invariant $\varphi:w\SL\to A$ of $M_w$, and a homomorphism $\psi:\dl{2}\to A$ which is invariant under Spun($M_c$) and such that $\psi\circ\text{Tube}=\varphi$. This gives $\text{Spun}(M_c)|_{2\!-\!r\SL}\subset\text{Tube}(M_w)$. We then check on broken surface diagrams that Tube($M_w$) can be performed using a Spun($M_c$) move, so that $\text{Spun}(M_c)|_{2\!-\!r\SL}\supset\text{Tube}(M_w)$.

\section{Three examples of residues}

\subsection{The SC and SV moves}

\subsubsection{Classification of the SC move}

We begin by introducing a c-classifying invariant of the $SC$ move, which was established by Habegger and Lin in \cite{HabLin} as a classification of string links up to link homotopy. \\

For a string link $L$, let $X_L$ denote the complement of an open tubular neighborhood of $L$ in $B^{2,1}$. For $\eps=0,1$, $\dd_{\eps}X_L$ is a disk with $n$ smaller and disjoint open disks removed. Hence $\pi_1(\dd_{\eps}X_L)\simeq F_n$, with generators $m_i^{(\eps)}$, called meridians, given by the positively oriented boundaries of these small disks, up to some given fixed path joining them to the basepoint.

\begin{lem}{\cite[Cor. 1.4]{HabLin}}
For $\eps=0,1$, the inclusion maps $\iota_{\eps}:\dd_{\eps}X_L\to X_L$ induce isomorphisms ${\iota_{\eps}^{\ast}:R\pi_1(\dd_{\eps}X_L)\to R\pi_1(X_L)}$.
\end{lem}

Using the identification $R\pi_1(\dd_{\eps}X_L)\simeq RF_n$, we can then define $\varphi_L:=(\iota_0^{\ast})^{-1}\circ\iota_1^{\ast}\in\text{Aut}(RF_n)$. As seen in the Wirtinger presentation of $L$, for each $i$ the meridians $m_i^{(0)}$ and $m_i^{(1)}$ are conjugates of each other, so $\varphi_L\in\text{Aut}_C(RF_n)$. Moreover, the product $m_1^{(0)}\cdots m_n^{(0)}$ is homotopic to $m_1^{(1)}\cdots m_n^{(1)}$ in $X_L$, so $\varphi_L\in\text{Aut}_C^0(RF_n)$.

\begin{prop}\label{phiSLSC}{\cite[Lem. 1.6]{HabLin}}
The map $\Phi_{\SL}:L\in\SL\mapsto\varphi_L\in\text{Aut}_C^0(RF_n)$ is a monoid homomorphism which is invariant under the $SC$ move.
\end{prop}

We obtain the following classification result:

\begin{prop}{\cite[Thm. 1.7]{HabLin}}
The homomorphism $\Phi_{\SL}$ c-classifies $SC$.
\end{prop}

\subsubsection{Classification of the SV move}

We now give a w-classifying invariant of the $SV$ move, which was established in \cite{RTwSL} by Audoux--Bellingeri--Meilhan--Wagner using the notion of coloring on Gauss diagrams. In order to facilitate the transition to string 2-links, we use the equivalent approach of virtual string link diagrams, rather than Gauss diagrams. The correspondence between the two consists in identifying arcs of virtual diagrams with what is referred to as \textquotedblleft tail intervals" of Gauss diagrams in \cite{RTwSL}. \\

We can generalize the Wirtinger presentation to virtual string link diagrams by associating a generator to each arc, and the usual conjugating relation at each classical crossing (and no relation at virtual crossings). For $D\in vSLD_n$, we denote by $\pi_1(D)$ the group given by this presentation.

\begin{de}
Let $\{x_1,\ldots,x_n\}$ be a generating set of $F_n$. If $y_i\in RF_n$ is a conjugate of $x_i$ for each $i$, a $(y_1,\ldots,y_n)$--coloring of a virtual string link diagram $D$ is a map from the arcs of $D$ to $RF_n$, which sends the $i^{th}$ bottom arc to $y_i$, and which satisfies the Wirtinger relation at each classical crossing. Equivalently, this last condition can be replaced by stating that the coloring induces a homomorphism from $\pi_1(D)$ to $RF_n$.
\end{de}

\begin{prop}\label{coloring}{\cite[Lem. 4.20]{RTwSL}}
If $D,D'\in vSLD_n$ are related by one of the wReid or $SV$ moves, then there exists a one-to-one correspondance between the $(y_1,\ldots,y_n)$--colorings of $D$ and $D'$, which preserves the image of the top arcs.
\end{prop}

It is clear that a virtual pure braid diagram (i.e. a diagam with monotone strands) admits a unique $(y_1,\ldots,y_n)$--coloring. Since up to $SV$, any virtual string link diagram is equivalent to a virtual pure braid diagram (see \cite[Thm. 4.12]{RTwSL}), it follows from Proposition~\ref{coloring} that any virtual string link diagram admits a unique $(y_1,\ldots,y_n)$--coloring. Hence for $L\in w\SL$ represented by a diagram $D$, we can define ${\psi_L\in\text{End}_C(RF_n)}$ by $\psi_L(x_i)=z_i$, where $z_i\in RF_n$ is the image of the $i^{th}$ top arc in the unique $(x_1,\ldots,x_n)$--coloring of $D$.

\begin{prop}\label{phiwSLhom}{\cite[Lem. 4.20]{RTwSL}}
The map $\Phi_{w\SL}:L\in w\SL\mapsto\psi_L\in\text{End}_C(RF_n)$ is a monoid homomorphism which is invariant under the $SV$ move. Moreover, we have $\psi_L\in\text{Aut}_C(RF_n)$.
\end{prop}

We obtain the following classification result:

\begin{prop}\label{invSV}{\cite[Thm. 4.17]{RTwSL}}
The homomorphism $\Phi_{w\SL}:w\SL\to\text{Aut}_C(RF_n)$ w-classifies $SV$.
\end{prop}

\paragraph{Remark:}
If $L\in\SL$ and $\pi_1(X_L)$ is given by the Wirtinger presentation associated to a diagram $D$ of $L$, then $(\iota_0^{\ast})^{-1}:R\pi_1(X_L)\simeq R\pi_1(D)\to R\pi_1(\dd_0X_L)\simeq RF_n$ gives an $(x_1,\ldots,x_n)$--coloring of $D$, and since $\iota_1^{\ast}$ sends $x_i$ to the generator of $\pi_1(X_L)$ associated to the $i^{th}$ top arc, we have $\psi_L(x_i)=(\iota_0^{\ast})^{-1}(\iota_1^{\ast}(x_i))=\varphi_L(x_i)$ for $\psi_L=\Phi_{w\SL}(L)$ and $\varphi_L=\Phi_{\SL}(L)$, so $\Phi_{w\SL}(L)=\Phi_{\SL}(L)$. \\

In particular, we obtain that the inclusion $SLD_n\to vSLD_n$ induces an injection $\SL^{SC}\to w\SL^{SV}$, so $SV$ is a welded extension of $SC$. We will now see that it is also a ribbon residue of $SC$.

\subsubsection{Extension to string 2--links}

Let $L\in\dl{2}$ be a string 2--link, $X_L$ the complement of an open tubular neighborhood of $L$, and $D$ a broken surface diagram of $L$. As described in \cite{S4S}, we get a Wirtinger presentation of $\pi_1(X_L)$ from $D$, with one generator for each connected component, called \emph{oversheet}, and a relation of the form $g_+=g_0^{-1}g_-g_0$ at a line of double points as indicated in the figure below:

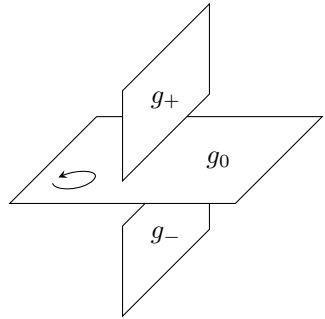
\begin{center}
\begin{tikzpicture}
\draw (0,-0.3,-1.5) -- (0,-0.3,1.5) -- (0,-1.5,1.5) -- (0,-1.5,-1.5) -- cycle ;
\draw [fill=white] (-1.5,0,-1.5) -- (-1.5,0,1.5) -- (1.5,0,1.5) -- (1.5,0,-1.5) -- cycle ;
\draw [fill=white] (0,0.3,-1.5) -- (0,0.3,1.5) -- (0,1.5,1.5) -- (0,1.5,-1.5) -- cycle ;
\draw (0,-1,0) node{$g_-$} ;
\draw (0.7,0,0) node{$g_0$} ;
\draw (0,0.8,0) node{$g_+$} ;
\draw (-1.2,0,0.8) .. controls +(0.15,0,0.3) and +(0,0,0.35) .. (-0.7,0,0.7) ;
\draw [>=stealth,->] (-0.7,0,0.7) .. controls +(0,0,-0.35) and +(0.15,0,-0.3) .. (-1.2,0,0.6) ;
\end{tikzpicture}
\end{center}

\begin{de}
If $y_i\in RF_n$ is a conjugate of $x_i$ for each $i$, a $(y_1,\ldots,y_n)$--coloring of a broken surface diagram $D$ is a map from the oversheets of $D$ to $RF_n$, which sends the $i^{th}$ bottom oversheet to $y_i$, and which satisfies the Wirtinger relation at each line of double points.
\end{de}

As in the case of string links, we have:

\begin{prop}{\cite[\S 4.1]{c2elh}}
For a string 2--link $L$, the inclusion maps $\iota_{\eps}:\dd_{\eps}X_L\to X_L$ induce isomorphisms $\iota_{\eps}^{\ast}:R\pi_1(\dd_{\eps}X_L)\simeq RF_n\to R\pi_1(X_L)$.
\end{prop}

\begin{prop}\label{spuntubecol}
For $L\in\SL$ (resp. $L\in w\SL$), there exists a one-to-one correspondence between the set of $(y_1,\ldots,y_n)$--colorings of $L$ and that of Spun($L$) (resp. of Tube($L$)), which preserves the image of the top and bottom arc/oversheet of each component.
\end{prop}

\paragraph{Proof:}
In the case of the Spun map, it follows directly from the fact that a broken surface diagram of Spun($L$) can be obtained from a diagram of $L$ by a rotation, which sends arcs to oversheets and crossings to lines of double points with the same Wirtinger relations. \\

In the case of the Tube map, we begin by taking a broken surface diagram $D$ of Tube($L$), obtained by the construction described after Definition~\ref{defTube}. There is a correspondence between arcs of $L$ and oversheets of $D$, except for one additional small disk in $D$ at each ribbon singularity. However it is easily seen that given the images of the other oversheets of $D$, the Wirtinger relations give a unique value for these disks, hence removing all ambiguity.
\demo \\

Using the invariance of the number of colorings up to link-homotopy and the fact that any string 2--link is link-homotopic to a ribbon one, it was proven in \cite{c2elh} that a string 2--link admits a unique $(y_1,\ldots,y_n)$--coloring. From this we can define $\theta_L\in\text{Aut}_C(RF_n)$ for $L\in\dl{2}$ in the same way as $\psi_L$ was for a welded string link $L$. Moreover, the map $\phidSL:L\in\dl{2}\mapsto\theta_L\in\text{Aut}_C(RF_n)$ is a monoid homomorphism.

\begin{prop}\label{phi2lh}
The homomorphism $\phidSL$ is invariant under the Spun($SC$) move.
\end{prop}

\paragraph{Proof:}
Let $L\in\dl{2}$ be a string 2--link, represented by a broken surface diagram $D$. Executing a Spun($SC$) move on $D$ yields a broken surface diagram $D'$, identical to $D$ outside of a solid torus $T$ in which the move occurs. Inside this torus, $D$ is subdivided in three annuli, which are all coming from the same component $S_i^{1,1}$ of $L$. We illustrate a slice of these annuli in the torus $T$ below, which highlights the Spun($SC$) move.

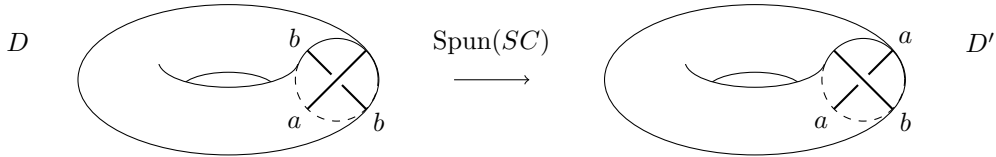
\begin{center}
\begin{tikzpicture}
\draw (0,0) ellipse (2 and 1) ;
\begin{scope}
\clip (-1.2,0.2) -- (-1.2,-0.5) -- (1.2,-0.5) -- (1.2,0.2) -- cycle ;
\draw (0,0.25) ellipse (0.93 and 0.35) ;
\end{scope}
\begin{scope}
\clip (0,0.25) ellipse (0.93 and 0.35) ;
\draw (0,-0.15) ellipse (0.65 and 0.25) ;
\end{scope}
\begin{scope}
\clip (0.5,0.15) -- (0.5,1) -- (2.2,1) -- (2.2,0) -- cycle ;
\draw (1.44,0) circle (0.55) ;
\end{scope}
\begin{scope}
\clip (0.5,0.15) -- (0.5,-1) -- (2.2,-1) -- (2.2,0) -- cycle ;
\draw [dashed] (1.44,0) circle (0.55) ;
\end{scope}
\draw [thick] (1.83,-0.39) -- (1.05,0.39) ;
\fill [white] (1.44,0) circle (0.1) ;
\draw [thick] (1.05,-0.39) -- (1.83,0.39) ;
\draw (0.88,-0.56) node{$a$} ;
\draw (2,-0.56) node{$b$} ;
\draw (0.88,0.56) node{$b$} ;
\draw (-2.8,0.5) node{$D$} ;

\begin{scope}[xshift=3.5cm]
\draw (0,0.5) node{Spun$(SC)$} ;
\draw [->] (-0.5,0) -- (0.5,0) ;
\end{scope}

\begin{scope}[xshift=7cm]
\draw (0,0) ellipse (2 and 1) ;
\begin{scope}
\clip (-1.2,0.2) -- (-1.2,-0.5) -- (1.2,-0.5) -- (1.2,0.2) -- cycle ;
\draw (0,0.25) ellipse (0.93 and 0.35) ;
\end{scope}
\begin{scope}
\clip (0,0.25) ellipse (0.93 and 0.35) ;
\draw (0,-0.15) ellipse (0.65 and 0.25) ;
\end{scope}
\begin{scope}
\clip (0.5,0.15) -- (0.5,1) -- (2.2,1) -- (2.2,0) -- cycle ;
\draw (1.44,0) circle (0.55) ;
\end{scope}
\begin{scope}
\clip (0.5,0.15) -- (0.5,-1) -- (2.2,-1) -- (2.2,0) -- cycle ;
\draw [dashed] (1.44,0) circle (0.55) ;
\end{scope}
\draw [thick] (1.05,-0.39) -- (1.83,0.39) ;
\fill [white] (1.44,0) circle (0.1) ;
\draw [thick] (1.83,-0.39) -- (1.05,0.39) ;
\draw (0.88,-0.56) node{$a$} ;
\draw (2,0.56) node{$a$} ;
\draw (2,-0.56) node{$b$} ;
\draw (3,0.5) node{$D'$} ;
\end{scope}
\end{tikzpicture}
\end{center}

The elements $a,b\in RF_n$ on the left are the images of these annuli in the unique $(x_1,\ldots,x_n)$--coloring of $D$. Since they come from the component $S_i^{1,1}$, $a$ and $b$ are both conjugates of $x_i$, which commute in $RF_n$. Hence the conjugating relation at the line of double points is trivial, and we obtain the same element of $RF_n$ on either side of this line of double points. Since the same holds for $D'$, we can obtain the $(x_1,\ldots,x_n)$--coloring of $D'$ by taking the coloring of $D$ outside of $T$, and extending it inside of $T$ according to its values on $\dd T$. In particular, we have $\theta_L(x_j)=\theta_{L'}(x_j)$ for all $j$, where $L'\in\dl{2}$ is represented by $D'$, and thus $\phidSL(L)=\phidSL(L')$.
\demo \\

As a direct corollary of Proposition~\ref{spuntubecol}, we obtain the following relations between the Spun and Tube maps and the invariants on string links, welded string links and string $2$--links described above:

\begin{prop}\label{comp}
We have $\phidSL\circ\text{Tube}=\Phi_{w\SL}$ and $\phidSL\circ\text{Spun}=\Phi_{\SL}$.
\end{prop}

We can now prove the main result of this section, i.e. the relation between Spun($SC$) and Tube($SV$).

\begin{theo}\label{SVresSC}
The $SV$ move is a ribbon residue of the $SC$ move.
\end{theo}

\paragraph{Proof:}
As indicated earlier, we first prove that $\text{Spun}(SC)|_{2\!-\!r\SL}\subset\text{Tube}(SV)$. Let $R_1,R_2\in\dr{2}$ be equivalent under Spun($SC$). Let $L_1\in\text{Tube}^{-1}(R_1)$ and $L_2\in\text{Tube}^{-1}(R_2)$. By Propositions~\ref{phi2lh} and~\ref{comp}, we have:
$$\Phi_{w\SL}(L_1)=\phidSL(R_1)=\phidSL(R_2)=\Phi_{w\SL}(L_2),$$
and since $\Phi_{w\SL}$ is a w-classifying invariant of $SV$, $L_1$ and $L_2$ are equivalent under $SV$. By definition, this implies that $R_1$ and $R_2$ are equivalent under Tube($SV$). \\

The other inclusion $\text{Spun}(SC)|_{2\!-\!r\SL}\supset\text{Tube}(SV)$ follows by Figure 10 when restricting ourselves to ribbon singularities of a component with itself.
\demo

\paragraph{Remark:}
Another consequence of Proposition~\ref{comp} is that if $L,L'\in w\SL$ have the same image by the Tube map, then $\Phi_{w\SL}(L)=\Phi_{w\SL}(L')$, so $L$ and $L'$ are equivalent under the $SV$ move. This is also true for any local move $M$ w-generating $SV$, which proves the remark given in Section 1.3.3.

\subsection{The $\Delta$ and F moves}

We begin by defining the linking numbers for string 2--links. This will be done using the $\phidSL$ map defined in the previous section, so let us first show how the virtual linking numbers on welded string links can be obtained via $\Phi_{w\SL}$. \\

For $1\leq j\leq n$, let $RF_{n-1}^{(j)}$ denote the subgroup of $RF_n$ generated by the $x_k$'s for $k\neq j$. For a conjugating automorphism $\varphi\in\text{Aut}_C(RF_n)$, let $\lbd_j\in RF_n$ be such that $\varphi(x_j)=\lbd_j^{-1}x_j\lbd_j$. As proven in {\cite[Lem. 4.25]{RTwSL}}, this $\lbd_j$ can be modified into an element of $RF_{n-1}^{(j)}$ by deleting occurrences of $x_j$ in it, and is then uniquely determined by $\varphi$. \\

Applying this to $\psi_L\in\text{Aut}_C(RF_n)$ for a welded string link $L$, we obtain unique elements $\lbd_j\in RF_{n-1}^{(j)}$ for $1\leq j\leq n$ such that $\psi_L(x_j)=\lbd_j^{-1}x_j\lbd_j$. Let $x_i^{\ast}:RF_n\to\ZZ$ be the group homomorphism defined by $x_i^{\ast}(x_i)=1$ and $x_i^{\ast}(x_k)=0$ for $k\neq i$. Keeping track of the undercrossings of $L$ on the $j^{th}$ strand, it is not difficult to check that $\text{vlk}_{ij}(L)=x_i^{\ast}(\lbd_j)$. \\

By making use of the invariant $\phidSL(L)\in\text{Aut}_C(RF_n)$ for $L\in\dl{2}$, this gives us a way of extending the notion of linking numbers to string 2--links.

\begin{de}\label{defLK}
For $i\neq j\in\{1,\ldots,n\}$, we define the \emph{linking number} $\text{LK}_{ij}:\dl{2}\to\ZZ$ by $\text{LK}_{ij}(L):=x_i^{\ast}(\lbd_j)$, where $\lbd_j$ is the unique element of $RF_{n-1}^{(j)}$ such that $\theta_L(x_j)=\lbd_j^{-1}x_j\lbd_j$ for $\theta_L=\phidSL(L)$.
\end{de}

\begin{prop}\label{LKDelta}
The map $\text{LK}_{ij}$ is a monoid homomorphism, which is invariant under the Spun($\Delta$) move.
\end{prop}

\paragraph{Proof:}
That $\text{LK}_{ij}$ is a monoid homomorphism follows from the fact that $L\in\dl{2}\mapsto\theta_L$ is a homomorphism: if $\theta_L(x_j)=\lbd_j^{-1}x_j\lbd_j$ and $\theta_{L'}(x_j)=\lbd_j'^{-1}x_j\lbd_j'$, then:
$$\theta_{L\bullet L'}(x_j)=(\theta_L\circ\theta_{L'})(x_j)=\lbd_j^{-1}\lbd_j'^{-1}x_j\lbd_j'\lbd_j,$$
so $\text{LK}_{ij}(L\bullet L')=x_i^{\ast}(\lbd_j'\lbd_j)=x_i^{\ast}(\lbd_j')+x_i^{\ast}(\lbd_j)=\text{LK}_{ij}(L')+\text{LK}_{ij}(L)$. \\

In order to show that $\text{LK}_{ij}$ is invariant under the Spun($\Delta$) move, we give a more geometric interpretation as follows: let $L\in\dl{2}$ be represented by a broken surface diagram $D$, with $D_1,\ldots,D_n\subset D$ denoting the $n$ components of the underlying string 2--link. Let $\gam_j$ be a path on $D_j$ from $\dd_0D_j$ to $\dd_1D_j$, having transverse intersections with lines of double points. Then $\text{LK}_{ij}(L)$ is the number of times (counted with a sign given by the orientation) $\gam_j$ crosses a line of double points where $D_j$ passes behind $D_i$. \\

From this interpretation, it follows that $\text{LK}_{ij}(L)$ only depends on the components $D_i$ and $D_j$. A Spun($\Delta$) move modifies the relative position of three components of the diagram, but not the relative position of any pair of components, so it does not affect the linking number of $L$.
\demo

\begin{prop}\label{LKprop}
For $i,j\in\{1,\ldots,n\}$, $i\neq j$, we have $\text{LK}_{ij}\circ\text{Tube}=\text{vlk}_{ij}$ and $\text{LK}_{ij}\circ\text{Spun}=\text{lk}_{ij}$.
\end{prop}

\paragraph{Proof:}
For the Tube map, it follows directly from Proposition~\ref{comp} and Definition~\ref{defLK}. For the Spun map, we also use the fact that $\Phi_{w\SL}|_{\SL}=\Phi_{\SL}$ and $\text{vlk}_{ij}|_{\SL}=\text{lk}_{ij}$.
\demo \\

We can now determine the relation between Spun($\Delta$) and Tube($F$).

\begin{theo}\label{FresDelta}
The $F$ move is a ribbon residue of the $\Delta$ move.
\end{theo}

\paragraph{Proof:}
Let $R_1,R_2\in\dr{2}$ be equivalent under Spun($\Delta$), and let $L_1\in\text{Tube}^{-1}(R_1)$, $L_2\in\text{Tube}^{-1}(R_2)$. By Propositions~\ref{LKDelta} and~\ref{LKprop}, we have $\text{vlk}_{ij}(L_1)=\text{LK}_{ij}(R_1)=\text{LK}_{ij}(R_2)=\text{vlk}_{ij}(L_2)$, and since $F$ is w-classified by the virtual linking numbers, $L_1$ and $L_2$ are equivalent under $F$. Hence $R_1$ and $R_2$ are equivalent under Tube($F$), and we have $\text{Spun}(\Delta)|_{2\!-\!r\SL}\subset\text{Tube}(F)$. \\

Since $F$ is w-equivalent to $UC$, we have $\text{Tube}(F)=\text{Tube}(UC)$, so to prove $\text{Spun}(\Delta)|_{2\!-\!r\SL}\supset\text{Tube}(F)$ it is enough to verify that we can perform a Tube($UC$) move using Spun($\Delta$) on a broken surface diagram. This is illustrated on Figure 11 at the end of the article, with the following conventions:
\begin{itemize}
\item it is to be read from left to right and top to bottom, beginning and ending with the images under Tube of the situations before and after a $UC$ move.
\item on the second row of this figure, we represent a portion of a broken surface diagram of a ribbon string 2--link by only drawing the intersection of the tubes with a median hyperplane, so as to make it more readable. In this representation, crossings correspond to lines of double points, which is why we first used Roseman moves to pull the diagonal blue tube \textquotedblleft inside" the vertical red tube. The represented strands are then treated as pieces of string links, keeping in mind that the moves we use must be compatible with Roseman moves on broken surface diagrams.
\item in the dotted rectangle, the surface is obtained by a rotation around an axis, so we can perform a Spun($\Delta$) move.
\item on the second to last picture, we come back to a proper representation of a broken surface diagram, skipping the intermediate step where the blue tube is pulled outside the red one.
\end{itemize}

This shows how to perform a Tube($UC$) move using Spun($\Delta$), at least for the case where the orientation of the strands matches the one in Figure 11. But it is not difficult to check that the $UC$ move with this specific orientation is w-equivalent to the other $UC$ moves, so there is no loss of generality.
\demo \\

Hence $F$ is a welded extension of $\Delta$ which is also a ribbon residue. From the remark following the definition of a residue, we obtain that the $F$ move is the only $SV$ w-generating welded extension of $\Delta$ which is also a residue; in particular, the $VC$ move is not a ribbon residue of $\Delta$.

\subsection{The BP move}

As we will see in the case of the $BP$ move, for some classical move $M$ the action of Spun($M$) on knotted surfaces can be drastically different from that of $M$ on classical string links, and produce a ribbon residue which fails to be a welded extension. \\

A natural candidate to be a ribbon residue of $BP$ is the $BV$ move, defined in Section 1.3.1. However, it turns out that the action of Spun($BP$) on knotted surfaces is much stronger than that of Tube($BV$). Let us consider the $DV$ move (for \textquotedblleft double virtualization\textquotedblright) defined as follows:

\begin{center}
\begin{tikzpicture}
\draw [thick] (-0.8,0) -- (0.8,0) ;
\fill [white] (-0.3,0) circle (0.15) ;
\fill [white] (0.3,0) circle (0.15) ;
\draw [thick] (-0.3,-0.5) -- (-0.3,0.5) ;
\draw [thick] (0.3,-0.5) -- (0.3,0.5) ;

\begin{scope}[xshift=2cm]
\draw [<->] (-0.4,0) -- (0.4,0) ;
\draw (0,0.3) node{DV} ;
\end{scope}

\begin{scope}[xshift=4cm]
\draw [thick] (-0.8,0) -- (0.8,0) ;
\draw (-0.3,0) circle (0.12) ;
\draw (0.3,0) circle (0.12) ;
\draw [thick] (-0.3,-0.5) -- (-0.3,0.5) ;
\draw [thick] (0.3,-0.5) -- (0.3,0.5) ;
\end{scope}
\end{tikzpicture}
\end{center}

\begin{prop}\label{DVprop}
The $DV$ move is w-equivalent to the $V$ move.
\end{prop}

\paragraph{Proof:}
The $DV$ move can be realized by two $V$ moves, and the $V$ move can be performed using the $DV$ move by introducing a second crossing with $R1$ and $vR1$ moves:

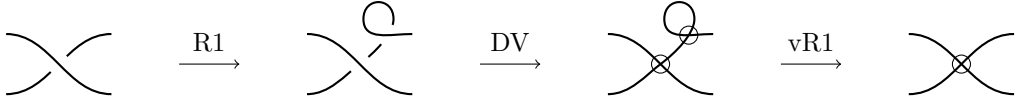
\begin{center}
\begin{tikzpicture}
\draw [thick] (-0.7,-0.4) .. controls +(0.6,0) and +(-0.6,0) .. (0.7,0.4) ;
\fill [white] (0,0) circle (0.15) ;
\draw [thick] (-0.7,0.4) .. controls +(0.6,0) and +(-0.6,0) .. (0.7,-0.4) ;

\begin{scope}[xshift=2cm]
\draw [->] (-0.4,0) -- (0.4,0) ;
\draw (0,0.3) node{R1} ;
\end{scope}

\begin{scope}[xshift=4cm]
\draw [thick] (-0.7,-0.4) .. controls +(0.3,0) and +(-0.2,-0.15) .. (0,0) ;
\draw [thick] (0,0) .. controls +(0.2,0.15) and +(0,-0.2) .. (0.45,0.6) ;
\draw [thick] (0.45,0.6) .. controls +(0,0.3) and +(0,0.3) .. (0.05,0.6) ;
\fill [white] (0.375,0.385) circle (0.15) ;
\draw [thick] (0.05,0.6) .. controls +(0,-0.3) and +(-0.3,0) .. (0.7,0.4) ;
\fill [white] (0,0) circle (0.15) ;
\draw [thick] (-0.7,0.4) .. controls +(0.6,0) and +(-0.6,0) .. (0.7,-0.4) ;
\end{scope}

\begin{scope}[xshift=6cm]
\draw [->] (-0.4,0) -- (0.4,0) ;
\draw (0,0.3) node{DV} ;
\end{scope}

\begin{scope}[xshift=8cm]
\draw [thick] (-0.7,-0.4) .. controls +(0.3,0) and +(-0.2,-0.15) .. (0,0) ;
\draw [thick] (0,0) .. controls +(0.2,0.15) and +(0,-0.2) .. (0.45,0.6) ;
\draw [thick] (0.45,0.6) .. controls +(0,0.3) and +(0,0.3) .. (0.05,0.6) ;
\draw (0.375,0.385) circle (0.12) ;
\draw [thick] (0.05,0.6) .. controls +(0,-0.3) and +(-0.3,0) .. (0.7,0.4) ;
\draw (0,0) circle (0.12) ;
\draw [thick] (-0.7,0.4) .. controls +(0.6,0) and +(-0.6,0) .. (0.7,-0.4) ;
\end{scope}

\begin{scope}[xshift=10cm]
\draw [->] (-0.4,0) -- (0.4,0) ;
\draw (0,0.3) node{vR1} ;
\end{scope}

\begin{scope}[xshift=12cm]
\draw [thick] (-0.7,-0.4) .. controls +(0.6,0) and +(-0.6,0) .. (0.7,0.4) ;
\draw (0,0) circle (0.12) ;
\draw [thick] (-0.7,0.4) .. controls +(0.6,0) and +(-0.6,0) .. (0.7,-0.4) ;
\end{scope}
\end{tikzpicture}
\end{center}
\demo

\begin{theo}\label{VresBP}
The $V$ move is a ribbon residue of the $BP$ move.
\end{theo}

\paragraph{Proof:}
Since the $V$ move identifies all welded string links on $n$ strands, we trivially have $\text{Tube}(V)\supset\text{Spun}(BP)|_{2\!-\!r\SL}$. By Proposition~\ref{DVprop} we have $\text{Tube}(V)=\text{Tube}(DV)$, so for the other inclusion we only need to check that a Tube($DV$) move can be performed using Spun($BP$), which is illustrated on Figure 12, where we use the same convention as in the proof of Theorem~\ref{FresDelta}. \\

To go from picture 3 to 4, we inflate the horizontal black tube and then push part of it inside itself, passing in front of the vertical red and blue tubes, by a finger-move using a Roseman (a) move (in {\cite[Fig. 1]{Roseman}}). This brings us in a position where we can perform a Spun($BP$) move. \\

On this figure, a choice has been made regarding the orientation of the vertical strands (note that the orientation of the horizontal strand is irrelevant). If they have the same orientation, an extra step is needed in between pictures 3 and 4, where we need to slide one of the vertical tubes into the other in order to obtain parallel vertical lines with the same relative position to the horizontal ones.
\demo \\

The action of Spun($BP$) on knotted surfaces is very different from that of $BP$ on classical string links, as it identifies all ribbon string 2--links. This is in part due to the nature of the line of double points involved in a Spun($BP$) move: these lines of double points are circles on annuli, so they can be of two types, essential or contractible. The occurrence of a line of double points whose preimages are both contractible is what can make the Spun of a local move much stronger than the move itself. Indeed, on Figure 12, in order to use a Spun($BP$) move we created lines of double points whose preimages are contractible on the horizontal black annulus \emph{and} the vertical blue and red annuli. This allowed us to go down from four strands involved in a $BP$ move to only three in a $DV$ move, which is much stronger than the $BV$ move involving four strands. \\

Even if it will eventually fail to produce a welded extension as a ribbon residue in this case, a strategy to counteract this feature of the Spun($BP$) move is to introduce some restrictions on the type of lines of double points involved.

\begin{de}
We define the Spun$^{\ast}$($BP$) move by imposing that preimages of the lines of double points inside the torus where the Spun($BP$) move occurs must all be essential on one of the two pairs of parallel annuli, and all be contractible on the other pair. As before, we also denote by Spun$^{\ast}$($BP$) the equivalence relation on $\dl{2}$ which identifies two string 2--links related by this move.
\end{de}

This new move cannot be used to perform a Tube($DV$) move but, as we will see, it is strong enough to obtain Tube($BV$). Let us first provide a classification of the $BV$ move.

\begin{lem}\label{BVwgen}
The $BV$ move w-generates the $F$, $SR$ and $VC$ moves.
\end{lem}

\paragraph{Proof:}
First, we show on virtual diagrams that $BV$ w-generates $F$, as illustrated below: \\

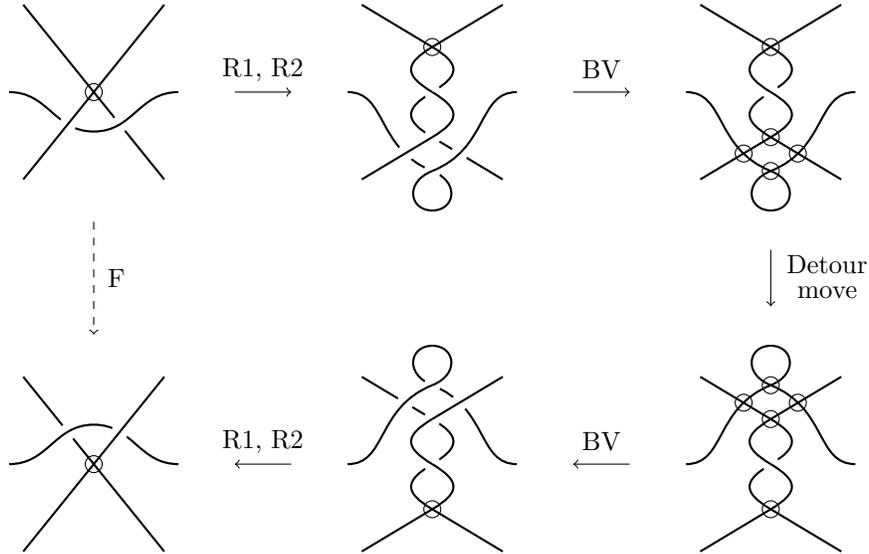
\begin{center}
\begin{tikzpicture}
\begin{scope}[scale=0.75]
\begin{scope}[xshift=-0.5cm]
\draw [thick] (1.25,-1.55) -- (-1.25,1.55) ;
\fill [white] (0.46,-0.57) circle (0.15) ;
\draw [thick] (-1.5,0) .. controls +(0.7,0) and +(-0.7,0) .. (0,-0.7) ;
\draw [thick] (0,-0.7) .. controls +(0.7,0) and +(-0.7,0) .. (1.5,0) ;
\fill [white] (-0.46,-0.57) circle (0.15) ;
\draw [thick] (-1.25,-1.55) -- (1.25,1.55) ;
\draw (0,0) circle (0.15) ;
\end{scope}

\begin{scope}[xshift=2.5cm]
\draw [->] (-0.5,0) -- (0.5,0) ;
\draw (0,0.4) node{R1, R2} ;
\end{scope}

\begin{scope}[xshift=5.5cm]
\draw [thick] (1.25,-1.55) -- (0,-0.8) ;
\draw [thick] (0,-0.8) .. controls +(-0.5,0.3) and +(-0.5,-0.3) .. (0,0) ;
\draw [thick] (0,0) .. controls +(0.5,0.3) and +(0.5,-0.3) .. (0,0.8) ;
\draw [thick] (0,0.8) -- (-1.25,1.55) ;
\fill [white] (0.48,-1.09) circle (0.15) ;
\fill [white] (0,-0.8) circle (0.15) ;
\fill [white] (0,0) circle (0.15) ;
\draw [thick] (-1.5,0) .. controls +(0.7,0) and +(-1,0.4) .. (0,-1.4) ;
\draw [thick] (0,-1.4) .. controls +(0.5,-0.2) and +(0.4,0) .. (0,-2.1) ;
\fill [white] (0,-1.4) circle (0.15) ;
\draw [thick] (0,-2.1) .. controls +(-0.4,0) and +(-0.5,-0.2) .. (0,-1.4) ;
\draw [thick] (0,-1.4) .. controls +(1,0.4) and +(-0.7,0) .. (1.5,0) ;
\fill [white] (-0.48,-1.09) circle (0.15) ;
\draw [thick] (-1.25,-1.55) -- (0,-0.8) ;
\draw [thick] (0,-0.8) .. controls +(0.5,0.3) and +(0.5,-0.3) .. (0,0) ;
\draw [thick] (0,0) .. controls +(-0.5,0.3) and +(-0.5,-0.3) .. (0,0.8) ;
\draw [thick] (0,0.8) -- (1.25,1.55) ;
\draw (0,0.8) circle (0.15) ;
\end{scope}

\begin{scope}[xshift=8.5cm]
\draw [->] (-0.5,0) -- (0.5,0) ;
\draw (0,0.4) node{BV} ;
\end{scope}

\begin{scope}[xshift=11.5cm]
\draw [thick] (1.25,-1.55) -- (0,-0.8) ;
\draw [thick] (0,-0.8) .. controls +(-0.5,0.3) and +(-0.5,-0.3) .. (0,0) ;
\draw [thick] (0,0) .. controls +(0.5,0.3) and +(0.5,-0.3) .. (0,0.8) ;
\draw [thick] (0,0.8) -- (-1.25,1.55) ;
\draw (0.48,-1.09) circle (0.15) ;
\draw (0,-0.8) circle (0.15) ;
\fill [white] (0,0) circle (0.15) ;
\draw [thick] (-1.5,0) .. controls +(0.7,0) and +(-1,0.4) .. (0,-1.4) ;
\draw [thick] (0,-1.4) .. controls +(0.5,-0.2) and +(0.4,0) .. (0,-2.1) ;
\draw (0,-1.4) circle (0.15) ;
\draw [thick] (0,-2.1) .. controls +(-0.4,0) and +(-0.5,-0.2) .. (0,-1.4) ;
\draw [thick] (0,-1.4) .. controls +(1,0.4) and +(-0.7,0) .. (1.5,0) ;
\draw (-0.48,-1.09) circle (0.15) ;
\draw [thick] (-1.25,-1.55) -- (0,-0.8) ;
\draw [thick] (0,-0.8) .. controls +(0.5,0.3) and +(0.5,-0.3) .. (0,0) ;
\draw [thick] (0,0) .. controls +(-0.5,0.3) and +(-0.5,-0.3) .. (0,0.8) ;
\draw [thick] (0,0.8) -- (1.25,1.55) ;
\draw (0,0.8) circle (0.15) ;
\end{scope}

\begin{scope}[xshift=11.5cm,yshift=-3.3cm]
\draw [->] (0,0.5) -- (0,-0.5) ;
\draw (1,0.25) node{Detour} ;
\draw (1,-0.25) node{move} ;
\end{scope}

\begin{scope}[xshift=11.5cm,yshift=-6.6cm]
\draw [thick] (1.25,-1.55) -- (0,-0.8) ;
\draw [thick] (0,-0.8) .. controls +(-0.5,0.3) and +(-0.5,-0.3) .. (0,0) ;
\draw [thick] (0,0) .. controls +(0.5,0.3) and +(0.5,-0.3) .. (0,0.8) ;
\draw [thick] (0,0.8) -- (-1.25,1.55) ;
\draw (0.48,1.09) circle (0.15) ;
\draw (0,-0.8) circle (0.15) ;
\fill [white] (0,0) circle (0.15) ;
\draw [thick] (-1.5,0) .. controls +(0.7,0) and +(-1,-0.4) .. (0,1.4) ;
\draw [thick] (0,1.4) .. controls +(0.5,0.2) and +(0.4,0) .. (0,2.1) ;
\draw (0,1.4) circle (0.15) ;
\draw [thick] (0,2.1) .. controls +(-0.4,0) and +(-0.5,0.2) .. (0,1.4) ;
\draw [thick] (0,1.4) .. controls +(1,-0.4) and +(-0.7,0) .. (1.5,0) ;
\draw (-0.48,1.09) circle (0.15) ;
\draw [thick] (-1.25,-1.55) -- (0,-0.8) ;
\draw [thick] (0,-0.8) .. controls +(0.5,0.3) and +(0.5,-0.3) .. (0,0) ;
\draw [thick] (0,0) .. controls +(-0.5,0.3) and +(-0.5,-0.3) .. (0,0.8) ;
\draw [thick] (0,0.8) -- (1.25,1.55) ;
\draw (0,0.8) circle (0.15) ;
\end{scope}

\begin{scope}[xshift=8.5cm,yshift=-6.6cm]
\draw [->] (0.5,0) -- (-0.5,0) ;
\draw (0,0.4) node{BV} ;
\end{scope}

\begin{scope}[xshift=5.5cm,yshift=-6.6cm]
\draw [thick] (1.25,-1.55) -- (0,-0.8) ;
\draw [thick] (0,-0.8) .. controls +(-0.5,0.3) and +(-0.5,-0.3) .. (0,0) ;
\draw [thick] (0,0) .. controls +(0.5,0.3) and +(0.5,-0.3) .. (0,0.8) ;
\draw [thick] (0,0.8) -- (-1.25,1.55) ;
\fill [white] (-0.48,1.09) circle (0.15) ;
\fill [white] (0,0.8) circle (0.15) ;
\fill [white] (0,0) circle (0.15) ;
\draw [thick] (1.5,0) .. controls +(-0.7,0) and +(1,-0.4) .. (0,1.4) ;
\draw [thick] (0,1.4) .. controls +(-0.5,0.2) and +(-0.4,0) .. (-0,2.1) ;
\fill [white] (0,1.4) circle (0.15) ;
\draw [thick] (0,2.1) .. controls +(0.4,0) and +(0.5,0.2) .. (0,1.4) ;
\draw [thick] (0,1.4) .. controls +(-1,-0.4) and +(0.7,0) .. (-1.5,0) ;
\fill [white] (0.48,1.09) circle (0.15) ;
\draw [thick] (-1.25,-1.55) -- (0,-0.8) ;
\draw [thick] (0,-0.8) .. controls +(0.5,0.3) and +(0.5,-0.3) .. (0,0) ;
\draw [thick] (0,0) .. controls +(-0.5,0.3) and +(-0.5,-0.3) .. (0,0.8) ;
\draw [thick] (0,0.8) -- (1.25,1.55) ;
\draw (0,-0.8) circle (0.15) ;
\end{scope}

\begin{scope}[xshift=2.5cm,yshift=-6.6cm]
\draw [->] (0.5,0) -- (-0.5,0) ;
\draw (0,0.4) node{R1, R2} ;
\end{scope}

\begin{scope}[xshift=-0.5cm,yshift=-6.6cm]
\draw [thick] (1.25,-1.55) -- (-1.25,1.55) ;
\fill [white] (-0.46,0.57) circle (0.15) ;
\draw [thick] (-1.5,0) .. controls +(0.7,0) and +(-0.7,0) .. (0,0.7) ;
\draw [thick] (0,0.7) .. controls +(0.7,0) and +(-0.7,0) .. (1.5,0) ;
\fill [white] (0.46,0.57) circle (0.15) ;
\draw [thick] (-1.25,-1.55) -- (1.25,1.55) ;
\draw (0,0) circle (0.15) ;
\end{scope}

\begin{scope}[xshift=-0.5cm,yshift=-3.3cm]
\draw [dashed,->] (0,1) -- (0,-1) ;
\draw (0.4,0) node{F} ;
\end{scope}
\end{scope}
\end{tikzpicture}
\end{center}

Since $F\overset{w}{\Leftrightarrow}UC$, we can use $UC$ to show that $BV$ w-generates $SR$. It is easily checked that, up to one $R2$ move, the move described below is w-equivalent to $SR$.

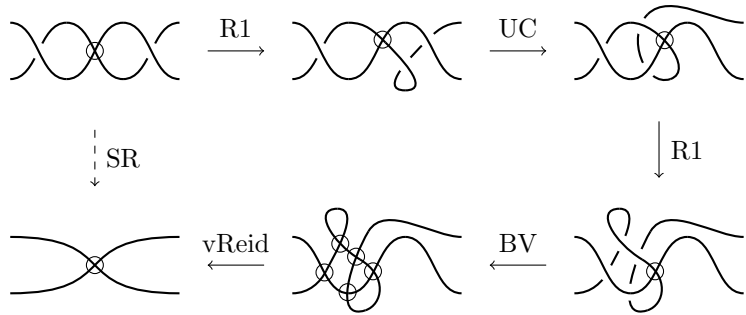
\begin{center}
\begin{tikzpicture}
\begin{scope}[scale=0.75]
\draw [thick] (-1.5,-0.5) .. controls +(0.5,0) and +(-0.5,0) .. (-0.5,0.5) ;
\draw [thick] (-0.5,0.5) .. controls +(0.5,0) and +(-0.5,0) .. (0.5,-0.5) ;
\draw [thick] (0.5,-0.5) .. controls +(0.5,0) and +(-0.5,0) .. (1.5,0.5) ;
\fill [white] (-1,0) circle (0.15) ;
\draw (0,0) circle (0.15) ;
\fill [white] (1,0) circle (0.15) ;
\draw [thick] (-1.5,0.5) .. controls +(0.5,0) and +(-0.5,0) .. (-0.5,-0.5) ;
\draw [thick] (-0.5,-0.5) .. controls +(0.5,0) and +(-0.5,0) .. (0.5,0.5) ;
\draw [thick] (0.5,0.5) .. controls +(0.5,0) and +(-0.5,0) .. (1.5,-0.5) ;

\begin{scope}[xshift=2.5cm]
\draw [->] (-0.5,0) -- (0.5,0) ;
\draw (0,0.4) node{R1} ;
\end{scope}

\begin{scope}[xshift=5cm]
\draw [thick] (-1.5,-0.5) .. controls +(0.5,0) and +(-0.5,0) .. (-0.5,0.5) ;
\draw [thick] (-0.5,0.5) .. controls +(0.3,0) and +(-0.2,0.2) .. (0.15,0.15) ;
\draw [thick] (0.5,-0.7) .. controls +(-0.5,0) and +(-0.2,-0.2) .. (0.85,0.15) ;
\fill [white] (0.5,-0.23) circle (0.15) ;
\draw [thick] (0.15,0.15) .. controls +(0.2,-0.2) and +(0.5,0) .. (0.5,-0.7) ;
\draw [thick] (0.85,0.15) .. controls +(0.2,0.2) and +(-0.3,0) .. (1.5,0.5) ;
\fill [white] (-1,0) circle (0.15) ;
\draw (0.1,0.2) circle (0.15) ;
\fill [white] (0.9,0.2) circle (0.15) ;
\draw [thick] (-1.5,0.5) .. controls +(0.5,0) and +(-0.5,0) .. (-0.5,-0.5) ;
\draw [thick] (-0.5,-0.5) .. controls +(0.5,0) and +(-0.5,0) .. (0.5,0.5) ;
\draw [thick] (0.5,0.5) .. controls +(0.5,0) and +(-0.5,0) .. (1.5,-0.5) ;
\end{scope}

\begin{scope}[xshift=7.5cm]
\draw [->] (-0.5,0) -- (0.5,0) ;
\draw (0,0.4) node{UC} ;
\end{scope}

\begin{scope}[xshift=10cm]
\draw [thick] (-1.5,-0.5) .. controls +(0.5,0) and +(-0.5,0) .. (-0.5,0.5) ;
\draw [thick] (0.2,-0.5) .. controls +(-0.45,-0.15) and +(-0.4,-0.4) .. (-0.2,0.65) ;
\fill [white] (-0.21,-0.34) circle (0.15) ;
\fill [white] (-0.31,0.44) circle (0.15) ;
\draw [thick] (-0.5,0.5) .. controls +(0.5,0) and +(0.45,0.15) .. (0.2,-0.5) ;
\draw [thick] (-0.2,0.65) .. controls +(0.4,0.4) and +(-0.5,0) .. (1.5,0.5) ;
\fill [white] (-1,0) circle (0.15) ;
\draw (0.1,0.19) circle (0.15) ;
\draw [thick] (-1.5,0.5) .. controls +(0.5,0) and +(-0.5,0) .. (-0.5,-0.5) ;
\draw [thick] (-0.5,-0.5) .. controls +(0.5,0) and +(-0.5,0) .. (0.5,0.5) ;
\draw [thick] (0.5,0.5) .. controls +(0.5,0) and +(-0.5,0) .. (1.5,-0.5) ;
\end{scope}

\begin{scope}[xshift=10cm,yshift=-1.75cm]
\draw [->] (0,0.5) -- (0,-0.5) ;
\draw (0.5,0) node{R1} ;
\end{scope}

\begin{scope}[xshift=10cm,yshift=-3.8cm]
\draw [thick] (-1.5,-0.5) .. controls +(0.3,0) and +(-0.1,-0.2) .. (-0.9,-0.1) ;
\draw [thick] (-0.9,-0.1) .. controls +(0.1,0.2) and +(0.45,0) .. (-0.7,1) ;
\fill [white] (-0.65,0.38) circle (0.15) ;
\draw [thick] (-0.2,-0.8) .. controls +(-0.6,-0.2) and +(-0.3,-0.5) .. (-0.15,0.6) ;
\fill [white] (-0.52,-0.48) circle (0.15) ;
\fill [white] (-0.37,0.15) circle (0.15) ;
\draw [thick] (-0.7,1) .. controls +(-0.45,0) and +(-0.5,0.3) .. (-0.3,0.1) ;
\draw [thick] (-0.3,0.1) .. controls +(0.5,-0.3) and +(0.3,0.1) .. (-0.2,-0.8) ;
\draw [thick] (-0.15,0.6) .. controls +(0.3,0.5) and +(-0.5,0) .. (1.5,0.5) ;
\fill [white] (-0.93,-0.13) circle (0.15) ;
\draw (-0.07,-0.11) circle (0.15) ;
\draw [thick] (-1.5,0.5) .. controls +(0.5,0) and +(-0.5,0) .. (-0.5,-0.5) ;
\draw [thick] (-0.5,-0.5) .. controls +(0.5,0) and +(-0.5,0) .. (0.5,0.5) ;
\draw [thick] (0.5,0.5) .. controls +(0.5,0) and +(-0.5,0) .. (1.5,-0.5) ;
\end{scope}

\begin{scope}[xshift=7.5cm,yshift=-3.8cm]
\draw [->] (0.5,0) -- (-0.5,0) ;
\draw (0,0.4) node{BV} ;
\end{scope}

\begin{scope}[xshift=5cm,yshift=-3.8cm]
\draw [thick] (-1.5,-0.5) .. controls +(0.3,0) and +(-0.1,-0.2) .. (-0.9,-0.1) ;
\draw [thick] (-0.9,-0.1) .. controls +(0.1,0.2) and +(0.45,0) .. (-0.7,1) ;
\draw (-0.65,0.38) circle (0.15) ;
\draw [thick] (-0.2,-0.8) .. controls +(-0.6,-0.2) and +(-0.3,-0.5) .. (-0.15,0.6) ;
\draw (-0.52,-0.48) circle (0.15) ;
\draw (-0.37,0.15) circle (0.15) ;
\draw [thick] (-0.7,1) .. controls +(-0.45,0) and +(-0.5,0.3) .. (-0.3,0.1) ;
\draw [thick] (-0.3,0.1) .. controls +(0.5,-0.3) and +(0.3,0.1) .. (-0.2,-0.8) ;
\draw [thick] (-0.15,0.6) .. controls +(0.3,0.5) and +(-0.5,0) .. (1.5,0.5) ;
\draw (-0.93,-0.13) circle (0.15) ;
\draw (-0.07,-0.11) circle (0.15) ;
\draw [thick] (-1.5,0.5) .. controls +(0.5,0) and +(-0.5,0) .. (-0.5,-0.5) ;
\draw [thick] (-0.5,-0.5) .. controls +(0.5,0) and +(-0.5,0) .. (0.5,0.5) ;
\draw [thick] (0.5,0.5) .. controls +(0.5,0) and +(-0.5,0) .. (1.5,-0.5) ;
\end{scope}

\begin{scope}[xshift=2.5cm,yshift=-3.8cm]
\draw [->] (0.5,0) -- (-0.5,0) ;
\draw (0,0.4) node{vReid} ;
\end{scope}

\begin{scope}[xshift=0cm,yshift=-3.8cm]
\draw [thick] (-1.5,-0.5) .. controls +(2,0) and +(-2,-0) .. (1.5,0.5) ;
\draw [thick] (-1.5,0.5) .. controls +(2,0) and +(-2,-0) .. (1.5,-0.5) ;
\draw (0,0) circle (0.15) ;
\end{scope}

\begin{scope}[xshift=0cm,yshift=-1.9cm]
\draw [dashed,->] (0,0.5) -- (0,-0.5) ;
\draw (0.5,0) node{SR} ;
\end{scope}
\end{scope}
\end{tikzpicture}
\end{center}

Hence the sign restriction on the $BV$ move on Gauss diagrams can be omitted. We can now show that $BV$ w-generates $VC$ on Gauss diagrams, without worrying about signs:

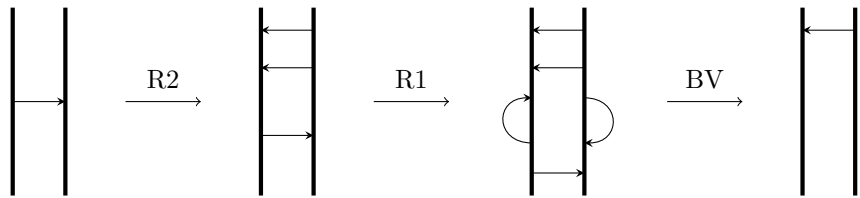
\begin{center}
\begin{tikzpicture}
\draw [ultra thick] (0,0) -- ++(0,2.5) ;
\draw [ultra thick] (0.7,0) -- ++(0,2.5) ;
\draw [>=stealth,->] (0,1.25) -- ++(0.7,0) ;

\begin{scope}[xshift=1.5cm,yshift=1.25cm]
\draw [->] (0,0) -- ++(1,0) ;
\draw (0.5,0.3) node{R2} ;
\end{scope}

\begin{scope}[xshift=3.3cm]
\draw [ultra thick] (0,0) -- ++(0,2.5) ;
\draw [ultra thick] (0.7,0) -- ++(0,2.5) ;
\draw [>=stealth,->] (0,0.8) -- ++(0.7,0) ;
\draw [>=stealth,->] (0.7,1.7) -- ++(-0.7,0) ;
\draw [>=stealth,->] (0.7,2.2) -- ++(-0.7,0) ;
\end{scope}

\begin{scope}[xshift=4.8cm,yshift=1.25cm]
\draw [->] (0,0) -- ++(1,0) ;
\draw (0.5,0.3) node{R1} ;
\end{scope}

\begin{scope}[xshift=6.9cm]
\draw [ultra thick] (0,0) -- ++(0,2.5) ;
\draw [ultra thick] (0.7,0) -- ++(0,2.5) ;
\draw [>=stealth,->] (0,0.3) -- ++(0.7,0) ;
\draw [>=stealth,->] (0,0.7) .. controls +(-0.5,0) and +(-0.5,0) .. (0,1.3) ;
\draw [>=stealth,->] (0.7,1.3) .. controls +(0.5,0) and +(0.5,0) .. (0.7,0.7) ;
\draw [>=stealth,->] (0.7,1.7) -- ++(-0.7,0) ;
\draw [>=stealth,->] (0.7,2.2) -- ++(-0.7,0) ;
\end{scope}

\begin{scope}[xshift=8.7cm,yshift=1.25cm]
\draw [->] (0,0) -- ++(1,0) ;
\draw (0.5,0.3) node{BV} ;
\end{scope}

\begin{scope}[xshift=10.5cm]
\draw [ultra thick] (0,0) -- ++(0,2.5) ;
\draw [ultra thick] (0.7,0) -- ++(0,2.5) ;
\draw [>=stealth,->] (0.7,2.2) -- ++(-0.7,0) ;
\end{scope}
\end{tikzpicture}
\end{center}

\demo \\

We have the following classification of the $BV$ move:

\begin{prop}\label{classBV}
The $BV$ move is w-classified by $(\text{vlk}_{i\ast}^{\text{mod}}+\text{vlk}_{\ast i}^{\text{mod}})_{1\leq i\leq n-1}:vSLD_n\to\ZZ_2^{n-1}$.
\end{prop}

\paragraph{Proof:}
This combination of virtual linking numbers gives the parity of the number of classical overcrossings and undercrossings involving each strand (except for the $n^{th}$ one). Note that we can include the self-crossings, since these count as two crossings on one strand. Since the $BV$ move deletes or adds an even number of classical crossings on each strand, the homomorphism defined above is invariant under $BV$. \\

For a given element $z=(z_1,\ldots,z_{n-1})\in\ZZ_2^{n-1}$, we define a Gauss diagram $G_z$ as follows: for $1\leq i\leq n-1$, we put a positive arrow from the $i^{th}$ strand to the $n^{th}$ strand if $z_i=1$, and no arrow if $z_i=0$. These arrows are taken to be horizontal, with the heads on the $n^{th}$ strand arranged in the order given by the index $i$. The figure below shows an example. By construction, we have $\text{vlk}_{i\ast}^{\text{mod}}(G_z)+\text{vlk}_{\ast i}^{\text{mod}}(G_z)=z_i$ for $1\leq i\leq n-1$.

\begin{center}
\begin{tikzpicture}
\draw [ultra thick,>=stealth,->] (0,0) -- (0,2) ;
\draw [ultra thick,>=stealth,->] (0.8,0) -- (0.8,2) ;
\draw [ultra thick,>=stealth,->] (1.6,0) -- (1.6,2) ;
\draw [ultra thick,>=stealth,->] (2.4,0) -- (2.4,2) ;
\draw [ultra thick,>=stealth,->] (3.2,0) -- (3.2,2) ;
\draw [>=stealth,->] (0,0.3) -- (3.2,0.3) ;
\draw [>=stealth,->] (1.6,0.9) -- (3.2,0.9) ;
\draw [>=stealth,->] (2.4,1.5) -- (3.2,1.5) ;
\draw (2.8,0.5) node{$+$} ;
\draw (2.8,1.1) node{$+$} ;
\draw (2.8,1.7) node{$+$} ;
\draw (1.6,-0.5) node{$G_{(1,0,1,1)}$} ;
\end{tikzpicture}
\end{center}

Using Lemma~\ref{BVwgen}, we can show that any Gauss diagram representing a welded string link $L$ with $\text{vlk}_{i\ast}^{\text{mod}}(L)+\text{vlk}_{\ast i}^{\text{mod}}(L)=z_i$ for $1\leq i\leq n-1$ is equivalent up to $BV$ to the diagram $G_z$ defined above. To achieve this, we first use $VC$ moves to transform arrows from $L_i$ to $L_1$ into arrows from $L_1$ to $L_i$ (where $L_k$ denotes the $k^{th}$ strand of $L$). Self-arrows of $L_1$ are deleted using $SV$ moves, which are allowed since $BV\impw F\impw SV$. The following move as indicated in the figure below allows us to turn arrows from $L_1$ to $L_i$ for $1<i<n$ into arrows from $L_1$ to $L_n$, by creating an additional arrow between $L_i$ and $L_n$. Since we have access to $SR$, signs are irrelevant, and are not displayed:

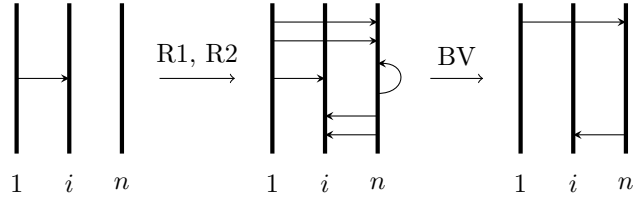
\begin{center}
\begin{tikzpicture}
\draw [ultra thick] (0,0) -- ++(0,2) ;
\draw [ultra thick] (0.7,0) -- ++(0,2) ;
\draw [ultra thick] (1.4,0) -- ++(0,2) ;
\draw [>=stealth,->] (0,1) -- ++(0.7,0) ;
\draw (0,-0.4) node{1} ;
\draw (0.7,-0.4) node{$i$} ;
\draw (1.4,-0.4) node{$n$} ;

\begin{scope}[xshift=1.9cm,yshift=1cm]
\draw [->] (0,0) -- ++(1,0) ;
\draw (0.5,0.3) node{R1, R2} ;
\end{scope}

\begin{scope}[xshift=3.4cm]
\draw [ultra thick] (0,0) -- ++(0,2) ;
\draw [ultra thick] (0.7,0) -- ++(0,2) ;
\draw [ultra thick] (1.4,0) -- ++(0,2) ;
\draw [>=stealth,->] (0,1) -- ++(0.7,0) ;
\draw [>=stealth,->] (0,1.5) -- ++(1.4,0) ;
\draw [>=stealth,->] (0,1.75) -- ++(1.4,0) ;
\draw [>=stealth,->] (1.4,0.8) .. controls +(0.4,0) and +(0.4,0) .. (1.4,1.2) ;
\draw [>=stealth,->] (1.4,0.5) -- ++(-0.7,0) ;
\draw [>=stealth,->] (1.4,0.25) -- ++(-0.7,0) ;
\draw (0,-0.4) node{1} ;
\draw (0.7,-0.4) node{$i$} ;
\draw (1.4,-0.4) node{$n$} ;
\end{scope}

\begin{scope}[xshift=5.5cm,yshift=1cm]
\draw [->] (0,0) -- ++(0.7,0) ;
\draw (0.35,0.3) node{BV} ;
\end{scope}

\begin{scope}[xshift=6.7cm]
\draw [ultra thick] (0,0) -- ++(0,2) ;
\draw [ultra thick] (0.7,0) -- ++(0,2) ;
\draw [ultra thick] (1.4,0) -- ++(0,2) ;
\draw [>=stealth,->] (0,1.75) -- ++(1.4,0) ;
\draw [>=stealth,->] (1.4,0.25) -- ++(-0.7,0) ;
\draw (0,-0.4) node{1} ;
\draw (0.7,-0.4) node{$i$} ;
\draw (1.4,-0.4) node{$n$} ;
\end{scope}
\end{tikzpicture}
\end{center}

Using $F$, $UC$ and $OC$ moves, we can regroup the arrows from $L_1$ to $L_n$, then asign them alternating signs using $SR$ in order to delete them by pairs using $R2$. After these steps, there is at most one arrow attached to $L_1$, which points to $L_n$, and can be made positive by $SR$. Moreover, its head can be positionned at the bottom of $L_n$. \\

We then iterate this process on the diagram obtained by ignoring the first strand. In the end, we are left with a diagram of the form $G_{z'}$. Since the process described above does not change the invariant $(\text{vlk}_{i\ast}^{\text{mod}}+\text{vlk}_{\ast i}^{\text{mod}})_{1\leq i\leq n-1}$, we have $z=z'$, which completes the proof.
\demo \\

The $BV$ move is a \textquotedblleft restricted ribbon residue\textquotedblright of the $BP$ move in the following sense:

\begin{theo}\label{BVpsresBP}
We have $\text{Spun}^{\ast}(BP)|_{2\!-\!r\SL}=\text{Tube}(BV)$.
\end{theo}

\paragraph{Proof:}
Using the geometric interpretation of the linking numbers $\text{LK}_{ij}$ given in the proof of Proposition~\ref{LKDelta}, it can be checked that for $1\leq i\leq n$, the homomorphism $\text{LK}_{i\ast}^{\text{mod}}+\text{LK}_{\ast i}^{\text{mod}}:\dl{2}\to\ZZ_2$ is invariant under a Spun$^{\ast}$($BP$) move. Indeed, the contribution to this combination of linking numbers of the lines of double points inside the torus where the move occurs is always zero before and after the move (note that it is not the case for the Spun($BP$) move used in the proof of Theorem~\ref{VresBP}). Using the classification of the $BV$ move given in Proposition~\ref{classBV}, we can conclude that $\text{Spun}^{\ast}(BP)|_{2\!-\!r\SL}\subset\text{Tube}(BV)$ with the same reasoning as in Theorem~\ref{FresDelta}. \\

For the other inclusion, we only need to check that a Tube($BV$) move can be performed using Spun$^{\ast}$($BP$), which is illustrated on Figure 13. As in the proof of Theorem~\ref{VresBP}, the orientation of the horizontal strands is irrelevant, but the orientation of the vertical strands matter, with the same step added in between pictures 3 and 4 if they have the same orientation. \\

The move used to go from picture 4 to 5 is indeed a Spun$^{\ast}$($BP$) move: the preimages of the lines of double points on the horizontal annuli are essential, and the ones on the vertical annuli are contractible.
\demo \\

As opposed to the ribbon residue $V$ of $BP$, the $BV$ move does not identify all welded string links on $n$ strands. However, because of the symmetry of linking numbers on $SLD_n$, its classifying invariant ${(\text{vlk}_{i\ast}^{\text{mod}}+\text{vlk}_{\ast i}^{\text{mod}})_{1\leq i\leq n-1}}$ vanishes on classical diagrams. Hence the $BV$ move still identifies all \emph{classical} string links on $n$ strands, and as a result this restricted ribbon residue is still not a welded extension of $BP$. \\

The $BP$ move illustrates the case of a classical move where no $SV$ w-genereating welded extension can be obtained by considering the action of its Spun on knotted surfaces, even after introducing some restrictions to compensate an artefact created by the Spun of the move.

\begin{figure}[p]
\centering
\begin{tikzpicture}
\draw [white] (0,0) -- (17,0) -- (17,-21) -- (0,-21) -- cycle ;

\begin{scope}[xshift=1cm,yshift=-3cm,scale=0.75]
\draw [thick,>=stealth,->] (0,0.5) -- (2,0.5) ;
\fill [white] (0.5,0.5) circle (0.15) ;
\fill [white] (1.5,0.5) circle (0.15) ;
\draw [thick,red,>=stealth,->] (0.5,0.75) -- (0.5,0) ;
\draw [thick,red] (0.5,0.75) .. controls +(0,0.75) and +(0,-0.75) .. (1.5,2) ;
\draw [thick,blue,>=stealth,->] (1.5,0.75) -- (1.5,0) ;
\draw [thick,blue] (1.5,0.75) .. controls +(0,0.75) and +(0,-0.75) .. (0.5,2) ;
\draw (1,1.375) circle (0.15) ;
\end{scope}

\begin{scope}[xshift=3.8cm,yshift=-2.25cm,scale=0.75]
\draw [thick,->] (0,0) -- ++(1,0) ;
\draw (0.5,0.5) node{Tube} ;
\end{scope}

\begin{scope}[xshift=6.25cm,yshift=-4.5cm,scale=0.75]
\draw [red] (1,0) -- (1,1) ;
\draw [red] (1,1.5) -- (1,6) ;
\draw [red] (2,0) -- (2,6) ;
\draw [red] (1.25,1.25) ellipse (0.1 and 0.4) ;
\fill [white] (1.1,1) -- (1.1,1.5) -- (1.25,1.5) -- (1.25,1) -- cycle ;
\draw [blue] (3.8,1.25) ellipse (0.1 and 0.4) ;
\fill [white] (3.65,1) -- (3.65,1.5) -- (3.8,1.5) -- (3.8,1) -- cycle ;
\draw (0,1.5) -- (1.32,1.5) ;
\draw [dashed] (1.32,1.5) -- (1.8,1.5) ;
\draw (2.2,1.5) -- (3.87,1.5) ;
\draw [dashed] (3.87,1.5) -- (4.3,1.5) ;
\draw (4.8,1.5) -- (6,1.5) ;
\draw (0,1) -- (1.32,1) ;
\draw [dashed] (1.32,1) -- (1.8,1) ;
\draw (2.2,1) -- (3.87,1) ;
\draw [dashed] (3.87,1) -- (4.3,1) ;
\draw (4.8,1) -- (6,1) ;
\draw [dashed] (4.3,1.25) ellipse (0.0625 and 0.25) ;
\draw (4.8,1.25) ellipse (0.0625 and 0.25) ;
\fill [white] (4.6,1.4) -- (4.6,0.9) -- (5.1,0.9) -- cycle ;
\draw [blue] (0,4.5) -- (1,3.5) ;
\draw [blue] (2,2.5) -- (3,1.5) ;
\draw [blue] (3.5,1) -- (4.5,0) ;
\draw [blue] (0,6) -- (1,5) ;
\draw [blue] (2,4) -- (6,0) ;
\draw [dashed] (1.8,1.25) ellipse (0.0625 and 0.25) ;
\draw (2.2,1.25) ellipse (0.0625 and 0.25) ;
\end{scope}

\begin{scope}[xshift=12.5cm,yshift=-4.5cm,scale=0.75]
\draw [red] (1,0) -- (1,1) ;
\draw [red] (1,1.5) -- (1,6) ;
\draw [red] (2,0) -- (2,6) ;
\draw [red] (1.25,1.25) ellipse (0.1 and 0.4) ;
\fill [white] (1.1,1) -- (1.1,1.5) -- (1.25,1.5) -- (1.25,1) -- cycle ;
\draw [blue] (3.8,1.25) ellipse (0.1 and 0.4) ;
\fill [white] (3.65,1) -- (3.65,1.5) -- (3.8,1.5) -- (3.8,1) -- cycle ;
\draw (0,1.5) -- (1.32,1.5) ;
\draw [dashed] (1.32,1.5) -- (1.8,1.5) ;
\draw (2.2,1.5) -- (3.87,1.5) ;
\draw [dashed] (3.87,1.5) -- (4.3,1.5) ;
\draw (4.8,1.5) -- (6,1.5) ;
\draw (0,1) -- (1.32,1) ;
\draw [dashed] (1.32,1) -- (1.8,1) ;
\draw (2.2,1) -- (3.87,1) ;
\draw [dashed] (3.87,1) -- (4.3,1) ;
\draw (4.8,1) -- (6,1) ;
\draw [dashed] (4.3,1.25) ellipse (0.0625 and 0.25) ;
\draw (4.8,1.25) ellipse (0.0625 and 0.25) ;
\fill [white] (4.6,1.4) -- (4.6,0.9) -- (5.1,0.9) -- cycle ;
\draw [blue] (0,4.5) -- (0.8,3.7) ;
\draw [blue] (0,6) -- (0.8,5.2) ;
\draw [blue] (2.2,3.05) ellipse (0.12 and 0.75) ;
\draw [blue,dashed] (1.2,4.05) ellipse (0.12 and 0.75) ;
\draw [blue,dashed] (1.2,3.3) -- (1.8,2.7) ;
\draw [blue,dashed] (1.2,4.8) -- (1.8,4.2) ;
\draw [blue,dashed] (1.8,3.45) ellipse (0.12 and 0.75) ;
\begin{scope}
\clip (0.8,3.7) -- (0.8,5.2) -- (0.6,5.2) -- (0.6,3.7) -- cycle ;
\draw [blue,dashed] (0.8,4.45) ellipse (0.12 and 0.75) ;
\end{scope}
\begin{scope}
\clip (0.8,3.7) -- (0.8,5.2) -- (1,5.2) -- (1,3.7) -- cycle ;
\draw [blue] (0.8,4.45) ellipse (0.12 and 0.75) ;
\end{scope}
\draw [blue] (2.2,2.3) -- (3,1.5) ;
\draw [blue] (2.2,3.8) -- (6,0) ;
\draw [blue] (3.5,1) -- (4.5,0) ;
\draw [dashed] (1.8,1.25) ellipse (0.0625 and 0.25) ;
\draw (2.2,1.25) ellipse (0.0625 and 0.25) ;
\end{scope}

\begin{scope}[xshift=0cm,yshift=-10.5cm,scale=0.75]
\draw [blue] (0,4.5) -- (4.5,0) ;
\fill [white] (1,3.5) circle (0.15) ;
\fill [white] (2,2.5) circle (0.15) ;
\fill [white] (3,1.5) circle (0.15) ;
\fill [white] (3.5,1) circle (0.15) ;
\draw (1.5,1.5) -- (6,1.5) ;
\draw (1.5,1) -- (6,1) ;
\fill [white] (2,1.5) circle (0.15) ;
\fill [white] (2,1) circle (0.15) ;
\fill [white] (4.5,1.5) circle (0.15) ;
\fill [white] (5,1) circle (0.15) ;
\draw [blue] (0,6) -- (6,0) ;
\fill [white] (1,5) circle (0.15) ;
\fill [white] (2,4) circle (0.15) ;
\draw [red] (1,0) -- (1,6) ;
\fill [white] (1,1.5) circle (0.15) ;
\fill [white] (1,1) circle (0.15) ;
\draw (0,1.5) -- (1.5,1.5) ;
\draw (0,1) -- (1.5,1) ;
\draw [red] (2,0) -- (2,6) ;
\end{scope}

\begin{scope}[xshift=6.25cm,yshift=-10.5cm,scale=0.75]
\draw [blue] (0,4.5) -- (4.5,0) ;
\fill [white] (1,3.5) circle (0.15) ;
\fill [white] (1.5,3) circle (0.15) ;
\fill [white] (2,2.5) circle (0.15) ;
\fill [white] (3.5,1) circle (0.15) ;
\draw (2.25,3.5) -- (6,3.5) ;
\draw (2.25,3) -- (6,3) ;
\fill [white] (2.5,3.5) circle (0.15) ;
\fill [white] (3.5,3.5) circle (0.15) ;
\fill [white] (3,3) circle (0.15) ;
\fill [white] (3.5,3) circle (0.15) ;
\draw [blue] (0,6) -- (6,0) ;
\fill [white] (2,4) circle (0.15) ;
\fill [white] (3.5,2.5) circle (0.15) ;
\draw [red] (2,0) -- (2,6) ;
\fill [white] (2,3.5) circle (0.15) ;
\fill [white] (2,3) circle (0.15) ;
\draw (0,3.5) -- (2.25,3.5) ;
\draw (0,3) -- (2.25,3) ;
\draw [red] (3.5,0) -- (3.5,6) ;
\end{scope}

\begin{scope}[xshift=13cm,yshift=-10.5cm,scale=0.75]
\draw [gray,dashed] (1.5,4.8) -- (2.7,4.8) -- (2.7,2.7) -- (1.5,2.7) -- cycle ;
\draw [blue] (0,5.5) .. controls +(0.4,-0.4) and +(0,0.6) .. (1,4) ;
\draw [blue] (1,4) -- (1,3.5) ;
\draw [blue] (1,3.5) .. controls +(0,-1.5) and +(-0.8,0.4) .. (2,1.25) ;
\draw [blue] (2,1.25) -- (4.5,0) ;
\fill [white] (1,4) circle (0.15) ;
\fill [white] (1,3.5) circle (0.15) ;
\fill [white] (2,1.25) circle (0.15) ;
\fill [white] (3,0.75) circle (0.15) ;
\draw (2.2,4) -- (4.5,4) ;
\draw (2.2,3.5) -- (4.5,3.5) ;
\fill [white] (2.4,4) circle (0.15) ;
\fill [white] (3,4) circle (0.15) ;
\fill [white] (2.4,3.5) circle (0.15) ;
\fill [white] (3,3.5) circle (0.15) ;
\draw [blue] (0.5,6) -- (2,4.5) ;
\draw [blue] (2,4.5) .. controls +(0.6,-0.6) and +(0.6,0.6) .. (2,3) ;
\draw [blue] (2,3) .. controls +(-0.3,-0.3) and +(-0.4,0.2) .. (2,2) ;
\draw [blue] (2,2) -- (4.5,0.75) ;
\fill [white] (2,4.5) circle (0.15) ;
\fill [white] (2,3) circle (0.15) ;
\fill [white] (2,2) circle (0.15) ;
\fill [white] (3,1.5) circle (0.15) ;
\draw [red] (2,0) -- (2,6) ;
\fill [white] (2,4) circle (0.15) ;
\fill [white] (2,3.5) circle (0.15) ;
\draw (0,4) -- (2.2,4) ;
\draw (0,3.5) -- (2.2,3.5) ;
\draw [red] (3,0) -- (3,6) ;
\end{scope}

\begin{scope}[xshift=0.45cm,yshift=-15.1cm,scale=0.8]
\draw [gray,dashed] (0.5,0) -- (0.5,4) -- (4,4) -- (4,0) -- cycle ;
\draw [gray,dashed] (0,2) -- (4.5,2) ;
\draw [gray] (0.1,2.125) .. controls +(-0.05,0.15) and +(0,0.25) .. (-0.15,2) ;
\draw [gray,>=stealth,->] (-0.15,2) .. controls +(0,-0.25) and +(-0.05,-0.15) .. (0.1,1.875) ;
\draw (2.3,2.75) -- (4,2.75) ;
\draw (2.3,1.25) -- (4,1.25) ;
\fill [white] (2.68,2.75) circle (0.15) ;
\fill [white] (2.68,1.25) circle (0.15) ;
\draw [blue] (1,4) .. controls +(0.5,-0.5) and +(0,1) .. (3,2) ;
\draw [blue] (1,0) .. controls +(0.5,0.5) and +(0,-1) .. (3,2) ;
\fill [white] (2,3.32) circle (0.15) ;
\fill [white] (2,0.68) circle (0.15) ;
\draw [red] (2,0) -- (2,4) ;
\fill [white] (2,2.75) circle (0.15) ;
\fill [white] (2,1.25) circle (0.15) ;
\draw (0.5,2.75) -- (2.3,2.75) ;
\draw (0.5,1.25) -- (2.3,1.25) ;
\end{scope}

\begin{scope}[xshift=4.85cm,yshift=-13.5cm,scale=0.8]
\draw [thick,->] (0,0) -- (1,0) ;
\draw (0.5,0.5) node{Spun($\Delta$)} ;
\end{scope}

\begin{scope}[xshift=6.7cm,yshift=-15.1cm,scale=0.8]
\draw [gray,dashed] (0.5,0) -- (0.5,4) -- (4,4) -- (4,0) -- cycle ;
\draw [gray,dashed] (0,2) -- (4.5,2) ;
\draw [gray] (0.1,2.125) .. controls +(-0.05,0.15) and +(0,0.25) .. (-0.15,2) ;
\draw [gray,>=stealth,->] (-0.15,2) .. controls +(0,-0.25) and +(-0.05,-0.15) .. (0.1,1.875) ;
\draw (0.5,3) -- (2,3) ;
\draw (0.5,1) -- (2,1) ;
\fill [white] (1.44,3) circle (0.15) ;
\fill [white] (1.44,1) circle (0.15) ;
\draw [blue] (1,4) .. controls +(0,-1.5) and +(0,0.5) .. (3,2) ;
\draw [blue] (1,0) .. controls +(0,1.5) and +(0,-0.5) .. (3,2) ;
\fill [white] (2.25,2.51) circle (0.15) ;
\fill [white] (2.25,1.49) circle (0.15) ;
\draw [red] (2.25,0) -- (2.25,4) ;
\fill [white] (2.25,3) circle (0.15) ;
\fill [white] (2.25,1) circle (0.15) ;
\draw (2,3) -- (4,3) ;
\draw (2,1) -- (4,1) ;
\end{scope}

\begin{scope}[xshift=12.95cm,yshift=-15.1cm,scale=0.8]
\draw [gray,dashed] (0.5,0) -- (0.5,4) -- (4,4) -- (4,0) -- cycle ;
\draw [gray,dashed] (0,2) -- (4.5,2) ;
\draw [gray] (0.1,2.125) .. controls +(-0.05,0.15) and +(0,0.25) .. (-0.15,2) ;
\draw [gray,>=stealth,->] (-0.15,2) .. controls +(0,-0.25) and +(-0.05,-0.15) .. (0.1,1.875) ;
\draw (0.5,3) -- (2,3) ;
\draw (0.5,1) -- (2,1) ;
\fill [white] (1.25,3) circle (0.15) ;
\fill [white] (1.25,1) circle (0.15) ;
\draw [blue] (1.25,4) -- (1.25,0) ;
\draw [red] (2.25,0) -- (2.25,4) ;
\fill [white] (2.25,3) circle (0.15) ;
\fill [white] (2.25,1) circle (0.15) ;
\draw (2,3) -- (4,3) ;
\draw (2,1) -- (4,1) ;
\end{scope}

\begin{scope}[xshift=0.5cm,yshift=-21cm,scale=0.75]
\draw [gray,dashed] (1.35,4.8) -- (2.5,4.8) -- (2.5,2.7) -- (1.35,2.7) -- cycle ;
\draw [blue] (0,5.5) .. controls +(0.4,-0.4) and +(0,0.6) .. (1,4) ;
\draw [blue] (1,4) -- (1,3.5) ;
\draw [blue] (1,3.5) .. controls +(0,-1.5) and +(-0.8,0.4) .. (2,1.25) ;
\draw [blue] (2,1.25) -- (4.5,0) ;
\fill [white] (1,4) circle (0.15) ;
\fill [white] (1,3.5) circle (0.15) ;
\fill [white] (2,1.25) circle (0.15) ;
\fill [white] (3,0.75) circle (0.15) ;
\draw (1.3,4) -- (1.85,4) ;
\draw (1.3,3.5) -- (1.85,3.5) ;
\fill [white] (1.7,4) circle (0.15) ;
\fill [white] (1.7,3.5) circle (0.15) ;
\draw [blue] (0.5,6) .. controls +(0.6,-0.6) and +(0,1) .. (1.7,4) ;
\draw [blue] (1.7,4) -- (1.7,3.5) ;
\draw [blue] (1.7,3.5) .. controls +(0,-0.5) and +(-0.4,0.2) .. (2,2) ;
\draw [blue] (2,2) -- (4.5,0.75) ;
\fill [white] (2,2) circle (0.15) ;
\fill [white] (3,1.5) circle (0.15) ;
\draw [red] (2,0) -- (2,6) ;
\fill [white] (2,4) circle (0.15) ;
\fill [white] (2,3.5) circle (0.15) ;
\draw (1.85,4) -- (4.5,4) ;
\draw (1.85,3.5) -- (4.5,3.5) ;
\fill [white] (3,4) circle (0.15) ;
\fill [white] (3,3.5) circle (0.15) ;
\draw (0,4) -- (1.3,4) ;
\draw (0,3.5) -- (1.3,3.5) ;
\draw [red] (3,0) -- (3,6) ;
\end{scope}

\begin{scope}[xshift=6.25cm,yshift=-21cm,scale=0.75]
\draw [red] (4,0) -- (4,4.5) ;
\draw [red] (4,5) -- (4,6) ;
\draw [red] (5,0) -- (5,6) ;
\draw [red] (4.25,4.75) ellipse (0.1 and 0.4) ;
\fill [white] (4.1,4.5) -- (4.1,5) -- (4.25,5) -- (4.25,4.5) -- cycle ;
\draw [blue] (1.8,4.75) ellipse (0.1 and 0.4) ;
\fill [white] (1.65,4.5) -- (1.65,5) -- (1.8,5) -- (1.8,4.5) -- cycle ;
\draw (0,5) -- (1.87,5) ;
\draw (0,4.5) -- (1.87,4.5) ;
\draw [dashed] (1.87,5) -- (2.3,5) ;
\draw [dashed] (1.87,4.5) -- (2.3,4.5) ;
\draw [dashed] (2.3,4.75) ellipse (0.0625 and 0.25) ;
\draw (2.8,4.75) ellipse (0.0625 and 0.25) ;
\draw (2.8,5) -- (4.32,5) ;
\draw (2.8,4.5) -- (4.32,4.5) ;
\fill [white] (2.6,4.9) -- (2.6,4.4) -- (3.1,4.4) -- cycle ;
\draw [dashed] (4.32,5) -- (4.8,5) ;
\draw [dashed] (4.32,4.5) -- (4.8,4.5) ;
\draw [dashed] (4.8,4.75) ellipse (0.0625 and 0.25) ;
\draw (5.2,4.75) ellipse (0.0625 and 0.25) ;
\draw (5.2,5) -- (6,5) ;
\draw (5.2,4.5) -- (6,4.5) ;
\draw [blue] (1.5,6) -- (4,3.5) ;
\draw [blue] (5,2.5) -- (6,1.5) ;
\draw [blue] (0,6) -- (1,5) ;
\draw [blue] (1.5,4.5) -- (4,2) ;
\draw [blue] (5,1) -- (6,0) ;
\end{scope}

\begin{scope}[xshift=12.45cm,yshift=-18.75cm,scale=0.75]
\draw [thick,->] (1,0) -- (0,0) ;
\draw (0.5,0.5) node{Tube} ;
\end{scope}

\begin{scope}[xshift=14.5cm,yshift=-19.5cm,scale=0.75]
\draw [thick,>=stealth,->] (0,1.5) -- (2,1.5) ;
\fill [white] (0.5,1.5) circle (0.15) ;
\fill [white] (1.5,1.5) circle (0.15) ;
\draw [thick,red] (1.5,1.25) -- (1.5,2) ;
\draw [thick,red,>=stealth,->] (1.5,1.25) .. controls +(0,-0.75) and +(0,0.75) .. (0.5,0) ;
\draw [thick,blue] (0.5,1.25) -- (0.5,2) ;
\draw [thick,blue,>=stealth,->] (0.5,1.25) .. controls +(0,-0.75) and +(0,0.75) .. (1.5,0) ;
\draw (1,0.625) circle (0.15) ;
\end{scope}
\end{tikzpicture}
\caption{Performing Tube($UC$) using Spun($\Delta$)}
\end{figure}
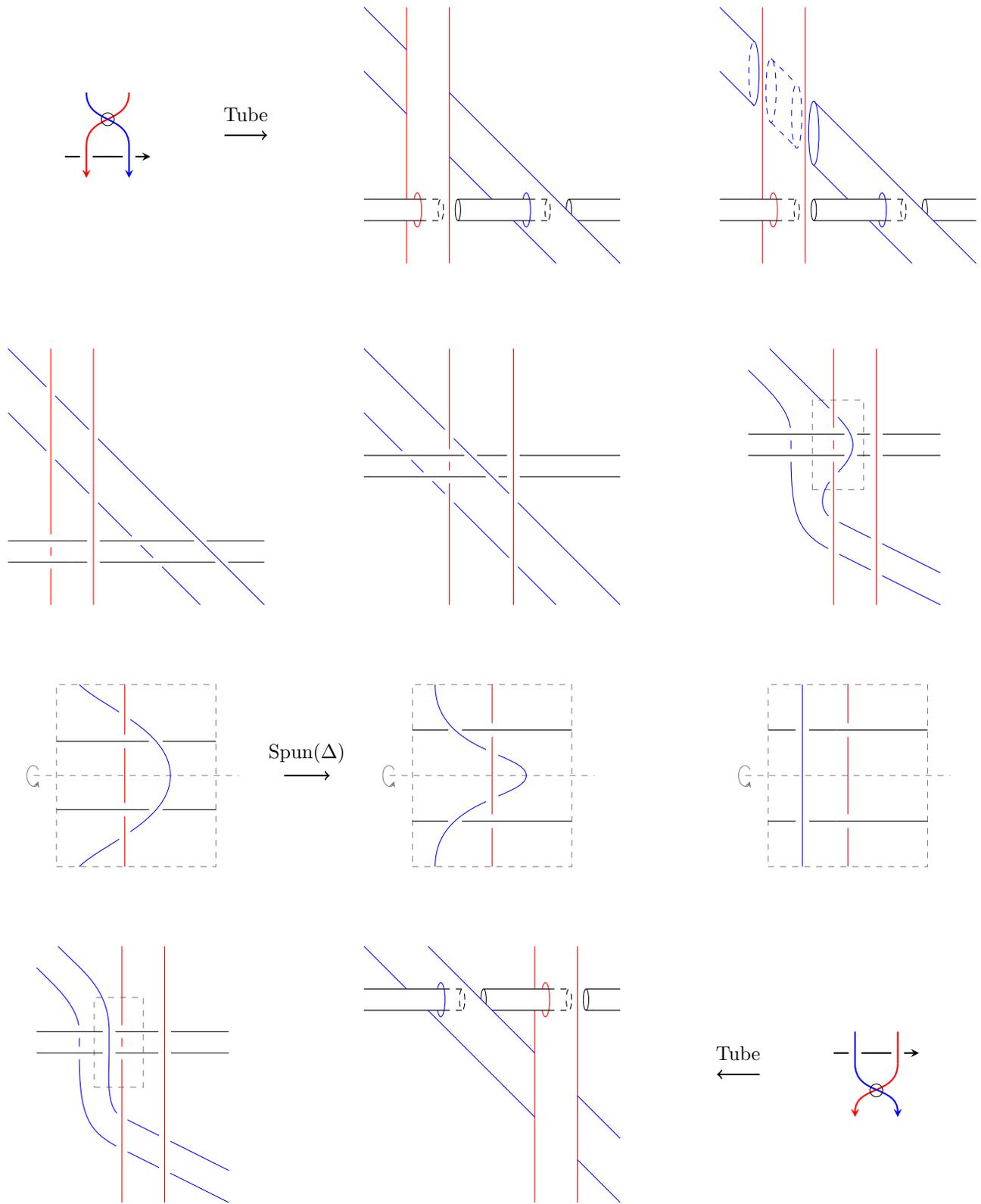

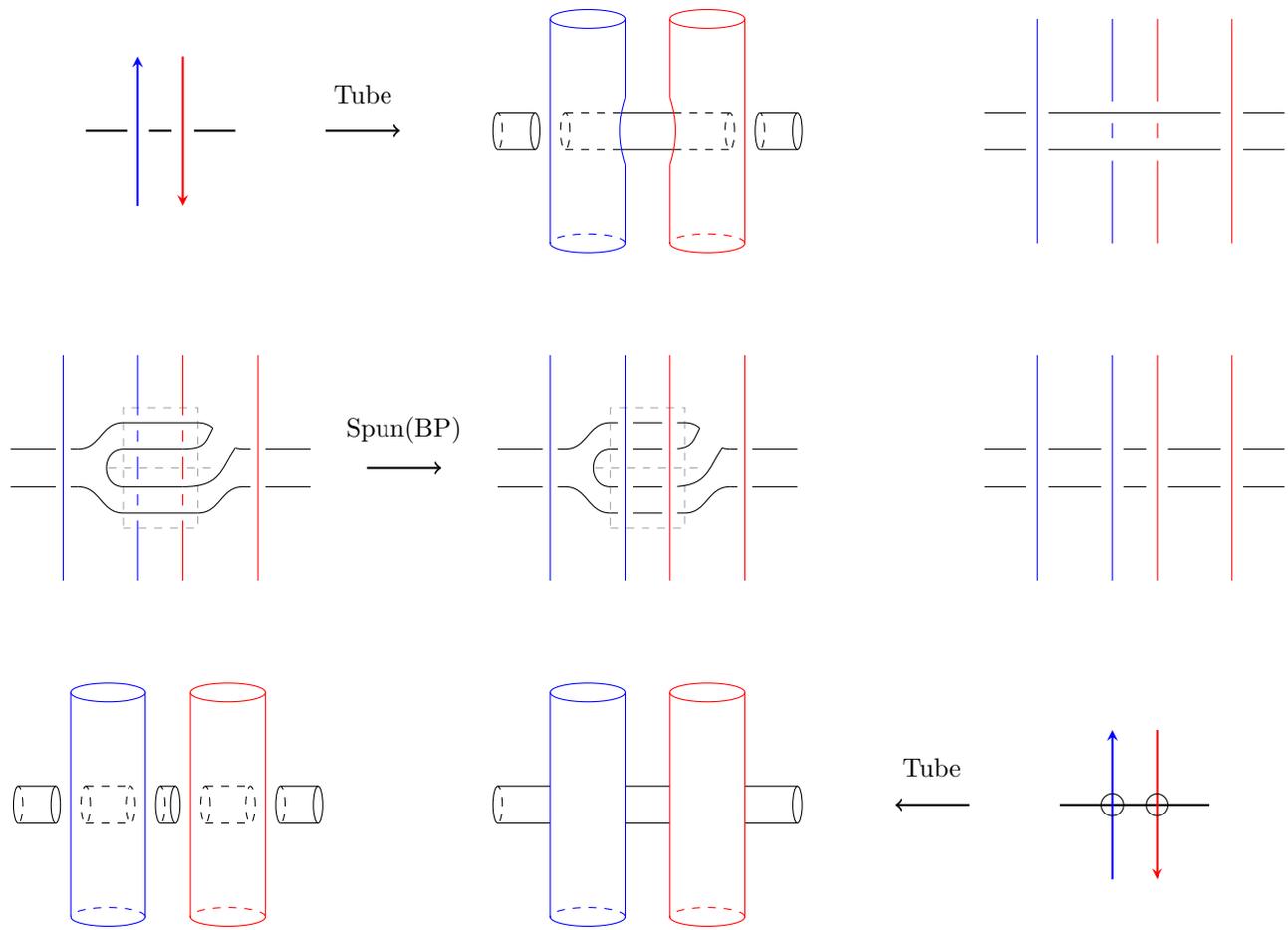
\begin{figure}[p]
\centering
\begin{tikzpicture}
\draw [white] (0,0) -- (17,0) -- (17,-13) -- (0,-13) -- cycle ;

\begin{scope}[xshift=2cm,yshift=-2cm]
\draw [thick] (-1,0) -- (1,0) ;
\fill [white] (-0.3,0) circle (0.15) ;
\fill [white] (0.3,0) circle (0.15) ;
\draw [thick,blue,>=stealth,->] (-0.3,-1) -- (-0.3,1) ;
\draw [thick,red,>=stealth,->] (0.3,1) -- (0.3,-1) ;
\end{scope}

\begin{scope}[xshift=4.7cm,yshift=-2cm]
\draw [thick,->] (-0.5,0) -- (0.5,0) ;
\draw (0,0.5) node{Tube} ;
\end{scope}

\begin{scope}[xshift=8.5cm,yshift=-2cm]
\begin{scope}
\clip (-2.5,-1) -- (-2,-1.5) -- (-2,1.5) -- (-2.5,1.5) -- cycle ;
\draw (-2,0) ellipse (0.0625 and 0.25) ;
\end{scope}
\begin{scope}
\clip (-1.5,-1) -- (-2,-1.5) -- (-2,1.5) -- (-1.5,1.5) -- cycle ;
\draw [dashed] (-2,0) ellipse (0.0625 and 0.25) ;
\end{scope}
\draw (-2,0.25) -- (-1.5,0.25) ;
\draw (-2,-0.25) -- (-1.5,-0.25) ;
\draw (-1.5,0) ellipse (0.0625 and 0.25) ;
\draw [dashed] (-1.1,0) ellipse (0.0625 and 0.25) ;
\draw [dashed] (-1.1,0.25) -- (-0.35,0.25) ;
\draw [dashed] (-1.1,-0.25) -- (-0.35,-0.25) ;
\draw (-0.35,0.25) -- (0.35,0.25) ;
\draw (-0.35,-0.25) -- (0.35,-0.25) ;
\draw [dashed] (0.35,0.25) -- (1.1,0.25) ;
\draw [dashed] (0.35,-0.25) -- (1.1,-0.25) ;
\draw [dashed] (1.1,0) ellipse (0.0625 and 0.25) ;
\begin{scope}
\clip (1,-1) -- (1.5,-1.5) -- (1.5,1.5) -- (1,1.5) -- cycle ;
\draw (1.5,0) ellipse (0.0625 and 0.25) ;
\end{scope}
\begin{scope}
\clip (2,-1) -- (1.5,-1.5) -- (1.5,1.5) -- (2,1.5) -- cycle ;
\draw [dashed] (1.5,0) ellipse (0.0625 and 0.25) ;
\end{scope}
\draw (1.5,0.25) -- (2,0.25) ;
\draw (1.5,-0.25) -- (2,-0.25) ;
\draw (2,0) ellipse (0.0625 and 0.25) ;
\draw [blue] (-0.8,1.5) ellipse (0.5 and 0.125) ;
\draw [blue] (-1.3,-1.5) -- (-1.3,1.5) ;
\draw [blue] (-0.3,-1.5) -- (-0.3,-0.45) ;
\draw [blue] (-0.3,-0.45) .. controls +(-0.1,0.3) and +(-0.1,-0.3) .. (-0.3,0.45) ;
\draw [blue] (-0.3,0.45) -- (-0.3,1.5) ;
\draw [red] (0.8,1.5) ellipse (0.5 and 0.125) ;
\draw [red] (1.3,-1.5) -- (1.3,1.5) ;
\draw [red] (0.3,-1.5) -- (0.3,-0.45) ;
\draw [red] (0.3,-0.45) .. controls +(0.1,0.3) and +(0.1,-0.3) .. (0.3,0.45) ;
\draw [red] (0.3,0.45) -- (0.3,1.5) ;
\begin{scope}
\clip (-1.3,-1.5) -- (1.3,-1.5) -- (1.3,-2) -- (-1.3,-2) -- cycle ;
\draw [blue] (-0.8,-1.5) ellipse (0.5 and 0.125) ;
\draw [red] (0.8,-1.5) ellipse (0.5 and 0.125) ;
\end{scope}
\begin{scope}
\clip (-1.3,-1.5) -- (1.3,-1.5) -- (1.3,-1) -- (-1.3,-1) -- cycle ;
\draw [blue,dashed] (-0.8,-1.5) ellipse (0.5 and 0.125) ;
\draw [red,dashed] (0.8,-1.5) ellipse (0.5 and 0.125) ;
\end{scope}
\end{scope}

\begin{scope}[xshift=15cm,yshift=-2cm]
\draw [blue] (-0.3,-1.5) -- (-0.3,1.5) ;
\draw [red] (0.3,-1.5) -- (0.3,1.5) ;
\fill [white] (-0.3,0.25) circle (0.15) ;
\fill [white] (-0.3,-0.25) circle (0.15) ;
\fill [white] (0.3,0.25) circle (0.15) ;
\fill [white] (0.3,-0.25) circle (0.15) ;
\draw (-2,0.25) -- (2,0.25) ;
\draw (-2,-0.25) -- (2,-0.25) ;
\fill [white] (-1.3,0.25) circle (0.15) ;
\fill [white] (-1.3,-0.25) circle (0.15) ;
\fill [white] (1.3,0.25) circle (0.15) ;
\fill [white] (1.3,-0.25) circle (0.15) ;
\draw [blue] (-1.3,-1.5) -- (-1.3,1.5) ;
\draw [red] (1.3,-1.5) -- (1.3,1.5) ;
\end{scope}

\begin{scope}[xshift=2cm,yshift=-6.5cm]
\draw [gray!70,dashed] (-0.5,0.8) -- (0.5,0.8) -- (0.5,-0.8) -- (-0.5,-0.8) -- cycle ;
\draw [gray!70,dashed] (-0.7,0) -- (0.7,0) ;
\draw [blue] (-0.3,-1.5) -- (-0.3,1.5) ;
\draw [red] (0.3,-1.5) -- (0.3,1.5) ;
\fill [white] (-0.3,-0.6) circle (0.1) ;
\fill [white] (-0.3,-0.25) circle (0.1) ;
\fill [white] (-0.3,0.25) circle (0.1) ;
\fill [white] (-0.3,0.6) circle (0.1) ;
\fill [white] (0.3,-0.6) circle (0.1) ;
\fill [white] (0.3,-0.25) circle (0.1) ;
\fill [white] (0.3,0.25) circle (0.1) ;
\fill [white] (0.3,0.6) circle (0.1) ;
\draw (-2,0.25) -- (-1.1,0.25) ;
\draw (-2,-0.25) -- (-1.1,-0.25) ;
\draw (-1.1,0.25) .. controls +(0.3,0) and +(-0.3,0) .. (-0.5,0.6) ;
\draw (-1.1,-0.25) .. controls +(0.3,0) and +(-0.3,0) .. (-0.5,-0.6) ;
\draw (-0.5,0.6) -- (0.5,0.6) ;
\draw (-0.5,-0.6) -- (0.5,-0.6) ;
\draw (0.5,0.6) .. controls +(0.3,0) and +(-0.3,0) .. (1.1,0.25) ;
\fill [white] (0.7,0.6) -- (1,0.6) -- (1,0.2) -- (0.7,0.2) -- cycle ;
\draw (0.7,0.53) .. controls +(-0.1,-0.2) and +(0.3,0) .. (0.3,0.25) ;
\draw (0.3,0.25) -- (-0.5,0.25) ;
\draw (-0.5,0.25) .. controls +(-0.3,0) and +(-0.3,0) .. (-0.5,-0.25) ;
\draw (-0.5,-0.25) -- (0.3,-0.25) ;
\draw (0.3,-0.25) .. controls +(0.5,0) and +(-0.2,-0.3) .. (1,0.27) ;
\draw (0.5,-0.6) .. controls +(0.3,0) and +(-0.3,0) .. (1.1,-0.25) ;
\draw (1.1,0.25) -- (2,0.25) ;
\draw (1.1,-0.25) -- (2,-0.25) ;
\fill [white] (-1.3,-0.25) circle (0.1) ;
\fill [white] (-1.3,0.25) circle (0.1) ;
\fill [white] (1.3,-0.25) circle (0.1) ;
\fill [white] (1.3,0.25) circle (0.1) ;
\draw [blue] (-1.3,-1.5) -- (-1.3,1.5) ;
\draw [red] (1.3,-1.5) -- (1.3,1.5) ;
\end{scope}

\begin{scope}[xshift=5.25cm,yshift=-6.5cm]
\draw [thick,->] (-0.5,0) -- (0.5,0) ;
\draw (0,0.5) node{Spun(BP)} ;
\end{scope}

\begin{scope}[xshift=8.5cm,yshift=-6.5cm]
\draw [gray!70,dashed] (-0.5,0.8) -- (0.5,0.8) -- (0.5,-0.8) -- (-0.5,-0.8) -- cycle ;
\draw [gray!70,dashed] (-0.7,0) -- (0.7,0) ;
\draw (-2,0.25) -- (-1.1,0.25) ;
\draw (-2,-0.25) -- (-1.1,-0.25) ;
\draw (-1.1,0.25) .. controls +(0.3,0) and +(-0.3,0) .. (-0.5,0.6) ;
\draw (-1.1,-0.25) .. controls +(0.3,0) and +(-0.3,0) .. (-0.5,-0.6) ;
\draw (-0.5,0.6) -- (0.5,0.6) ;
\draw (-0.5,-0.6) -- (0.5,-0.6) ;
\draw (0.5,0.6) .. controls +(0.3,0) and +(-0.3,0) .. (1.1,0.25) ;
\fill [white] (0.7,0.6) -- (1,0.6) -- (1,0.2) -- (0.7,0.2) -- cycle ;
\draw (0.7,0.53) .. controls +(-0.1,-0.2) and +(0.3,0) .. (0.3,0.25) ;
\draw (0.3,0.25) -- (-0.5,0.25) ;
\draw (-0.5,0.25) .. controls +(-0.3,0) and +(-0.3,0) .. (-0.5,-0.25) ;
\draw (-0.5,-0.25) -- (0.3,-0.25) ;
\draw (0.3,-0.25) .. controls +(0.5,0) and +(-0.2,-0.3) .. (1,0.27) ;
\draw (0.5,-0.6) .. controls +(0.3,0) and +(-0.3,0) .. (1.1,-0.25) ;
\draw (1.1,0.25) -- (2,0.25) ;
\draw (1.1,-0.25) -- (2,-0.25) ;
\fill [white] (-0.3,-0.6) circle (0.1) ;
\fill [white] (-0.3,-0.25) circle (0.1) ;
\fill [white] (-0.3,0.25) circle (0.1) ;
\fill [white] (-0.3,0.6) circle (0.1) ;
\fill [white] (0.3,-0.6) circle (0.1) ;
\fill [white] (0.3,-0.25) circle (0.1) ;
\fill [white] (0.3,0.25) circle (0.1) ;
\fill [white] (0.3,0.6) circle (0.1) ;
\fill [white] (-1.3,-0.25) circle (0.1) ;
\fill [white] (-1.3,0.25) circle (0.1) ;
\fill [white] (1.3,-0.25) circle (0.1) ;
\fill [white] (1.3,0.25) circle (0.1) ;
\draw [blue] (-0.3,-1.5) -- (-0.3,1.5) ;
\draw [red] (0.3,-1.5) -- (0.3,1.5) ;
\draw [blue] (-1.3,-1.5) -- (-1.3,1.5) ;
\draw [red] (1.3,-1.5) -- (1.3,1.5) ;
\end{scope}

\begin{scope}[xshift=15cm,yshift=-6.5cm]
\draw (-2,0.25) -- (2,0.25) ;
\draw (-2,-0.25) -- (2,-0.25) ;
\fill [white] (-0.3,0.25) circle (0.15) ;
\fill [white] (-0.3,-0.25) circle (0.15) ;
\fill [white] (0.3,0.25) circle (0.15) ;
\fill [white] (0.3,-0.25) circle (0.15) ;
\fill [white] (-1.3,0.25) circle (0.15) ;
\fill [white] (-1.3,-0.25) circle (0.15) ;
\fill [white] (1.3,0.25) circle (0.15) ;
\fill [white] (1.3,-0.25) circle (0.15) ;
\draw [blue] (-0.3,-1.5) -- (-0.3,1.5) ;
\draw [red] (0.3,-1.5) -- (0.3,1.5) ;
\draw [blue] (-1.3,-1.5) -- (-1.3,1.5) ;
\draw [red] (1.3,-1.5) -- (1.3,1.5) ;
\end{scope}

\begin{scope}[xshift=2.1cm,yshift=-11cm]
\begin{scope}
\clip (-2.1,-1) -- (-2,-1.5) -- (-2,1.5) -- (-2.1,1.5) -- cycle ;
\draw (-2,0) ellipse (0.0625 and 0.25) ;
\end{scope}
\begin{scope}
\clip (-1.9,-1) -- (-2,-1.5) -- (-2,1.5) -- (-1.9,1.5) -- cycle ;
\draw [dashed] (-2,0) ellipse (0.0625 and 0.25) ;
\end{scope}
\draw (-2,0.25) -- (-1.5,0.25) ;
\draw (-2,-0.25) -- (-1.5,-0.25) ;
\draw (-1.5,0) ellipse (0.0625 and 0.25) ;
\draw [dashed] (-1.1,0) ellipse (0.0625 and 0.25) ;
\draw [dashed] (-1.1,0.25) -- (-0.5,0.25) ;
\draw [dashed] (-1.1,-0.25) -- (-0.5,-0.25) ;
\draw [dashed] (-0.5,0) ellipse (0.0625 and 0.25) ;
\begin{scope}
\clip (-1,-1) -- (-0.1,-1.5) -- (-0.1,1.5) -- (-1,1.5) -- cycle ;
\draw (-0.1,0) ellipse (0.0625 and 0.25) ;
\end{scope}
\begin{scope}
\clip (-0.1,-1) -- (0,-1.5) -- (0,1.5) -- (-0.1,1.5) -- cycle ;
\draw [dashed] (-0.1,0) ellipse (0.0625 and 0.25) ;
\end{scope}
\draw (-0.1,0.25) -- (0.1,0.25) ;
\draw (-0.1,-0.25) -- (0.1,-0.25) ;
\draw (0.1,0) ellipse (0.0625 and 0.25) ;
\draw [dashed] (0.5,0) ellipse (0.0625 and 0.25) ;
\draw [dashed] (0.5,0.25) -- (1.1,0.25) ;
\draw [dashed] (0.5,-0.25) -- (1.1,-0.25) ;
\draw [dashed] (1.1,0) ellipse (0.0625 and 0.25) ;
\begin{scope}
\clip (1,-1) -- (1.5,-1.5) -- (1.5,1.5) -- (1,1.5) -- cycle ;
\draw (1.5,0) ellipse (0.0625 and 0.25) ;
\end{scope}
\begin{scope}
\clip (2,-1) -- (1.5,-1.5) -- (1.5,1.5) -- (2,1.5) -- cycle ;
\draw [dashed] (1.5,0) ellipse (0.0625 and 0.25) ;
\end{scope}
\draw (1.5,0.25) -- (2,0.25) ;
\draw (1.5,-0.25) -- (2,-0.25) ;
\draw (2,0) ellipse (0.0625 and 0.25) ;
\draw [blue] (-0.8,1.5) ellipse (0.5 and 0.125) ;
\draw [blue] (-1.3,-1.5) -- (-1.3,1.5) ;
\draw [blue] (-0.3,-1.5) -- (-0.3,1.5) ;
\draw [red] (0.8,1.5) ellipse (0.5 and 0.125) ;
\draw [red] (1.3,-1.5) -- (1.3,1.5) ;
\draw [red] (0.3,-1.5) -- (0.3,1.5) ;
\begin{scope}
\clip (-1.3,-1.5) -- (1.3,-1.5) -- (1.3,-2) -- (-1.3,-2) -- cycle ;
\draw [blue] (-0.8,-1.5) ellipse (0.5 and 0.125) ;
\draw [red] (0.8,-1.5) ellipse (0.5 and 0.125) ;
\end{scope}
\begin{scope}
\clip (-1.3,-1.5) -- (1.3,-1.5) -- (1.3,-1) -- (-1.3,-1) -- cycle ;
\draw [blue,dashed] (-0.8,-1.5) ellipse (0.5 and 0.125) ;
\draw [red,dashed] (0.8,-1.5) ellipse (0.5 and 0.125) ;
\end{scope}
\end{scope}

\begin{scope}[xshift=8.5cm,yshift=-11cm]
\begin{scope}
\clip (-2.5,-1) -- (-2,-1.5) -- (-2,1.5) -- (-2.5,1.5) -- cycle ;
\draw (-2,0) ellipse (0.0625 and 0.25) ;
\end{scope}
\begin{scope}
\clip (-1.5,-1) -- (-2,-1.5) -- (-2,1.5) -- (-1.5,1.5) -- cycle ;
\draw [dashed] (-2,0) ellipse (0.0625 and 0.25) ;
\end{scope}
\draw (-2,0.25) -- (-1.3,0.25) ;
\draw (-2,-0.25) -- (-1.3,-0.25) ;
\draw (-0.3,0.25) -- (0.3,0.25) ;
\draw (-0.3,-0.25) -- (0.3,-0.25) ;
\draw (1.3,0.25) -- (2,0.25) ;
\draw (1.3,-0.25) -- (2,-0.25) ;
\draw (2,0) ellipse (0.0625 and 0.25) ;
\draw [blue] (-0.8,1.5) ellipse (0.5 and 0.125) ;
\draw [blue] (-1.3,-1.5) -- (-1.3,1.5) ;
\draw [blue] (-0.3,-1.5) -- (-0.3,1.5) ;
\draw [red] (0.8,1.5) ellipse (0.5 and 0.125) ;
\draw [red] (1.3,-1.5) -- (1.3,1.5) ;
\draw [red] (0.3,-1.5) -- (0.3,1.5) ;
\begin{scope}
\clip (-1.3,-1.5) -- (1.3,-1.5) -- (1.3,-2) -- (-1.3,-2) -- cycle ;
\draw [blue] (-0.8,-1.5) ellipse (0.5 and 0.125) ;
\draw [red] (0.8,-1.5) ellipse (0.5 and 0.125) ;
\end{scope}
\begin{scope}
\clip (-1.3,-1.5) -- (1.3,-1.5) -- (1.3,-1) -- (-1.3,-1) -- cycle ;
\draw [blue,dashed] (-0.8,-1.5) ellipse (0.5 and 0.125) ;
\draw [red,dashed] (0.8,-1.5) ellipse (0.5 and 0.125) ;
\end{scope}
\end{scope}

\begin{scope}[xshift=12.3cm,yshift=-11cm]
\draw [thick,->] (0.5,0) -- (-0.5,0) ;
\draw (0,0.5) node{Tube} ;
\end{scope}

\begin{scope}[xshift=15cm,yshift=-11cm]
\draw [thick] (-1,0) -- (1,0) ;
\draw [thick,blue,>=stealth,->] (-0.3,-1) -- (-0.3,1) ;
\draw [thick,red,>=stealth,->] (0.3,1) -- (0.3,-1) ;
\draw (-0.3,0) circle (0.15) ;
\draw (0.3,0) circle (0.15) ;
\end{scope}
\end{tikzpicture}
\caption{Performing Tube($DV$) using Spun($BP$)}
\end{figure}

\begin{figure}[p]
\centering
\begin{tikzpicture}
\draw [white] (0,0) -- (17,0) -- (17,-13) -- (0,-13) -- cycle ;

\begin{scope}[xshift=2cm,yshift=-2cm]
\draw [thick,blue] (-1,0.3) -- (1,0.3) ;
\draw [thick,red] (-1,-0.3) -- (1,-0.3) ;
\fill [white] (-0.3,0.3) circle (0.15) ;
\fill [white] (0.3,0.3) circle (0.15) ;
\fill [white] (-0.3,-0.3) circle (0.15) ;
\fill [white] (0.3,-0.3) circle (0.15) ;
\draw [thick,>=stealth,->] (-0.3,-1) -- (-0.3,1) ;
\draw [thick,>=stealth,->] (0.3,1) -- (0.3,-1) ;
\end{scope}

\begin{scope}[xshift=4.7cm,yshift=-2cm]
\draw [thick,->] (-0.5,0) -- (0.5,0) ;
\draw (0,0.5) node{Tube} ;
\end{scope}

\begin{scope}[xshift=8.5cm,yshift=-2cm]
\begin{scope}
\clip (-2.5,-1) -- (-2,-1.5) -- (-2,1.5) -- (-2.5,1.5) -- cycle ;
\draw [blue] (-2,0.75) ellipse (0.0625 and 0.25) ;
\draw [red] (-2,-0.75) ellipse (0.0625 and 0.25) ;
\end{scope}
\begin{scope}
\clip (-1.5,-1) -- (-2,-1.5) -- (-2,1.5) -- (-1.5,1.5) -- cycle ;
\draw [blue,dashed] (-2,0.75) ellipse (0.0625 and 0.25) ;
\draw [red,dashed] (-2,-0.75) ellipse (0.0625 and 0.25) ;
\end{scope}
\draw [blue] (-2,1) -- (-1.5,1) ;
\draw [blue] (-2,0.5) -- (-1.5,0.5) ;
\draw [blue] (-1.5,0.75) ellipse (0.0625 and 0.25) ;
\draw [red] (-2,-1) -- (-1.5,-1) ;
\draw [red] (-2,-0.5) -- (-1.5,-0.5) ;
\draw [red] (-1.5,-0.75) ellipse (0.0625 and 0.25) ;
\draw [blue,dashed] (-1.1,0.75) ellipse (0.0625 and 0.25) ;
\draw [red,dashed] (-1.1,-0.75) ellipse (0.0625 and 0.25) ;
\draw [blue,dashed] (-1.1,1) -- (-0.35,1) ;
\draw [blue,dashed] (-1.1,0.5) -- (-0.35,0.5) ;
\draw [red,dashed] (-1.1,-1) -- (-0.35,-1) ;
\draw [red,dashed] (-1.1,-0.5) -- (-0.35,-0.5) ;
\draw [blue] (-0.35,1) -- (0.35,1) ;
\draw [blue] (-0.35,0.5) -- (0.35,0.5) ;
\draw [red] (-0.35,-1) -- (0.35,-1) ;
\draw [red] (-0.35,-0.5) -- (0.35,-0.5) ;
\draw [blue,dashed] (0.35,1) -- (1.1,1) ;
\draw [blue,dashed] (0.35,0.5) -- (1.1,0.5) ;
\draw [blue,dashed] (1.1,0.75) ellipse (0.0625 and 0.25) ;
\draw [red,dashed] (0.35,-1) -- (1.1,-1) ;
\draw [red,dashed] (0.35,-0.5) -- (1.1,-0.5) ;
\draw [red,dashed] (1.1,-0.75) ellipse (0.0625 and 0.25) ;
\begin{scope}
\clip (1,-1) -- (1.5,-1.5) -- (1.5,1.5) -- (1,1.5) -- cycle ;
\draw [blue] (1.5,0.75) ellipse (0.0625 and 0.25) ;
\draw [red] (1.5,-0.75) ellipse (0.0625 and 0.25) ;
\end{scope}
\begin{scope}
\clip (2,-1) -- (1.5,-1.5) -- (1.5,1.5) -- (2,1.5) -- cycle ;
\draw [blue,dashed] (1.5,0.75) ellipse (0.0625 and 0.25) ;
\draw [red,dashed] (1.5,-0.75) ellipse (0.0625 and 0.25) ;
\end{scope}
\draw [blue] (1.5,1) -- (2,1) ;
\draw [blue] (1.5,0.5) -- (2,0.5) ;
\draw [blue] (2,0.75) ellipse (0.0625 and 0.25) ;
\draw [red] (1.5,-1) -- (2,-1) ;
\draw [red] (1.5,-0.5) -- (2,-0.5) ;
\draw [red] (2,-0.75) ellipse (0.0625 and 0.25) ;
\draw (-0.8,1.5) ellipse (0.5 and 0.125) ;
\draw (-1.3,-1.5) -- (-1.3,1.5) ;
\draw (-0.3,-1.5) -- (-0.3,-1.2) ;
\draw (-0.3,-1.2) .. controls +(-0.1,0.3) and +(-0.1,-0.3) .. (-0.3,-0.3) ;
\draw (-0.3,-0.3) -- (-0.3,0.3) ;
\draw (-0.3,0.3) .. controls +(-0.1,0.3) and +(-0.1,-0.3) .. (-0.3,1.2) ;
\draw (-0.3,1.2) -- (-0.3,1.5) ;
\draw (0.8,1.5) ellipse (0.5 and 0.125) ;
\draw (1.3,-1.5) -- (1.3,1.5) ;
\draw (0.3,-1.5) -- (0.3,-1.2) ;
\draw (0.3,-1.2) .. controls +(0.1,0.3) and +(0.1,-0.3) .. (0.3,-0.3) ;
\draw (0.3,-0.3) -- (0.3,0.3) ;
\draw (0.3,0.3) .. controls +(0.1,0.3) and +(0.1,-0.3) .. (0.3,1.2) ;
\draw (0.3,1.2) -- (0.3,1.5) ;
\begin{scope}
\clip (-1.3,-1.5) -- (1.3,-1.5) -- (1.3,-2) -- (-1.3,-2) -- cycle ;
\draw (-0.8,-1.5) ellipse (0.5 and 0.125) ;
\draw (0.8,-1.5) ellipse (0.5 and 0.125) ;
\end{scope}
\begin{scope}
\clip (-1.3,-1.5) -- (1.3,-1.5) -- (1.3,-1) -- (-1.3,-1) -- cycle ;
\draw [dashed] (-0.8,-1.5) ellipse (0.5 and 0.125) ;
\draw [dashed] (0.8,-1.5) ellipse (0.5 and 0.125) ;
\end{scope}
\end{scope}

\begin{scope}[xshift=15cm,yshift=-2cm]
\draw (-0.3,-1.5) -- (-0.3,1.5) ;
\draw (0.3,-1.5) -- (0.3,1.5) ;
\fill [white] (-0.3,1) circle (0.15) ;
\fill [white] (-0.3,0.5) circle (0.15) ;
\fill [white] (-0.3,-1) circle (0.15) ;
\fill [white] (-0.3,-0.5) circle (0.15) ;
\fill [white] (0.3,1) circle (0.15) ;
\fill [white] (0.3,0.5) circle (0.15) ;
\fill [white] (0.3,-1) circle (0.15) ;
\fill [white] (0.3,-0.5) circle (0.15) ;
\draw [blue] (-2,1) -- (2,1) ;
\draw [blue] (-2,0.5) -- (2,0.5) ;
\draw [red] (-2,-1) -- (2,-1) ;
\draw [red] (-2,-0.5) -- (2,-0.5) ;
\fill [white] (-1.3,1) circle (0.15) ;
\fill [white] (-1.3,0.5) circle (0.15) ;
\fill [white] (-1.3,-1) circle (0.15) ;
\fill [white] (-1.3,-0.5) circle (0.15) ;
\fill [white] (1.3,1) circle (0.15) ;
\fill [white] (1.3,0.5) circle (0.15) ;
\fill [white] (1.3,-1) circle (0.15) ;
\fill [white] (1.3,-0.5) circle (0.15) ;
\draw (-1.3,-1.5) -- (-1.3,1.5) ;
\draw (1.3,-1.5) -- (1.3,1.5) ;
\end{scope}

\begin{scope}[xshift=2cm,yshift=-6.5cm]
\draw [gray!70,dashed] (-0.6,0.65) -- (0.6,0.65) -- (0.6,-0.65) -- (-0.6,-0.65) -- cycle ;
\draw [gray!70,dashed] (-0.8,0) -- (0.8,0) ;
\draw (-0.3,-1.5) -- (-0.3,1.5) ;
\draw (0.3,-1.5) -- (0.3,1.5) ;
\fill [white] (-0.3,0.45) circle (0.1) ;
\fill [white] (-0.3,0.15) circle (0.1) ;
\fill [white] (-0.3,-0.45) circle (0.1) ;
\fill [white] (-0.3,-0.15) circle (0.1) ;
\fill [white] (0.3,0.45) circle (0.1) ;
\fill [white] (0.3,0.15) circle (0.1) ;
\fill [white] (0.3,-0.45) circle (0.1) ;
\fill [white] (0.3,-0.15) circle (0.1) ;
\draw [red] (-2,-1) -- (-1.1,-1) ;
\draw [red] (-2,-0.5) -- (-1.1,-0.5) ;
\draw [red] (-1.1,-1) .. controls +(0.5,0) and +(-0.3,0) .. (-0.5,-0.45) ;
\draw [red] (-1.1,-0.5) .. controls +(0.3,0) and +(-0.5,0) .. (-0.5,0.45) ;
\draw [red] (-0.5,0.45) -- (0.5,0.45) ;
\draw [red] (-0.5,-0.45) -- (0.5,-0.45) ;
\draw [red] (0.5,-0.45) .. controls +(0.3,0) and +(-0.5,0) .. (1.1,-1) ;
\draw [red] (0.5,0.45) .. controls +(0.5,0) and +(-0.3,0) .. (1.1,-0.5) ;
\draw [red] (1.1,-1) -- (2,-1) ;
\draw [red] (1.1,-0.5) -- (2,-0.5) ;
\fill [white] (-0.83,0.15) circle (0.1) ;
\fill [white] (-0.695,0.38) circle (0.1) ;
\fill [white] (0.83,0.15) circle (0.1) ;
\fill [white] (0.695,0.38) circle (0.1) ;
\draw [blue] (-2,1) -- (-1.1,1) ;
\draw [blue] (-2,0.5) -- (-1.1,0.5) ;
\draw [blue] (-1.1,1) .. controls +(0.4,0) and +(-0.3,0) .. (-0.5,0.15) ;
\draw [blue] (-1.1,0.5) .. controls +(0.3,0) and +(-0.4,0) .. (-0.5,-0.15) ;
\draw [blue] (-0.5,0.15) -- (0.5,0.15) ;
\draw [blue] (-0.5,-0.15) -- (0.5,-0.15) ;
\draw [blue] (0.5,0.15) .. controls +(0.3,0) and +(-0.4,0) .. (1.1,1) ;
\draw [blue] (0.5,-0.15) .. controls +(0.4,0) and +(-0.3,0) .. (1.1,0.5) ;
\draw [blue] (1.1,1) -- (2,1) ;
\draw [blue] (1.1,0.5) -- (2,0.5) ;
\fill [white] (-1.3,1) circle (0.15) ;
\fill [white] (-1.3,0.5) circle (0.15) ;
\fill [white] (-1.3,-1) circle (0.15) ;
\fill [white] (-1.3,-0.5) circle (0.15) ;
\fill [white] (1.3,1) circle (0.15) ;
\fill [white] (1.3,0.5) circle (0.15) ;
\fill [white] (1.3,-1) circle (0.15) ;
\fill [white] (1.3,-0.5) circle (0.15) ;
\draw (-1.3,-1.5) -- (-1.3,1.5) ;
\draw (1.3,-1.5) -- (1.3,1.5) ;
\end{scope}

\begin{scope}[xshift=5.25cm,yshift=-6.5cm]
\draw [thick,->] (-0.5,0) -- (0.5,0) ;
\draw (0,0.5) node{Spun$^{\ast}$(BP)} ;
\end{scope}

\begin{scope}[xshift=8.5cm,yshift=-6.5cm]
\draw [gray!70,dashed] (-0.6,0.65) -- (0.6,0.65) -- (0.6,-0.65) -- (-0.6,-0.65) -- cycle ;
\draw [gray!70,dashed] (-0.8,0) -- (0.8,0) ;
\draw [red] (-2,-1) -- (-1.1,-1) ;
\draw [red] (-2,-0.5) -- (-1.1,-0.5) ;
\draw [red] (-1.1,-1) .. controls +(0.5,0) and +(-0.3,0) .. (-0.5,-0.45) ;
\draw [red] (-1.1,-0.5) .. controls +(0.3,0) and +(-0.5,0) .. (-0.5,0.45) ;
\draw [red] (-0.5,0.45) -- (0.5,0.45) ;
\draw [red] (-0.5,-0.45) -- (0.5,-0.45) ;
\draw [red] (0.5,-0.45) .. controls +(0.3,0) and +(-0.5,0) .. (1.1,-1) ;
\draw [red] (0.5,0.45) .. controls +(0.5,0) and +(-0.3,0) .. (1.1,-0.5) ;
\draw [red] (1.1,-1) -- (2,-1) ;
\draw [red] (1.1,-0.5) -- (2,-0.5) ;
\fill [white] (-0.83,0.15) circle (0.1) ;
\fill [white] (-0.695,0.38) circle (0.1) ;
\fill [white] (0.83,0.15) circle (0.1) ;
\fill [white] (0.695,0.38) circle (0.1) ;
\draw [blue] (-2,1) -- (-1.1,1) ;
\draw [blue] (-2,0.5) -- (-1.1,0.5) ;
\draw [blue] (-1.1,1) .. controls +(0.4,0) and +(-0.3,0) .. (-0.5,0.15) ;
\draw [blue] (-1.1,0.5) .. controls +(0.3,0) and +(-0.4,0) .. (-0.5,-0.15) ;
\draw [blue] (-0.5,0.15) -- (0.5,0.15) ;
\draw [blue] (-0.5,-0.15) -- (0.5,-0.15) ;
\draw [blue] (0.5,0.15) .. controls +(0.3,0) and +(-0.4,0) .. (1.1,1) ;
\draw [blue] (0.5,-0.15) .. controls +(0.4,0) and +(-0.3,0) .. (1.1,0.5) ;
\draw [blue] (1.1,1) -- (2,1) ;
\draw [blue] (1.1,0.5) -- (2,0.5) ;
\fill [white] (-0.3,0.45) circle (0.1) ;
\fill [white] (-0.3,0.15) circle (0.1) ;
\fill [white] (-0.3,-0.45) circle (0.1) ;
\fill [white] (-0.3,-0.15) circle (0.1) ;
\fill [white] (0.3,0.45) circle (0.1) ;
\fill [white] (0.3,0.15) circle (0.1) ;
\fill [white] (0.3,-0.45) circle (0.1) ;
\fill [white] (0.3,-0.15) circle (0.1) ;
\fill [white] (-1.3,1) circle (0.15) ;
\fill [white] (-1.3,0.5) circle (0.15) ;
\fill [white] (-1.3,-1) circle (0.15) ;
\fill [white] (-1.3,-0.5) circle (0.15) ;
\fill [white] (1.3,1) circle (0.15) ;
\fill [white] (1.3,0.5) circle (0.15) ;
\fill [white] (1.3,-1) circle (0.15) ;
\fill [white] (1.3,-0.5) circle (0.15) ;
\draw (-0.3,-1.5) -- (-0.3,1.5) ;
\draw (0.3,-1.5) -- (0.3,1.5) ;
\draw (-1.3,-1.5) -- (-1.3,1.5) ;
\draw (1.3,-1.5) -- (1.3,1.5) ;
\end{scope}

\begin{scope}[xshift=15cm,yshift=-6.5cm]
\draw [blue] (-2,1) -- (2,1) ;
\draw [blue] (-2,0.5) -- (2,0.5) ;
\draw [red] (-2,-1) -- (2,-1) ;
\draw [red] (-2,-0.5) -- (2,-0.5) ;
\fill [white] (-0.3,1) circle (0.15) ;
\fill [white] (-0.3,0.5) circle (0.15) ;
\fill [white] (-0.3,-1) circle (0.15) ;
\fill [white] (-0.3,-0.5) circle (0.15) ;
\fill [white] (0.3,1) circle (0.15) ;
\fill [white] (0.3,0.5) circle (0.15) ;
\fill [white] (0.3,-1) circle (0.15) ;
\fill [white] (0.3,-0.5) circle (0.15) ;
\fill [white] (-1.3,1) circle (0.15) ;
\fill [white] (-1.3,0.5) circle (0.15) ;
\fill [white] (-1.3,-1) circle (0.15) ;
\fill [white] (-1.3,-0.5) circle (0.15) ;
\fill [white] (1.3,1) circle (0.15) ;
\fill [white] (1.3,0.5) circle (0.15) ;
\fill [white] (1.3,-1) circle (0.15) ;
\fill [white] (1.3,-0.5) circle (0.15) ;
\draw (-0.3,-1.5) -- (-0.3,1.5) ;
\draw (0.3,-1.5) -- (0.3,1.5) ;
\draw (-1.3,-1.5) -- (-1.3,1.5) ;
\draw (1.3,-1.5) -- (1.3,1.5) ;
\end{scope}

\begin{scope}[xshift=2.1cm,yshift=-11cm]
\begin{scope}
\clip (-2.1,-1) -- (-2,-1.5) -- (-2,1.5) -- (-2.1,1.5) -- cycle ;
\draw [blue] (-2,0.75) ellipse (0.0625 and 0.25) ;
\draw [red] (-2,-0.75) ellipse (0.0625 and 0.25) ;
\end{scope}
\begin{scope}
\clip (-1.9,-1) -- (-2,-1.5) -- (-2,1.5) -- (-1.9,1.5) -- cycle ;
\draw [blue,dashed] (-2,0.75) ellipse (0.0625 and 0.25) ;
\draw [red,dashed] (-2,-0.75) ellipse (0.0625 and 0.25) ;
\end{scope}
\draw [blue] (-2,1) -- (-1.5,1) ;
\draw [blue] (-2,0.5) -- (-1.5,0.5) ;
\draw [blue] (-1.5,0.75) ellipse (0.0625 and 0.25) ;
\draw [red] (-2,-1) -- (-1.5,-1) ;
\draw [red] (-2,-0.5) -- (-1.5,-0.5) ;
\draw [red] (-1.5,-0.75) ellipse (0.0625 and 0.25) ;
\draw [blue,dashed] (-1.1,0.75) ellipse (0.0625 and 0.25) ;
\draw [red,dashed] (-1.1,-0.75) ellipse (0.0625 and 0.25) ;
\draw [blue,dashed] (-1.1,1) -- (-0.5,1) ;
\draw [blue,dashed] (-1.1,0.5) -- (-0.5,0.5) ;
\draw [blue,dashed] (-0.5,0.75) ellipse (0.0625 and 0.25) ;
\draw [red,dashed] (-1.1,-1) -- (-0.5,-1) ;
\draw [red,dashed] (-1.1,-0.5) -- (-0.5,-0.5) ;
\draw [red,dashed] (-0.5,-0.75) ellipse (0.0625 and 0.25) ;
\begin{scope}
\clip (-1,-1) -- (-0.1,-1.5) -- (-0.1,1.5) -- (-1,1.5) -- cycle ;
\draw [blue] (-0.1,0.75) ellipse (0.0625 and 0.25) ;
\draw [red] (-0.1,-0.75) ellipse (0.0625 and 0.25) ;
\end{scope}
\begin{scope}
\clip (-0.1,-1) -- (0,-1.5) -- (0,1.5) -- (-0.1,1.5) -- cycle ;
\draw [blue,dashed] (-0.1,0.75) ellipse (0.0625 and 0.25) ;
\draw [red,dashed] (-0.1,-0.75) ellipse (0.0625 and 0.25) ;
\end{scope}
\draw [blue] (-0.1,1) -- (0.1,1) ;
\draw [blue] (-0.1,0.5) -- (0.1,0.5) ;
\draw [blue] (0.1,0.75) ellipse (0.0625 and 0.25) ;
\draw [red] (-0.1,-1) -- (0.1,-1) ;
\draw [red] (-0.1,-0.5) -- (0.1,-0.5) ;
\draw [red] (0.1,-0.75) ellipse (0.0625 and 0.25) ;
\draw [blue,dashed] (0.5,0.75) ellipse (0.0625 and 0.25) ;
\draw [red,dashed] (0.5,-0.75) ellipse (0.0625 and 0.25) ;
\draw [blue,dashed] (0.5,1) -- (1.1,1) ;
\draw [blue,dashed] (0.5,0.5) -- (1.1,0.5) ;
\draw [blue,dashed] (1.1,0.75) ellipse (0.0625 and 0.25) ;
\draw [red,dashed] (0.5,-1) -- (1.1,-1) ;
\draw [red,dashed] (0.5,-0.5) -- (1.1,-0.5) ;
\draw [red,dashed] (1.1,-0.75) ellipse (0.0625 and 0.25) ;
\begin{scope}
\clip (1,-1) -- (1.5,-1.5) -- (1.5,1.5) -- (1,1.5) -- cycle ;
\draw [blue] (1.5,0.75) ellipse (0.0625 and 0.25) ;
\draw [red] (1.5,-0.75) ellipse (0.0625 and 0.25) ;
\end{scope}
\begin{scope}
\clip (2,-1) -- (1.5,-1.5) -- (1.5,1.5) -- (2,1.5) -- cycle ;
\draw [blue,dashed] (1.5,0.75) ellipse (0.0625 and 0.25) ;
\draw [red,dashed] (1.5,-0.75) ellipse (0.0625 and 0.25) ;
\end{scope}
\draw [blue] (1.5,1) -- (2,1) ;
\draw [blue] (1.5,0.5) -- (2,0.5) ;
\draw [blue] (2,0.75) ellipse (0.0625 and 0.25) ;
\draw [red] (1.5,-1) -- (2,-1) ;
\draw [red] (1.5,-0.5) -- (2,-0.5) ;
\draw [red] (2,-0.75) ellipse (0.0625 and 0.25) ;
\draw (-0.8,1.5) ellipse (0.5 and 0.125) ;
\draw (-1.3,-1.5) -- (-1.3,1.5) ;
\draw (-0.3,-1.5) -- (-0.3,1.5) ;
\draw (0.8,1.5) ellipse (0.5 and 0.125) ;
\draw (1.3,-1.5) -- (1.3,1.5) ;
\draw (0.3,-1.5) -- (0.3,1.5) ;
\begin{scope}
\clip (-1.3,-1.5) -- (1.3,-1.5) -- (1.3,-2) -- (-1.3,-2) -- cycle ;
\draw (-0.8,-1.5) ellipse (0.5 and 0.125) ;
\draw (0.8,-1.5) ellipse (0.5 and 0.125) ;
\end{scope}
\begin{scope}
\clip (-1.3,-1.5) -- (1.3,-1.5) -- (1.3,-1) -- (-1.3,-1) -- cycle ;
\draw [dashed] (-0.8,-1.5) ellipse (0.5 and 0.125) ;
\draw [dashed] (0.8,-1.5) ellipse (0.5 and 0.125) ;
\end{scope}
\end{scope}

\begin{scope}[xshift=8.5cm,yshift=-11cm]
\begin{scope}
\clip (-2.5,-1) -- (-2,-1.5) -- (-2,1.5) -- (-2.5,1.5) -- cycle ;
\draw [blue] (-2,0.75) ellipse (0.0625 and 0.25) ;
\draw [red] (-2,-0.75) ellipse (0.0625 and 0.25) ;
\end{scope}
\begin{scope}
\clip (-1.5,-1) -- (-2,-1.5) -- (-2,1.5) -- (-1.5,1.5) -- cycle ;
\draw [blue,dashed] (-2,0.75) ellipse (0.0625 and 0.25) ;
\draw [red,dashed] (-2,-0.75) ellipse (0.0625 and 0.25) ;
\end{scope}
\draw [blue] (-2,1) -- (-1.3,1) ;
\draw [blue] (-2,0.5) -- (-1.3,0.5) ;
\draw [red] (-2,-1) -- (-1.3,-1) ;
\draw [red] (-2,-0.5) -- (-1.3,-0.5) ;
\draw [blue] (-0.3,1) -- (0.3,1) ;
\draw [blue] (-0.3,0.5) -- (0.3,0.5) ;
\draw [red] (-0.3,-1) -- (0.3,-1) ;
\draw [red] (-0.3,-0.5) -- (0.3,-0.5) ;
\draw [blue] (1.3,1) -- (2,1) ;
\draw [blue] (1.3,0.5) -- (2,0.5) ;
\draw [blue] (2,0.75) ellipse (0.0625 and 0.25) ;
\draw [red] (1.3,-1) -- (2,-1) ;
\draw [red] (1.3,-0.5) -- (2,-0.5) ;
\draw [red] (2,-0.75) ellipse (0.0625 and 0.25) ;
\draw (-0.8,1.5) ellipse (0.5 and 0.125) ;
\draw (-1.3,-1.5) -- (-1.3,1.5) ;
\draw (-0.3,-1.5) -- (-0.3,1.5) ;
\draw (0.8,1.5) ellipse (0.5 and 0.125) ;
\draw (1.3,-1.5) -- (1.3,1.5) ;
\draw (0.3,-1.5) -- (0.3,1.5) ;
\begin{scope}
\clip (-1.3,-1.5) -- (1.3,-1.5) -- (1.3,-2) -- (-1.3,-2) -- cycle ;
\draw (-0.8,-1.5) ellipse (0.5 and 0.125) ;
\draw (0.8,-1.5) ellipse (0.5 and 0.125) ;
\end{scope}
\begin{scope}
\clip (-1.3,-1.5) -- (1.3,-1.5) -- (1.3,-1) -- (-1.3,-1) -- cycle ;
\draw [dashed] (-0.8,-1.5) ellipse (0.5 and 0.125) ;
\draw [dashed] (0.8,-1.5) ellipse (0.5 and 0.125) ;
\end{scope}
\end{scope}

\begin{scope}[xshift=12.3cm,yshift=-11cm]
\draw [thick,->] (0.5,0) -- (-0.5,0) ;
\draw (0,0.5) node{Tube} ;
\end{scope}

\begin{scope}[xshift=15cm,yshift=-11cm]
\draw [thick,blue] (-1,0.3) -- (1,0.3) ;
\draw [thick,red] (-1,-0.3) -- (1,-0.3) ;
\draw [thick,>=stealth,->] (-0.3,-1) -- (-0.3,1) ;
\draw [thick,>=stealth,->] (0.3,1) -- (0.3,-1) ;
\draw (-0.3,0.3) circle (0.15) ;
\draw (0.3,0.3) circle (0.15) ;
\draw (-0.3,-0.3) circle (0.15) ;
\draw (0.3,-0.3) circle (0.15) ;
\end{scope}
\end{tikzpicture}
\caption{Performing Tube($BV$) using Spun$^{\ast}$($BP$)}
\end{figure}

\newpage
\bibliographystyle{plain}
\bibliography{biblio}

\end{document}